\documentclass{amsart}

  \ifx\LabelFigloaded\MYundefined\relax
  \else
   \fi

  \def\LabelFigloaded{\relax}

  \chardef\LabelFigCatAt\the\catcode`\@
  \catcode`\@=11

 \let\LabelFigwlog@ld\wlog
 \def\wlog#1{\relax}

 \ifx\\\MYundefined@
    \let\\\relax
 \fi

  \def\ms@g{\immediate\write16}

 \def\N@wif{\csname newif\endcsname }
 \def\Temp@ {\N@wif\ifIN@}
 \ifx\INN@\MYundefined@
    \else \let\Temp@\relax
 \fi
 \Temp@

  \def\IN@{\expandafter\INN@\expandafter}
  \long\def\INN@0#1@#2@{\long\def\NI@##1#1##2##3\ENDNI@
    {\ifx\m@rker##2\IN@false\else\IN@true\fi}%
     \expandafter\NI@#2@@#1\m@rker\ENDNI@}
  \def\m@rker{\m@@rker}
 
  \newtoks\Initialtoks@  \newtoks\Terminaltoks@
  \def\SPLIT@{\expandafter\SPLITT@\expandafter}
  \def\SPLITT@0#1@#2@{\def\TTILPS@##1#1##2@{%
     \Initialtoks@{##1}\Terminaltoks@{##2}}\expandafter\TTILPS@#2@}

 \def\Shifted@@#1#2#3{\setbox0=\hbox{#3}%
   \raise -\dp0\vbox {\kern-#2%
       \hbox {\kern#1\unhbox0\kern-#1}%
           \kern#2}}

 \newcount\gridcount
 \newbox\auxGridbox@ \newbox\hGridbox@ \newbox\vGridbox@
 \newbox\Labelbox@ \newbox\auxLabelbox@
 \newbox\Coordinatebox@
 \newtoks\Labeltoks@
 \newdimen\Wdd@ \newdimen\Htt@
 \newdimen\Wddd@ \newdimen\Httt@
 
 \def\Wr@{\immediate\write16}

 \newdimen\GL@wd
 \def\GridLineWidth#1{\GL@wd=#1}

 \def\gobble#1{}
 \def\EdgeErr@{\Wr@{}%
      \Wr@{\string\Edges\space argument
      1, 10, 100 or 1000 please\string!}%
      }

 \newcount\Edgect@

 \def\Sweepup#1\endSweepup{}

 \def\SetEdges@{%
    \edef\Zr@@s{\expandafter\gobble\number\Edgect@\empty}%
        \count255=0\Zr@@s\relax
        \ifnum\count255=\z@\else\EdgeErr@\show\tailtest\fi
        \count255=1\Zr@@s\relax
        \ifnum\count255=\Edgect@\relax\else\EdgeErr@\show\leadtest\fi
    \EdgGl@b\edef\Zr@s{\expandafter\gobble\Zr@@s\empty}
    \ifnum\Edgect@>\@ne\relax\EdgGl@b\let\L@Dc\empty
        \else\EdgGl@b\edef\L@Dc{\string.}\fi
    \ifnum\Edgect@>\@ne\relax
        \EdgGl@b\edef\Edgescale@##1{\divide##1 by \Edgect@}%
        \else\EdgGl@b\edef\Edgescale@##1{}\fi
    }

 \def\Edges#1{\Edgect@=#1\relax
     \let\EdgGl@b\global \SetEdges@}

 \Edges{1}

 \def\hhrule{\hrule height \GL@wd\vskip-.\GL@wd}

 \def\hRule@{%
   \advance\gridcount -2%
   \vfil\hhrule\vfil
   \llap{\smash{\raise -2.5pt
     \hbox{\L@Dc\number\gridcount\Zr@s\kern2pt}}}%
   \hhrule
   }

\def\vvrule{\vrule width \GL@wd \kern-\GL@wd}

 \def\vRule@{\advance\gridcount 2%
   \hfil\vvrule\hfil
   \setbox\auxGridbox@=\vbox to 0pt
      {\vskip \Htt@\vskip 2pt
        \hbox to 0pt{\hss\L@Dc\number\gridcount\Zr@s\hss}\vss}%
      \wd\auxGridbox@=0pt \box\auxGridbox@
   \vvrule
   }

 \def\PlaceGrid@@{\gridcount=10 
  \setbox\hGridbox@=\hbox{%
        \hbox{%
             \hskip-.4pt\vrule
             \vbox to \Htt@{%
               \offinterlineskip\parindent=\z@\relax
               \hbox to \Wdd@{\hfil}
               \hRule@\hRule@\hRule@\hRule@
               \vfil\hhrule\vfil}%
             \vrule\hskip-.4pt}
    }%
  \gridcount=0%
  \setbox\vGridbox@=\hbox{%
      \vbox{\offinterlineskip\parindent=0pt\hsize=0pt
         \vskip-.4pt\hrule%
         \hbox to \Wdd@{%
                 \vtop to \Htt@{\vfil}%
                 \vRule@\vRule@\vRule@\vRule@
                 \hfil\vvrule\hfil}%
         \hrule\vskip-.4pt}}%
  \wd\hGridbox@=0pt\ht\hGridbox@=0pt
  \wd\vGridbox@=0pt\ht\vGridbox@=0pt
  \hbox{\box\hGridbox@\box\vGridbox@}%
  }

 \def\LabelsGlobal{\def\LabGl@b{\global}}
 \def\LabelsLocal{\def\LabGl@b{}}
 \LabelsGlobal 

 \def\SetLabels#1\endSetLabels{%
   \LabGl@b\Labeltoks@={#1()\\}%
   }

 \LabGl@b\Labeltoks@={()\\}

 \def\ShowGrid{\LabGl@b\let\PlaceGrid@\PlaceGrid@@}
 \def\HideGrid{\LabGl@b\let\PlaceGrid@\relax}
 \def\Grids{\ShowGrid\LabGl@b\let\GridSwitch@\ShowGrid}
 \def\noGrids{\HideGrid\LabGl@b\let\GridSwitch@\HideGrid}

 \noGrids

 \def\bAdjust@@{%
     \setbox\auxLabelbox@=\hbox{\raise \dp\auxLabelbox@
            \box\auxLabelbox@}}
 \def\bAdjust@{\let\vAdjust@\bAdjust@@}

 \def\eAdjust@@{\dimen0=-.5\ht\auxLabelbox@
     \advance\dimen0 by .5\dp\auxLabelbox@
     \setbox\auxLabelbox@=
            \hbox{\raise\dimen0\box\auxLabelbox@}}
 \def\eAdjust@{\let\vAdjust@\eAdjust@@}

 \def\tAdjust@@{%
     \setbox\auxLabelbox@=\hbox{\raise-\ht\auxLabelbox@
            \box\auxLabelbox@}}
 \def\tAdjust@{\let\vAdjust@\tAdjust@@}

 \let\vAdjust@\relax

 \def\lAdjust@{\let\hAdjust@\rlap}
 \def\rAdjust@{\let\hAdjust@\llap}

 \let\hAdjust@\relax\let\vAdjust@\relax

 \def\FetchLabel@#1(#2)#3\\{%
     \IN@0#2@@\ifIN@
        \setbox0=\hbox{\ignorespaces#1#3\unskip}%
        \ifdim\wd0>0pt
           \ms@g{}%
           \ms@g{ !!! Bad label(s)? !!!}%
           \message{ #1(#2)#3}%
        \fi
        \def\LabelMole@##1\endFetchLabel@{%
            \IN@0()\\@##1@%
            \ifIN@\def\Temp@{\FetchLabel@##1\endFetchLabel@}%
            \else\def\Temp@{}%
            \fi
            \Temp@
           }%
     \else
       \ignorespaces#1\unskip
       \setbox\auxLabelbox@=%
         \hbox to 0pt{\hss\ignorespaces\hAdjust@
          {\ignorespaces#3\unskip}\hss}%
       \vAdjust@
       \let\hAdjust@\relax\let\vAdjust@\relax
       \AugmentLabelBox@@{#2}%
       \ht\Labelbox@=0pt\dp\Labelbox@=0pt
       \let\LabelMole@\FetchLabel@%
     \fi\LabelMole@}

 \newtoks\XYSep@ 
 \def\SetXYSeparator#1{%
     \IN@0#1@@\ifIN@\XYSep@{*}%
     \else
     \XYSep@{#1}%
     \fi
     }

 \SetXYSeparator*

 \def\AugmentLabelBox@@#1{%
     \IN@0\the\XYSep@ @#1@\ifIN@
       \SPLIT@0\the\XYSep@ @#1@%
       \setbox\Labelbox@=\hbox to 0pt{%
         \unhbox\Labelbox@
         \Shifted@@{\the\Initialtoks@\Wddd@}%
         {\the\Terminaltoks@\Httt@}%
         {\box\auxLabelbox@}}%
     \else
         \ms@g{}%
         \ms@g{ !!! Bad insertion point. !!!}%
         \message{ (#1\ this point was rejected.)}%
     \fi
    }

 \def\FetchOption@#1[#2]#3\endFetchOption@{%
    \def\temp{#1}
    \ifx\temp\empty
       \Edgect@=#2\relax
       \let\EdgGl@b\relax
       \SetEdges@
       \Cleaner@#3%
    \fi}

 \def\Cleaner@#1[@]{\Labeltoks@{#1}}
     
 \def\PlaceLabels@@{\mathsurround=0pt
     \def\Cr@{\\}%
     \let\L\lAdjust@\let\R\rAdjust@
     \let\B\bAdjust@\let\E\eAdjust@\let\T\tAdjust@
     \expandafter\FetchOption@\the\Labeltoks@[@]\endFetchOption@
     \Wddd@=\Wdd@ \Edgescale@\Wddd@ 
     \Httt@=\Htt@ \Edgescale@\Httt@
     \expandafter\FetchLabel@\the\Labeltoks@\endFetchLabel@
     \box\Labelbox@
     }%

 \let \PlaceLabels@\PlaceLabels@@

 \def\AffixLabels#1{\setbox\Coordinatebox@=\hbox{#1}%
      \Wdd@=\wd\Coordinatebox@ \Htt@=\ht\Coordinatebox@
      \advance\Htt@ \dp\Coordinatebox@
      \hbox{\copy\Coordinatebox@\kern-\Wdd@ 
           \Shifted@@{0pt}{-\dp\Coordinatebox@}%
           {\PlaceLabels@\PlaceGrid@}%
           \kern\Wdd@}%
      \GridSwitch@ 
      \LabGl@b\Labeltoks@{()\\}%
      }
 
   \let\wlog\LabelFigwlog@ld   
   \catcode`\@=\LabelFigCatAt  

                                By

              Raymond S\'eroul <A18645@FRCCSC21.BITNET>
                                and 
              Laurent Siebenmann <lcs@topo.math.u-psud.fr>
    
              VERSIONS: July 1991, Oct 1991, Jan 1992, July 1992

INTRODUCTION

      This labelling package is intended for TeX users who
rely on non-TeX sources for for their graphics inserts.  It
provides means for adding TeX labels to such inserts with a
minimum of fuss. 

       For most labels, TeX users have in the past found it
reasonably convenient to rely on non-TeX sources. Typical
occasions when an inescapable need for TeX labels seemed to
arise are

 (a) when the graphics program lacks certain exotic or complex
mathematical symbols

 (b) when the very highest typographical quality is wanted for the
labels

 (c) when labels included with the graphics fail to print, 
 and you cannot figure out why (cf. boxedeps.doc).  The labels
 provided by labelfig.tex are 100

       Since this package first appeared, many users, who in the
past scarcely dreamed of using TeX labels, have come to use
nothing but.  So it is now appropriate to add

Intoxication Warning:  TeX labels may be addictive and expensive. 

     If you have a fast preview you may disagree, and even find
that this package provides an agreeable paste-up environment; see
extra applications at end.

     Note to publishers: It is possible and convenient to ultimately
export the TeX labels produced by labelfig.tex to become an integral
part of the EPS file. This is often desired by a publisher who typically
uses an "upmarket" graphics or page layout program, with which the
staff is skilled in perfecting figures.  See Appendix I for
a recipe.

     The authors are grateful to Patrick Ion of Math Reviews for
helpful comments and encouragement.

BASIC INSTRUCTIONS

    After reading in the macro file using

preview or proof your figure with a coordinate grid printed on
top, by typing the following:

    \ShowGrid  
    \AffixLabels{<the graphics insertion>}

Here <the graphics insertion> is what you would type to insert
the graphics object alone without the grid.  This must provide
for the space around it. For example <the graphics insertion>
might well be \BoxedEPSF{MyFigure scaled 700} using the
boxedeps.tex macro package (from same source); this provides a
TeX box containing the encapsulated PostScript insert specified by
the file MyFigure. \AffixLabels{...} provides the grid (supposing
\ShowGrid is present) and later, once you have specified labels
using the grid, it will "tack on" the labels.

     The grid is a sort of (usually elongated) checkerboard of
ten rows and ten columns and its (internal) partitions are by
default numbered  .1, ... ,.9  both horizontally (X-coordinate
running left to right) and vertically (Y-coordinate running bottom
to top).  Thus the points enclosed by the grid correspond to the
points of the unit square in the cartesian "X-Y" plane, the lower
left corner corresponding to the origin (0,0).  By extrapolation,
the full page corresponds to a larger rectangle in the plane.

     These coordinates serve to position labels as follows.
Before the \AffixLabels{...} command type label specifications:

  \SetLabels
   (<X-coordinate>*<Y-coordinate>) <first label> \\
   .
   .
   .
   (<X-coordinate>*<Y-coordinate>)  <last label> \\
  \endSetLabels

Each row specifies one label and is terminated by \\.  In each
row, the position indicator comes first; it is written as a
standard cartesian point except that the X- and Y- coordinates
are separated by * rather than a comma because TeX allows a
comma as decimal point. There are no dimension units to specify
as the unit is the grid itself.

     By default, this cartesian point specifies where the middle
of the baseline of the label will be located.  However if you precede
the point by \L [or \R] the left [or right] edge of the baseline will
be located there. Similarly you may also precede the point by \T, \E,
or \B to vertically align the top equator or bottom of the label box
at the specified point.  This gives nine standard positions of
the label with respect to the insertion point --- corresponding to
the eight principle points of the compas and the center

                     \L\T     \T      \R\T

                     \L\E     \E      \R\E

                     \L\B     \B      \R\B

But this neglects the default "baseline" level of TeX,
giving potentially three more positions

                     \L    <no tag>   \R

For text, the baseline level is often the preferred. Its relation to
the others is variable. It will often coincide with the bottom level,
as happens for "X".  But it is often distinct, as for "g", in which
case you have in all 12 distinct positions rather than 9.

     It is convenient to think of this specification of label
position as attaching the label by a thumb-tack to the coordinate
grid. There are up to twelve positions of the thumb-tack on the
label, while the position of the thumb-tack on the coordinate grid is
arbitrary.  Normally, one choses the position of the thumb-tack on
the label to be the one that is the closest to the item being
labeled.  There are good reasons for this "rule of thumb":

   (a)  It facilitates correct positioning at first try.

   (b)  If the scale of the figure must be altered after labels
have been affixed, the labels have a good chance of remaining well
positioned.

   (c)  The visible grid need not extend beyond the "bounding box"
for the figure, because the best preferred position is always
(at least almost) within the bounding box .

The second reason is particularly important. Indeed it often
happens that scale has to be altered after labelling begins, in
order to either provide space for the labels, or to adjust
proportions between the labels and the figure.  (The size of labels
is unaffected by scaling.)

     Here is an artificial but self-contained test which uses
TeX rules to make a graphics object.

TEST

    Do not skip this!

 \def\FrameIt#1{\hbox{\vrule$\vcenter {\hrule\kern3pt%
             \hbox {\kern3pt #1\kern3pt}%
               \kern3pt\hrule}$\relax\vrule}}

 \def\Caption#1#2{\FrameIt{%
       \vtop {\hsize=#1\relax \parindent=0pt
         \leftskip=0pt \rightskip=0pt plus15pt
         \parfillskip=0pt
         \lineskip=1pt\baselineskip=0pt
         #2}}}

 \def\FirstQuadrant{\hbox to 100pt{\vrule\vbox to 100pt{%
        \hbox to 100pt{\hfil}\vfil\hrule}\hss}}

  \SetLabels
    \R(.5*.2) $\zeta\,\cdot$\\
    (.9*-.10) $\xi$\\
    \R(-.03*.9) $\eta$\\
    \T(.5*.9) \Caption{70pt}{%
          \it The norm of
          $g(\xi+i\eta)$ is indicated on
          contours of this invisible surface.}\\
  \endSetLabels

  \AffixLabels{\FirstQuadrant}

  \end

  Note that the coordinates to use for labels are indicated on the
edges of the grid (when visible) corresponding to the conventional
x- and y- axes of the Cartesian plane. By default the grid is
1-by-1. However, by the command \Edges{100}, you can change this
to 100-by-100 and many users find this alternative most
convenient. Place the command \Edges{...} in your style file (or
header) since its effect is is global. Other possible edge values
are 10 and 1000.

  If you use the command \Edges{...} at all, do so with care.  For
if you accidentally delete an \Edges{...} command your labels will
abruptly be badly misplaced and may logically but mysteriously
generate "dimension too big" errors under TeX and "off page" errors
under your driver.  

  You can dictate the edgescale for an individual figure by giving
the scale in brackets immediately after \SetLabels.  Thus, to
import into an article using say \Edge{100} a figure labelled using
another edgescale, say the original 1-by-1 default, you can use
\SetLabels[1]...\endSetLabels.

GETTING IT DOWN PAT

     Complicated labeling deserves the same respect as
complicated mathematics.  Do not expect it to come out perfect the
first time!  What is needed in either case is a mechanism to
repeatedly typeset troublesome pieces.

     One mechanism is always available.  One does complicated
labelling in a separate "test" file involving just the figure being
labelled;  a texpert will know how to \dump TeX's current state as
a temporary format that restarts rapidly at each retry.  Usually,
one then pastes the completed labelled figure back into the main
TeX file, but, of course, one can also \input it as an auxiliary
file.

     If you do not have a TeXpert at handy, here is a first
approximation to an efficient setup. By deletions reduce a copy
of your article to just a few lines before and after the figure.
Now label the figure, and finally, copy and paste the labelled
figure to the original article. Then copy the next figure to label
into this testbed and repeat. The TeXpert can improve the  speed
at which TeX starts up, by compiling a format specifically for
your article; just one caution: best NOT include in the format
ephemeral details of setup like \Set<mydriver>ArtSpecials (from
boxedeps.tex because this reads  figure dimensions which you may
change during your work session.

     An improved mechanism to repeatedly typeset troublesome
pieces is now available on the Macintosh; it is called LinoTeX;
see the same ftp sources.  It could be set up on many types
of computer.

     Before using labelfig.tex to attach labels to a graphics
object inserted using boxedeps.tex or BoxedArt.tex, make it a
firm rule to carefully adjust the bounding box using the trimming
commands of these packages, and also at least tentatively scale
and position the object. Beware of changing the grid inadvertently
after the labels have been positioned.  For example, correcting
the bounding box of a PostScript graphics object can foul up the
labels by changing the coordinate grid to which the labels are
attached. This is particularly true for the trimming  commands of
boxedeps.tex and BoxedArt.tex. However, as noted already, change
of scale is much less disruptive, and modest adjustments should be
well tolerated.

     Sometimes the labels protrude so far from the bounding box
of a figure that the figure has to be repositioned.  Best do this
by ad hoc spacing, say using \hglue and \vglue; altering the
bounding box would create a vicious circle.

     Remember that you are responsible for preventing labels
from overlapping. You are responsible for all label typography
including size and style. A label is really just about anything
that can be put in a TeX box. Note that spaces at the beginning
and end of labels will normally be suppressed; if you really want
them you must protect them with TeX braces.

     This package temporarily sets the \mathsurround parameter
of TeX to zero  while the labels are being affixed. This is done
because nonzero \mathsurround space would influence the position
of left and right aligned labels; then, when a texpert or printer
modifies mathsurround, diagram labeling might be disastrously
altered. There is a small price to pay involving labels that are
formatted as caption boxes including mathematics: you  may want or
need to specify an explicit mathsurround space within the caption
box; it will not influence anything outside.

     Those hostile to the use of * as separator between
the X and Y coordinates of label insertion points, are free to
impose another using \SetXYSeparator{<the new separator>}.  
Americans may prefer "," to "*" since they never use a 
comma as a decimal point; on the other hand, * may be more visible.

APPENDIX (I)  MERGING labelfig.tex LABELS INTO AN EPSF GRAPHICS OBJECT.

     As promised in the introduction, here is a recipe useful for
publishers. It works at least on Macintosh and at least for vectorized
graphics and Adobe type1 fonts.  (There is surely a similar recipe for
PCs under MSWindows.)

 (a)  Use boxedeps.tex utility to integrate the figure given by the eps
file, "x.eps" say, with a visible frame around it.  See
\ShowDisplacementBoxes command in boxedeps.tex.  To get precise results
automatically it is important to use the \Trim... commands of
boxedeps.tex making the "DisplacementBox" neatly fit the figure.

 (b)  Use the TeX printer driver and LaserWriter (versions >= 8.1.1) to
export to an EPSF the DVI page containing the integrated, labelled
figure. You now have an EPS file  "xx.eps"  that contains too much, and at
the wrong scale, and at wrong position.

 (c)  Convert the EPSF to an Adode Illustrator format EPSF using
the shareware utility called epsConvert by Sam Weiss
1993-- (currently $25).

 (d)  In Illustrator (or a compatible program), group the labels and the
"DisplacementBox"; copy them to the clipboard and paste them into "x.ps".
This step requires that all the label fonts be "visible to the Macintosh.

 (e)  Translate and scale the pasted group consisting of the labels plus
the "DisplacementBox" so as to make the "DisplacementBox" the bounding
box of (labelless) figure represented by "x.eps".  At this point the
labels will be correctly placed on the figure "x.eps".

 (f)  Ungroup and delete the "DisplacementBox".  The result is the
desired single EPS file, "x+.eps" say, It contains the original figure
plus its labels.  

     Using grouping and ungrouping appropriately in "x+.eps", a
publisher's staff can very efficiently improve label positions etc.

APPENDIX II)  SOME EXOTIC APPLICATIONS

     The grid of labelfig.tex is analogous to a light-table in
classical page makeup with wax or latex glue.  In principle, you
can use it to compose any page from its indivisible parts.  This
even has some of the artisanal charm of classical paste-up
provided you have a fast screen preview to make the process
"interactive".

     In practice labelfig.tex is a tool for nonstandard jobs.
Here are a few going beyond the labelling already discussed.

(I)  GRAPHICS INTEGRATION.

     This is accomplished by treating the imported graphics
objects as labels.  The underlying graphics object is then
typically an empty  \vbox to <dimension>{\vfill} in a TeX
\midinsert...\endinsert construction.  A label line
might be of the form

   (.1*.1) \\

The exact form of the special command varies from driver to
driver.  However, in the case of encapsulated PostScript graphics
(EPSF norm), by relying on boxedeps.tex, one can have the
following standard syntax (independant of driver  (see
boxedeps.doc for details.
  
  (.1*.1) \BoxedEPSF{MyFigure scaled <scale in mils>}\\

This may be slow since it requires TeX to read the PostScript
file to read bounding box using many complex macros.  So you
may want to try

  (.1*.1) \EPSFSpecial{MyFigure}{<scale in mils>}\\

which is fast and driver independant, but it squashes the
bounding box, normally to its lower left corner.

     Similarly for graphics of the Macintosh PICT norm ---
using BoxedArt.tex (same sources) in place of boxedeps.tex.

     This approach to integration is to be recommended when
one is assembling a composite graphics object.

 (II)  COMMUTATIVE DIAGRAM ENHANCEMENT

     Commutative diagrams or arrays of mathematical objects
connected by arrows of various sorts are common in mathematics.
The mathematical objects require the use of TeX.  Recently TeX
acquired a good collection of arrows of all slopes --- that of
LamSTeX --- plus pwerful macros to build the diagrams.

     However, even the LamSTeX collection is often
inadequate; it lacks for example double shafted arrows, dotted
arrows and curved arrows. Fortunately it is possible to produce
such arrows on an individual basis using sophisticated graphics
programs such as Illustrator and AldusFreehand (both serving
the EPSF norm) or using Metafont (with its public domain norm).
Since the creation of each new arrow is a work of love, you
probably want to limit the number of arrows by using LamSTeX
for most arrows. The 40K commutative diagram module of LamSTeX
has been adapted to work with AmSTeX and a copy may be posted
with LabelFig and related files. Unfortunately no one has yet
offered a version that works with Plain TeX or LaTeX.

       Suffice it here to say that when the exotic arrow has
been somehow imported into TeX, labelfig.tex treats it as a
label that one affixes to the commutative diagram.  Two other
steps will be treated in separate notes, namely the matter of
extracting the dimension specifications for the arrow and the
construction of the arrow --- for these steps are far from
unique and often depend intimately on your computer environment. 
Notes for the Macintosh-Textures-Illustrator combination are
found in the file ExoticArrows.doc.

 (III) NESTING 

Ingenuity pays off in exploiting labelfig.tex. One can
mix graphics and typography quite freely.  labelfig.tex is good
for freeform or overlapping arrangements, while boxedeps.tex (or
BoxedArt.tex) is best for regimented non-overlapping
arrangements --- and the two can be combined.

     The default behavior of labelfig.tex is not ideal 
for nesting objects, because to prevent trouble for beginners
the register for labels is globally cleared when \AffixLabels
concludes.  But there are switches available

      \LabelsGlobal      \LabelsLocal

which change this.  To understand this, extend the above test 
by something like:

 \LabelsLocal

 \SetLabels
    (.5*.5) AAA\\
 \endSetLabels

 {
 \SetLabels
    (.5*.5) ZZZ\\
 \endSetLabels
   \AffixLabels{\FirstQuadrant}
 }

   \AffixLabels{\FirstQuadrant}

     There are however potential pitfalls.  Neither
labelfig.tex nor boxedeps.tex has been tested under extreme
conditions. Problems may occur if their procedures are
indiscriminately nested. For boxedeps.tex (not labelfig.tex)
there is a precise cause for worry, namely many of its
variables are "global", which means that TeX braces will not
provide the protection one might expect.

COMMAND SUMMARY FOR labelfig.tex

  Here [...] means optional (one or zero)
       [...]* means any number of such constructs

  \SetLabels
    [[<P>](<X><Sep><Y>) <label> \\]*
  \endSetLabels
  \ShowGrid  
  \AffixLabels{<the figure>}

   --- <P> is tack position, one of eleven or empty
              order irrelevant

                   \L\T      \T      \R\T

                   \L\E      \E      \R\E

                     \L               \R

                   \L\B      \B      \R\B

   --- (<X><Sep><Y>) insertion point;
  <Sep> is separator, = * by default;
  \SetXYSeparator{<Sep>} changes it.
   <X> and <Y> are real numbers

  --- <label> a label to attach 

  --- <the figure> the figure to label 

  \GlobalLabels (default)     
  \LocalLabels  setting for nested constructs.

 \Grids makes ALL grids appear; \HideGrid then makes just next disappear.
 \noGrids returns to default.  The commands are always global.

 \GridLineWidth{<dimension>} adjusts width of grid lines. Default is very
small, to give "hairline" effect. If your grid lines are missing try
setting \GridLineWidth{1pt}.

 \Edges#1 globally changes the edge size of all grids to the numerical 
value #1, which must be 1, 10, 100, or 1000.  The default is 1.

VERSION HISTORY.
 --- Jan 1993: \Edges#1 and [??] option after \SetLabels
 --- July 1992: \Grids, \noGrids, \HideGrid;
       Gridlines become hairlines; \GridLineWidth{<dimension>}.
 --- Oct 1991, Jan 1992: \SetXYSeparator{<Sep>},  \LabelsGlobal,
       \LabelsLocal.
 --- July 1991: first release

Address for bugs and other feedback:

        Raymond S\'eroul
        IREM and Lab. de Typographie Informatise
        Univ. Rene Descartes
        Strasbourg

    Tel 33-88-41-63-45
    Email:  A18645@FRCCSC21.BITNET

        Laurent Siebenmann
        Mathematique, Bat. 425,
        Univ de Paris-Sud,
        91405-Orsay,
        France

    Tel 33-1-6941-7949; 
    Email: lcs@topo.math.u-psud.fr

\input xy
\xyoption {all}

\usepackage{graphicx}
\usepackage{latexsym}
\usepackage{amssymb} 
\usepackage{color} 
\usepackage[frame,cmtip,arrow,matrix,line,graph]{xy}

\usepackage{amsmath} 
\usepackage{amsthm} 
\usepackage{amsfonts} 
\usepackage{amscd}

\usepackage{enumerate}
 
\usepackage{epsfig} 

\usepackage{caption2}
\usepackage{mathrsfs}

\newtheorem{thm}{Theorem}[section] 
\newtheorem{lem}[thm]{Lemma} 
\newtheorem{prop}[thm]{Proposition} 
\newtheorem{rem}[thm]{Remark} 
\newtheorem{cor}[thm]{Corollary} 
\newtheorem{Def}[thm]{Definition} 
 
\newtheorem{ex}[thm]{Example}
\newtheorem{conj}[thm]{Conjecture}

\newcommand{\Dif}{{\rm{Diff}} } 
\newcommand{\Conf}{{\rm{Conf}} } 
\newcommand{\Aut}{\operatorname{Aut}}
\newcommand{\Out}{\operatorname{Out}}

\newcommand{\link}{\operatorname{Link}}
\newcommand{\St}{\operatorname{St}}
\newcommand{\SAut}{\Sigma{\rm{Aut}}}
\newcommand{\SHAut}{\Sigma_H\mkern-1mu{\rm{Aut}}}

\newcommand{\col}{\colon\!}

\newcommand{\Ñ}{\kern-.1pt\vrule height2.8ptwidth3ptdepth-2.3pt\kern1pt}
\newcommand{\D}{\mathcal{D}}
\newcommand{\Sph}{\mathcal{S}}

\newcommand{\Z}{\mathbb{Z}} 
\newcommand{\R}{\mathbb{R}} 
\newcommand{\T}{\mathbb{T}} 

\newcommand{\G}{\Gamma}
\newcommand{\del}{\partial} 

\newcommand{\cdotss}{\mathinner{\mkern.5mu\cdotp\mkern-2mu\cdotp\mkern-2mu\cdotp}}
\newcommand{\tild}{\widetilde}

\newcommand{\cln}{\colon\!}

\newcommand{\emp}{\varnothing}
\newcommand{\nat}{\,\natural\,}

\newcommand{\rar}{\rightarrow} 
\newcommand{\inc}{\hookrightarrow}

\newcommand{\lgl}{\langle}
\newcommand{\rgl}{\rangle}

\newcommand{\s}{\sigma}

\newcommand{\Si}{\Sigma}
\newcommand{\Ga}{\Gamma}

\newcommand{\Pts}{\Lambda}

\newcommand{\x}{\!\times\!}

\newcounter{samcounter}

\begin{document}

\title{Stabilization for mapping class groups of 3-manifolds}

\author{Allen Hatcher}
\author{Nathalie Wahl}


\begin{abstract}
We prove that the homology of the mapping class group of any 3-manifold stabilizes under connected sum and boundary connected sum with an arbitrary 3-manifold when both manifolds are compact and orientable. The stabilization also holds for the quotient group by twists along spheres and disks, and includes as particular cases homological stability for symmetric automorphisms of free groups, automorphisms of certain free products, and handlebody mapping class groups. Our methods also apply to manifolds of other dimensions in the case of stabilization by punctures. 
\end{abstract}


\maketitle

The main result of this paper is a homological stability theorem for mapping class groups of 3-manifolds, where the stabilization is by
connected sum with an arbitrary $3$-manifold. More precisely, 
we show that given any two compact, connected, oriented 3-manifolds $N$ and $P$ with $\del N\neq \emp$, 
the homology group $H_i(\pi_0\Dif(N\# P\#\cdots\# P\ \textrm{rel}\ \del N); \Z)$ is independent of the number $n$ of
copies of $P$ in the connected sum, as long as $n\ge 2i+2$, i.e. each homology group stabilizes with $P$. 
We also prove an analogous result for boundary connected sum, and a version for the quotient group of the mapping class group
by twists along spheres and disks, a group closely related to the automorphism group of the fundamental group of the manifold.

Homological stability theorems were first found in the sixties for symmetric groups by Nakaoka \cite{Na} and linear
groups by Quillen,  
and now form the foundation of modern algebraic K-theory (see for example
\cite[Part IV]{LS} and \cite{vdK}). 
Stability theorems for mapping class groups of surfaces were obtained in the eighties by Harer and
Ivanov \cite{harer,I} and recently turned out to be a key ingredient to a solution of  
the Mumford conjecture about the homology of the Riemann moduli space \cite{MW}. 
The other main examples of families of groups for which stability has been known are braid groups
\cite{A} and automorphism groups of
free groups \cite{H3,HV1}. 

The present paper extends significantly the class of groups for which homological stability is known to hold. 
It suggests that it is a widespread phenomenon among families of groups containing enough `symmetries'. In addition to the already
mentioned stability theorems  for mapping class groups of 3-manifolds, corollaries of our main result 
include stability for handlebody subgroups of surface mapping class groups,  symmetric
automorphism groups of free groups, and automorphism groups of free products $*_nG$ for
many groups $G$. 
Using similar techniques we obtain stability results also for mapping class groups
$\pi_0\Dif(M\!-\!\{n\ {\rm points}\}\ \textrm{rel}\ \del M)$ for $M$ any $m$-dimensional manifold with boundary, $m\ge 2$ (even the case $m=2$ is new here), as well as for the fundamental group $\pi_1\Conf(M,n)$ of the
configuration space of $n$ unordered points in $M$.
Our paper
thus also unifies previous known results as we recover stability for braid groups (as $\pi_1\Conf(D^2,n)$),
symmetric groups (as $\pi_1\Conf(D^3,n)$) and automorphism groups
of free groups (as $\Aut(*_n\Z)$).

 For mapping class groups of $3$-manifolds, we show moreover that, when $P=S^1\times S^2$, the group 
$H_i(\pi_0\Dif(N\# P\#\cdots\# P\ \textrm{rel}\ \del N);\Z)$ is independent of the number of
 boundary spheres and tori coming from $D^3$ and
$S^1\times D^2$ summands in $N$ when there are at least $2i+4$ copies of $P$.  This is to be compared to the dimension
 2 case, where the same result holds for $P=S^1\times S^1$ with respect to $D^2$ summands in
 $N$, i.e. the corresponding homology group is independent of the number of boundary components fixed by the mapping class group when
 the genus of the surface is large enough \cite{harer,I}.  
  This type of result does not hold for a general $3$-manifold $P$, however, as we show by an example, although a version of it does hold for boundary connected sum with $P=S^1\times D^2$.

\smallskip

Following a standard strategy in geometric proofs of stability, we prove our results by building highly connected simplicial complexes on which the
groups act, and deduce stability from the spectral sequence associated to the action. 
To be able to consider many families of groups at once, we have
axiomatized a large part of the process, in particular the spectral sequence argument. 
The core of the proof of such a stability theorem is showing that the complexes are indeed highly connected. 
For this, we introduce a new combinatorial technique 
which shows that certain complexes
obtained from a Cohen-Macaulay complex by adding labeling data on its vertices are highly connected.

Our paper is concerned with the groups of components of the diffeomorphism groups of $3$-manifolds. 
 The diffeomorphism group of a $3$-manifold does not
in general have contractible components, though contractibility holds for most prime $3$-manifolds. 
For example, the diffeomorphism group $\Dif(\mathcal{H}_g)$ of a handlebody of genus $g\ge 2$ has contractible components, and
 our result gives stability for the homology of the classifying space $B\Dif(\mathcal{H}_g)$ with respect to genus.  
A question left open by the results in this paper
is whether stability holds in general for the full diffeomorphism group of a $3$-manifold, not just the group of components.

\smallskip

The second author would like to thank Andr\'e Henriques, Marc Lackenby and Thomas Schick
for stimulating conversations at the Oberwolfach Topology Meeting 2006, after which the scope of this paper grew tremendously. The authors would furthermore like to thank Karen Vogtmann for many helpful conversations.
The second author was partially supported by the NSF grant DMS-05044932, the Mittag-Leffler Institute and the Danish Natural Science Research Council.

\section{statement of results}\label{results}

Throughout the paper we use results about $2$-spheres in $3$-manifolds from our earlier paper \cite{HW}. These results require that the $3$-manifolds in question, as well as their universal covers, contain no connected summands which are counterexamples to the Poincar\'e conjecture. As this conjecture has now been proved, we will make no further mention of this underlying hypothesis. The more general geometrization conjecture is also used, in a less essential way, in Section~\ref{mcg} (see Proposition~\ref{kernel}).

\smallskip

Let $M$ be a compact, connected, orientable 3-manifold with a boundary component $\del_0M$, and let 
$$\G(M)=\G(M,\del_0M)=\pi_0\Dif(M \ \textrm{rel}\ \del_0M\,)$$ 
denote the mapping class group of $M$, the group of path-components of the space of orientation-preserving diffeomorphisms of $M$ that restrict to the identity on $\del_0M$. Taking the induced automorphisms of $\pi_1M$ using a basepoint in $\del_0M$ gives a canonical homomorphism $$\Phi\col\G(M)\to\Aut(\pi_1M).$$  The following facts about $\Phi$ are explained in more detail in Section~2. The kernel of $\Phi $ contains the subgroup generated by twists along embedded spheres and properly embedded disks. 
This subgroup is normal, so we can form the quotient group $A(M)$, with an induced homomorphism $A(M)\to\Aut(\pi_1M)$.  This last map is injective when $\del M - \del_0M$ has no sphere components, and in many interesting cases it is also surjective.

\smallskip

Suppose that $M$ is the connected sum $N_n^P=N\# P\#\cdotss \# P$ of a manifold $N$ with $n$ copies of a manifold $P$. We assume $\del N\ne\emp$ and we choose a component $\del_0N$ of $\del N$.  Let $R$ be any compact subsurface of $\del N$ that contains $\del_0N$. We denote by 
$$\G_n^P(N,R )=\G(N_n^P,R )=\pi_0\Dif(N_n^P \ \textrm{rel}\ R\,)$$  
the mapping class group of $N_n^P$ fixing $R$, and by 
$$A_n^P(N,R,T)=\G_n^P(N,R)/K(T)$$
the quotient group by the subgroup $K(T)$ of $\G_n^P(N,R)$ generated by twists along spheres and disks with boundary in  a possibly empty compact subsurface $T$ of $\del N^P_n$ disjoint from $R$ and 
invariant under $\G^P_n(N,R)$.

There is an inclusion $N^P_n \to N^P_{n+1}$ obtained by gluing on a copy of the manifold $P^0$ obtained from $P$ by deleting the interior of a ball, where the gluing is done by identifying a disk in the resulting boundary sphere $\del_0 P^0=\del P^0 -\del P$ with a disk in $\del_0N$. By extending diffeomorphisms of $N^P_n$ to diffeomorphisms of $N^P_{n+1}$ via the identity on the adjoined $P^0$ we obtain maps
$$
\G^P_n(N,R)\to\G^P_{n+1}(N,R) \qquad {\rm and} \qquad A^P_n(N,R,T)\to A^P_{n+1}(N,R,T)
$$
where $T$ is extended to $N_{n+1}^P$ by invariance under the $\G_{n+1}^P(N,R)$-action.

Here is our main result:
\begin{thm}\label{main}
{\rm (i)} For any compact, connected, oriented 3-manifolds $N$ and 
$P$ with subsurfaces $R$ and $T$ of $\del N^P_n$ as above, 
the induced stabilization maps
$$
 H_i(\G^P_n(N,R)) \to H_i(\G^P_{n+1}(N,R)) \ \ \textrm{and}\ \ 
H_i(A^P_n(N,R,T)) \to H_i(A^P_{n+1}(N,R,T))
$$
are isomorphisms when $n\ge 2i+2$ and surjections when $n= 2i+1$. \vskip2pt\noindent
{\rm (ii)} When $P=S^1\x S^2$ and $n\ge 2i+4$, the groups $H_i(\G_n^P(N,R))$ and $H_i(A_n^P(N,R,T))$ are moreover independent of the number of $D^3$ and $S^1\times D^2$ summands of $N$ whose boundary spheres and tori are contained in $R$.
 \end{thm}

Here and throughout the paper we use untwisted $\Z$ coefficients for homology. 
The case where $P$ is a prime $3$-manifold is of particular interest in examples. 
Note that the theorem for $P$ prime implies that the same result holds for $P$ not prime. 

In statement (ii) of the theorem the isomorphisms are induced by filling in boundary spheres with balls and filling in boundary tori of $S^1 \times D^2$ summands with solid tori to produce trivial $S^3$ summands. The number of $D^3$ or $S^1\times D^2$ summands with boundary in $R$ can be zero (and in this case $R$ is specifically allowed to be empty), but if this extreme case is avoided then the inequality $n\ge 2i+4$ can be improved to $n\ge 2i+2$.

The special case $P=S^1\x S^2$ and $N=D^3$ recovers homological stability for $\Aut(F_n)$ and $\Out(F_n)$, the automorphism and outer automorphism groups of free groups, as proved in \cite{H3,HV1,HV3,HVW}. The proof given here is somewhat simpler. More generally, if $N$ is 
the connected sum of $s\ge 1$ balls and $k\ge 0$ solid tori (and possibly some extra summands), with $R=\del N$ and $T=\emp$, then $A^P_n(N,R,T)$ is the group denoted $A^s_{n,k}$ in \cite{HW}, and the theorem recovers the main stability results of \cite{HW} with an improvement in the stable dimension range.

Another example can be obtained by taking $P=S^1\x D^2$, $N=D^3$, $R=\del N$, and $T=\del N^P_n - \del N$.  
Then $\pi_1N^P_n$ is the free group $F_n$ and $A^P_n(N,R,T)$ is the symmetric automorphism group $\SAut(F_n)$, the subgroup of
$\Aut(F_n)$ generated by the automorphisms that conjugate one basis element by another, send a basis element to its inverse, or permute
basis elements. 
This group is also known as the `circle-braid' group, the fundamental group of the space of configurations of $n$ disjoint unknotted, unlinked circles in $3$-space, studied in \cite{D,BWC,JMM,Ru,BH}.
Thus we have:

\begin{cor}\label{SymAut}
The stabilization $\SAut(F_n)\to\SAut(F_{n+1})$ induces an isomorphism on $H_i$ for $n\ge 2i+2$ and a surjection for $n=2i+1$.
\end{cor}

More generally, for a free product $*_iG_i$ of a finite collection of groups $G_i$, the symmetric automorphism group $\SAut(*_iG_i)$ consists of the automorphisms that take each $G_i$ onto a conjugate of a $G_j$. When no $G_i$ is $\Z$ or a free product, this subgroup is all of $\Aut(*_iG_i)$ (see for example \cite{MM2}). Consider the case where $*_iG_i=*_nG$ is the free product of $n$ copies of the same group $G$. 
Then for any subgroup $H$ of $\Aut(G)$ containing the inner automorphisms, we can look at the subgroup $\SHAut(*_nG)$ of $\SAut(*_nG)$ consisting of automorphisms taking each factor $G$ to a conjugate of another $G$ factor via the composition of an automorphism in $H$ and conjugation in $*_nG$. The following result is then a corollary of the proof of our main theorem rather than of the theorem itself:

\begin{cor}\label{free}
Let $P$ be a prime compact orientable $3$-manifold with $G=\pi_1(P)$ not a free product, and 
let $H\ge \operatorname{Inn}(G)$ be a subgroup of $\Aut(G)$. If
all elements of $H$ are realized by orientation-preserving diffeomorphisms of $P$, then the stabilization $\SHAut(*_nG)\to\SHAut(*_{n+1}G)$ induces an isomorphism on $H_i$ for $n\ge 2i+2$ and a surjection for $n=2i+1$.
\end{cor}

The group $H$ can always be taken to be just the inner automorphism group of $G$.  When $G=\Z$, the group $\SHAut(*_nG)$ is then the subgroup of $\SAut(F_n)$ consisting of automorphisms taking each basis element to a conjugate of a basis element. 

At the other extreme, there are many cases when we can take $H=\Aut(G)$, so $\SHAut(*_nG)=\Aut(*_nG)$ if $G\ne \Z$. Some examples of groups $G$ to which this applies are:

\begin{enumerate}[{(a)}]
\item $\Z/2$, $\Z/3$, $\Z/4$, $\Z/6$, with $P$ a lens space. (For a general lens space not all automorphisms of $\pi_1$ are realized by diffeomorphisms --- see for example \cite[Table~3]{Mc-iso}.)

\item The fundamental group of a closed orientable surface of positive genus, with $P$ the product of this surface and an interval. (The full mapping class group of a closed surface is the outer automorphism group of its fundamental group.)

\item The fundamental group of a hyperbolic $3$-manifold $M$ of finite volume and no orientation-reversing isometries. (Here every automorphism of $\pi_1M$ can be realized by a diffeomorphism of $M$, in fact by an isometry, by Mostow rigidity.) 

\end{enumerate}

\noindent
Based on these examples, we make the following conjecture: 

\begin{conj}
For any group $G$, the map $\Aut(*_n G)\to \Aut(*_{n+1}G)$ induces an isomorphism on $H_i$ for $n\ge 2i+2$. 
\end{conj}

\medskip

The main theorem can be specialized to the case $P =D^3$ when connected sum with $P$ just adds a puncture to the manifold without changing the fundamental group. The proof in this case is a lot simpler and yields a slightly better stable range. It also generalizes to manifolds of any dimension $m\ge 2$ in the following way. Let $M$ be a connected $m$-dimensional manifold without any conditions on orientability or compactness. Let $\del_0M$ be a boundary component of $M$ and let 
$\Pts_n=\{p_1,\cdotss,p_n\}$ be a set of $n$ distinct points in the interior of $M$. We denote by $\Ga(M,n,R)$ the group of components of the space of diffeomorphisms of $M$ which map $\Pts_n$ to itself and fix a submanifold $R$ of $\del M$, with $\del_0M\subset R$. When $m=2,3,4$, the group $\Ga(M,n,R)$ is isomorphic to $\Ga(M_n^{D^m},R)$, where $M_n^{D^m}$ denotes $M\# D^m\#\cdots\#D^m$ as before. This follows from the fact that $\Dif(S^{m-1})$ has the homotopy type of the orthogonal group $O(m)$ when $m\le 4$.

There is a map $\Ga (M,n,R)\to \Ga (M,{n+1},R)$ induced by gluing onto $\del_0M$ a ball $D^m$ containing a point $p_{n+1}$ in its interior, identifying a disk $D^{m-1}$ in $\del D^m$ with a similar disk in $\del_0M$. 

\begin{prop}\label{higherm}
For any connected manifold $M$ of dimension $m\ge 2$, the stabilization map
$\Ga (M,n,R)\to \Ga (M,{n+1},R)$ induces an isomorphism on $H_i$ when $n\ge 2i+1$, and a surjection when $n=2i$. 
\end{prop}

When the manifold $M$ has dimension $2$, this result gives a stabilization for mapping class groups of surfaces different
from the classical one considered by Harer \cite{harer}. (A partial result in this direction can be found in \cite{Han}.) 
 Together with the classical genus stability, it shows that mapping class groups of orientable surfaces stabilize with respect to connected sum with any surface. This holds for nonorientable surfaces as well by \cite{NW}. These general stabilization results are the analogs for surfaces of the main stabilization result in this paper for $3$-manifolds.  

\smallskip

Proposition~\ref{higherm} also holds if $\Ga (M,n,R)$ is replaced by $\pi_1\Conf(M,n)$, the fundamental group of the
configuration space of $n$ unordered points in $M$. (This is a subgroup of $\Ga (M,n,R)$ when $\pi_1\Dif(M)=0$.) 
It gives a variant of the classical stability for configuration spaces, namely stability for the homology of the space $\Conf(M,n)$ \cite[Append.]{S}.

When $m\ge 3$, $\pi_1\Conf(M,n)$ is isomorphic to $\pi_1(M)\wr \Si_n$, the wreath product of $\pi_1(M)$ with the
symmetric group. 
Our techniques here apply to the general case
$G\wr \Si_n$ for any group $G$ and we obtain
the following completely algebraic statement:

\begin{prop}\label{wreath}
For any group $G$, the inclusion map $G\wr \Si_n\to G\wr \Si_{n+1}$
induces an isomorphism on $H_i$ when $n\ge 2i+1$ and a surjection when $n=2i$. 
\end{prop}

This seems to have been known for a long time and can be recovered from \cite{Be} (see also \cite{Han}). In the case where $G$ is the trivial group, i.e. for the symmetric groups, Nakaoka \cite{Na} showed that the stabilization map is an isomorphism on $H_i$ for $n\ge 2i$ and is injective for all $i$ and $n$. Granting the injectivity result, his stable range thus agrees with ours.

\smallskip

When $M=D^2$, $\Ga (D^2,n,\del D^2)\cong \pi_1\Conf(D^2,n)$ is the braid group $B_n$, and Proposition~\ref{higherm} recovers
the homological stability of the braid groups \cite{A}. For an arbitrary $S$ with $\del S \ne \emp$, 
$\Conf(S,n)$ is a $K(\pi,1)$ and 
the stability for the surface braid group $B^S_n=\pi_1\Conf(S,n)$ recovers the dimension 2 case of
\cite{S}, Proposition~A.1 (see also \cite{Nap}). The proof also extends easily to wreath products $G\wr B^S_n$, formed using the natural map $B_n^S\to \Si_n$.  

\begin{prop}\label{braidwreath}
For any group $G$, the inclusion map 
$G\wr B^S_n\to G\wr B^S_{n+1}$
induces an isomorphism on $H_i$ when $n\ge 2i+1$ and a surjection when $n=2i$. 
\end{prop}

\medskip

Finally we consider stabilization with respect to boundary connected sum. 
The operation of boundary connected sum for compact, connected orientable $3$-manifolds $M_1$ and $M_2$ consists of identifying a disk $D_1 \subset \del M_1$ with a disk $D_2\subset\del M_2$, producing a manifold $M_1\nat M_2$.  As with ordinary connected sum, $M_1\nat M_2$ depends on choosing orientations for $M_1$ and $M_2$, so we assume this has been done, but $M_1\nat M_2$ also depends on choosing a component $\del_0 M_i$ of $\del M_i$ to contain the disk $D_i$. Hence we consider now manifolds $M$ with a chosen component $\del_0 M$ of $\del M$, whether we indicate this in the notation or not. For simplicity we will also restrict attention here to irreducible manifolds.

Let $M=N\nat P\nat \cdotss \nat P$ be the boundary connected sum of a manifold $N$ and $n$ copies of a manifold $P$. 
Let $R$  be a nonempty union of disjoint disks in $\del_0N$. This plays the role of a set of boundary circles in the stabilization theory for surfaces, or of boundary spheres in the $3$-manifold stabilization described earlier. 
We then have a mapping class group $\Gamma_n^P(N,R)$ consisting of isotopy classes of diffeomorphisms of $M$ fixed on $R$. There is a stabilization $\Gamma_n^P(N,R)\to\Gamma_{n+1}^P(N,R)$ obtained by attaching another copy of $P$ by identifying half of a disk in $\del_0 P$ with half of a disk of a chosen component $R_0$ of $R$. 

\begin{thm}\label{disk-stabilization}
{\rm (i)} For any compact, connected, oriented, irreducible $3$-manifolds $N$ and $P$ with chosen boundary components $\del_0N$ and $\del_0P$, 
the stabilization map $H_i(\G_n^P(N,R)) \to H_i(\G_{n+1}^P(N,R))$ induced by boundary connected sum with $P$ is an isomorphism when $n\ge 2i+2$ and a surjection when $n= 2i+1$. 

\noindent
{\rm (ii)} When $P=S^1\times D^2$ the group $H_i(\G_n^P(N,R))$ is independent of the number of disks in $R$ if $n\ge 2i+2$, and the map $ \G_n^P(N,R)\to \G_n^P(N)$ forgetting $R$ induces an isomorphism on $H_i$ when $n\ge 2i+4$. 
\end{thm}

Every compact connected oriented irreducible manifold $M$ has a decomposition as a sum $P_1\nat \cdotss \nat P_n$ where each $P_i$ is prime with respect to the sum operation, and such a decomposition is unique up to order and
insertion or deletion of trivial summands $D^3$.  This result, which is less standard
than the corresponding result for ordinary connected sum, is a direct corollary of the existence and uniqueness of compression bodies in irreducible manifolds \cite[Thm.~2.1]{Bo}.
As with the usual connected sum, the theorem for $P$ prime with respect to $\nat$-sum implies the result for $P$ not prime.

For $N=D^3$ and $P=S^1\x D^2$, the 
group $\G_g^P(N)$ is the mapping class group $\G(\mathcal{H}_g)$ of a handlebody $\mathcal{H}_g$ of genus $g$. This is a subgroup of the mapping class group of $\del\mathcal{H}_g$. 

\begin{cor}
The homology group  $H_i(\G(\mathcal{H}_g))$ is independent of $g$ when $g\ge 2i+4$. 
\end{cor}

\bigskip

\noindent {\em Organization of the paper\/}:  Section~\ref{mcg} gives the information we will need about mapping class groups of reducible $3$-manifolds, including information on the map $\G(M)\to \Aut(\pi_1(M))$ and the fact that stabilization $M\to M\#P$ induces an injection on mapping class groups.  The proof of the main theorem, Theorem~\ref{main}, is spread over Sections~\ref{connectivity} to \ref{correction}. 
Section~\ref{connectivity} studies a class of simplicial complexes that we call join complexes and gives our main technical result, Theorem~\ref{XtoY}. Section~\ref{complexes} defines simplicial complexes on which mapping class groups of 3-manifolds act and proves high connectivity of these complexes using the results of the previous section. Section~\ref{spectralseq} gives the axiomatization of the spectral sequence arguments and applies this, along with connectivity results from the previous section, to prove part~(i) of Theorem~\ref{main}. Corollaries~\ref{SymAut} and \ref{free} follow using Section~\ref{mcg}. Section~\ref{correction} proves part~(ii) of Theorem~\ref{main}. Section~\ref{otherdim} considers manifolds of any dimension $m\ge 2$ and proves Propositions~\ref{higherm}, \ref{wreath} and \ref{braidwreath}. Finally Section~\ref{boundary} is concerned with boundary connected sums and proves Theorem~\ref{disk-stabilization}.

\section{mapping class groups of nonprime 3-manifolds}\label{mcg}

For a compact connected orientable $3$-manifold $M$ with a basepoint $x$ in the interior of $M$, let $\G(M,x)$ denote the group of path components of the group $\Dif\mkern1mu^+(M,x)$ of orientation-preserving diffeomorphisms of $M$ that fix $x$.  There is then a homomorphism $\Phi\cln\G(M,x)\to\Aut(\pi_1(M,x))$. When $M$ is prime there has been much work done on determining the extent to which $\Phi$ is an isomorphism. We begin this section by describing what is known about $\Phi$, or easily deducible from known results, for a general nonprime $M$. The end of the section is concerned with injectivity of the stabilization map $\G(M,R)\to \G(M\# P,R)$.

\medskip

There are four obvious types of diffeomorphisms that give elements of the kernel of the homomorphism $\Phi\cln\G(M,x)\to\Aut(\pi_1(M,x))$:

\begin{enumerate}[{(1)}]
\item Twists along spheres, supported in a product $S^2\x I \subset M$, rotating the slices $S^2\x \{t\}$ according to a loop generating $\pi_1SO(3) = \Z/2$.

\item Twists along properly embedded disks, supported in a product $D^2\x I \subset M$, rotating the slices $D^2\x \{t\}$ according to a loop generating $\pi_1SO(2) = \Z$.

\item Permuting two boundary spheres of $M$ by a diffeomorphism supported in a neighborhood of the union of the two spheres and an arc joining them.

\item Sliding a boundary sphere of $M$ around a loop in $M$ with its endpoints on the sphere. This gives a diffeomorphism supported in a neighborhood of the union of the loop and the sphere.

\end{enumerate}

\noindent
The following proposition sums up various results spread in the literature over the last 40 years.

\begin{prop}\label{kernel}
The kernel of $\Phi\cln\G(M,x)\to\Aut(\pi_1(M))$ is generated by diffeomorphisms of the types \hbox{\rm (1)-(4)}.  
In particular, $\Phi$ is injective for irreducible manifolds with {\rm (}possibly empty{\rm \/)} incompressible boundary.
\end{prop}

In our applications, we actually consider the map $\Phi'\cln\G(M,B)\to\Aut(\pi_1(M,x))$ with domain the mapping class group of $M$ fixing a ball $B$ around $x$. The result for $\Phi$ is easily seen to imply the corresponding result for $\Phi'$ using the short exact sequence $K\to \G(M,B)\to \G(M,x)$, where the kernel is generated by the (possibly trivial) twist along $\del B$.

\begin{proof}
Consider first the case that $M$ is irreducible. For Haken manifolds that are closed or have incompressible boundary this is a classical result of Waldhausen \cite[Thm.~7.1]{W}. The case of irreducible manifolds with nonempty compressible boundary was shown in \cite[Thm.~6.2.1]{MM}. For closed (non-Haken) hyperbolic manifolds, it is a corollary of the main theorems of \cite{G,GMT} as such manifolds are $K(\pi,1)$'s. According to the geometrization conjecture, this leaves only closed non-Haken Seifert manifolds. 
Those with infinite fundamental group are $K(\pi,1)$'s, and the injectivity of $\Phi$ follows from \cite[Thm.~3]{BO}. Those with finite fundamental group are spherical manifolds and injectivity of $\Phi$ is shown in the proof of \cite[Thm.~3.1]{Mc-iso}.

For reducible manifolds, the result follows now from \cite[Thm.~1.5]{Mc-top}. An alternative 
proof is given in Section~\ref{proof}. 
\end{proof}

There are a few simple but useful observations that can be made about twists along spheres and disks.

\begin{figure}[htp]
\vskip-6pt
\includegraphics[width=0.7\textwidth]{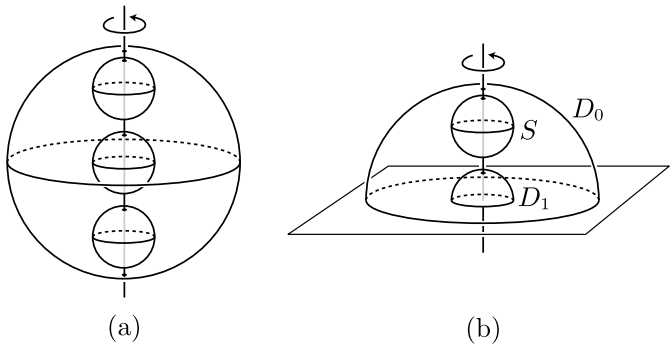}
\vskip-10pt
\caption{Twists along spheres and disks}\label{twists}
\vskip-15pt
\end{figure}

\smallskip
\noindent
(I) \enskip If $D^3_n$ is a $3$-ball $D^3$ with $n$ disjoint sub-balls removed, then for any embedding $D^3_n\subset M$ the composition of the $n+1$ twists along the boundary spheres of $D^3_n$ is isotopic to the identity. To construct an isotopy from the identity to this composition of twists, align the sub-balls along an axis in $D^3$  as in Fig.~\ref{twists}(a) 
and then rotate each point $x$ in the region between the inner and outer boundary spheres of $D^3_n$ by an angle $t\theta(x)$ about the axis, where $t\in I$ is the isotopy parameter and $\theta(x)$ goes from $0$ to $2\pi$ as $x$ varies across an $\epsilon$-neighborhood of the boundary spheres, with $\theta(x)=0$ on the spheres and $\theta(x)=2\pi$ outside the neighborhood.
 From this observation it follows that if $S_1,\cdotss,S_k$ is any maximal sphere system in $M$ then every twist along a sphere in $M$ is isotopic to a composition of twists along a subset of the $S_i$'s, as a twist around a sphere $S$ that intersects $S_1,\cdotss,S_k$ is equal to the product of twists around spheres obtained by surgering $S$ along the intersections.
In particular the subgroup of $\G(M,x)$ generated by twists along spheres is finitely generated, a product of $\Z/2$'s. We will say more about the number of $\Z/2$ factors at the end of this section.

\smallskip
\noindent
(II) \enskip In a similar fashion, if a disk $D_0$ is surgered to produce a disk $D_1$ and a sphere $S$ then a twist along $D_0$ is isotopic to the composition of twists along $D_1$ and $S$ (see Fig.~\ref{twists}(b)). Hence any composition of twists along disks is isotopic to a composition of twists along disks that are disjoint from a given maximal sphere system and twists along spheres in the maximal system. 

We turn now to the image of $\Phi$. If $M$ is the connected sum of prime manifolds $P_i$ then $\pi_1M$ is the free product of the groups $\pi_1P_i$, so we will be concerned with the automorphism group of a free product.  For a free product $G=G_1*\cdotss *G_n$ where no $G_i$ can be decomposed further as a free product, $\Aut(G)$ is generated by the following types of automorphisms (see for example \cite{Gi}):

\begin{enumerate}[{(1)}]
\item Automorphisms of the individual factors $G_i$.
\item Permutations of isomorphic factors.
\item Partial conjugations, in which one factor $G_i$ is conjugated by an arbitrary element $x_j\in G_j$ for some $j$ and the other factors $G_k$, $k\ne i$, are fixed. In particular, inner automorphisms of $G$ can be realized by compositions of such partial conjugations.
\item In case some $G_i$ is infinite cyclic with generator $g_i$ there is an automorphism which sends $g_i$ to $g_ig_j$ or $g_jg_i$ for $g_j$ an arbitrary element of $G_j$, $j\ne i$, and which fixes all other factors $G_k$, $k\ne i$.
\end{enumerate}
The first three types of automorphisms generate the symmetric automorphism group $\SAut(G)$. This can also be described as the subgroup of $\Aut(G)$ consisting of automorphisms that take each $G_i$ to a conjugate of a $G_j$. If no $G_i$ is $\Z$ then $\SAut(G)=\Aut(G)$.

If the prime factors $P_i$ of $M$ have incompressible boundary, then their fundamental groups $\pi_1P_i$ are not decomposable as free products, by the version of Kneser's theorem for manifolds with boundary (see for example \cite[Thm.~7.1]{Hem}). Conversely, if $P_i$ has compressible boundary then either $\pi_1P_i$ splits as a free product or $P_i=S^1\x D^2$.

\begin{prop}\label{image}
The map $\Phi\cln \Gamma(M,x)\to\Aut(\pi_1(M,x))$ is surjective provided that this is true for each prime factor $P_i$ of $M$ and provided also that the following two conditions are satisfied\/{\rm :}  

\begin{enumerate}[{\rm (i)}]
\item  Each $P_i$ which is not closed has incompressible boundary. 

\item  All factors $P_i$ whose fundamental groups are isomorphic are diffeomorphic via orientation-preserving diffeomorphisms.
\end{enumerate}
\end{prop}

Concerning condition (ii), most closed orientable prime $3$-manifolds are determined up to diffeomorphism by their fundamental group, the only exceptions being lens spaces. If nonempty boundaries are allowed there are many more examples, such as products of a circle with surfaces of the same Euler characteristic. If orientations are taken into account the situation is more subtle. There are many orientable prime $3$-manifolds, both with and without boundary, which have no orientation-reversing diffeomorphisms. If $M$ contains copies of these with both orientations, then there will not exist orientation-preserving diffeomorphisms of $M$ that permute all the prime factors with the same fundamental group, even though there are automorphisms of $\pi_1(M)$ that permute the corresponding factors.

\begin{proof}
We can build $M$ as the connected sum of its prime summands $P_1,\cdotss,P_n$ in the following way. Remove the interiors of $n$ disjoint balls from a sphere to produce a manifold $S^3_n$ with boundary spheres $S^2_i$, and remove the interior of a ball from $P_i$ to produce a manifold $P^0_i$ with a new boundary sphere $\del_0P^0_i$. Then glue each $P^0_i$ to $S^3_n$ by identifying $\del_0P^0_i$ with $S^2_i$. We may assume the basepoint $x$ lies in $S^3_n$. 

Automorphisms of types (1) and (2) are realized by elements of $\Gamma(M,x)$ by hypothesis. Type (3) automorphisms can be realized by dragging $P^0_i$ around a loop in $P^0_j \cup S^3_n$, producing a diffeomorphism supported in a neighborhood of the union of this loop and $S^2_i$. Type (4) automorphisms arise when $\pi_1P_i=\Z$, which means that $P_i=S^1\x S^2$ as it cannot be $S^1\x D^2$ by assumption. To realize a type (4) automorphism in that case, observe that the connected sum with $S^1\x S^2$ obtained by attaching $P_i^0$ to $S^3_n$ can be achieved instead by replacing $S^3_n$ by $S^3_{n+1}$ so that the boundary sphere $S^2_i$ is replaced by two boundary spheres which are identified to produce the connected sum with $S^1\x S^2$. Then the automorphism in (4) can be obtained by dragging one of these two spheres around a loop in $P^0_j \cup S^3_{n+1}$.
\end{proof}

When $P_i = S^1 \x D^2$, a type (4) automorphism cannot be realized by a diffeomorphism. To see why, note first that the generator $g_i$ of $\pi_1(P_i)$ can be realized by a loop consisting of a circle in the torus $\del P_i$ joined to the basepoint in $S^3_n$ by an arc. A diffeomorphism of $M$ must take $\del P_i$ to a torus $\del P_j$ for a summand $P_j=S^1\x D^2$, so the diffeomorphism must take the loop representing $g_i$ to a similar loop representing a conjugate of $g_j$ or $g_j^{-1}$.

When all the summands $P_i$ of $M$ are $S^1\x D^2$, the image of $\Phi$ is contained in $\SAut(F_n)$ by the argument in the preceding paragraph, and in fact the image of $\Phi$ is all of $\SAut(F_n)$ since the nontrivial automorphism of $\pi_1P_i=\Z$ can be realized by an orientation-preserving diffeomorphism of $S^1 \x D^2$.

\bigskip 
\noindent
{\bf Injectivity of Stabilization.} 

\smallskip
In the preceding section we defined stabilization maps $\G^P_n(N,R)\to\G^P_{n+1}(N,R)$ and $ A^P_n(N,R,T)\to A^P_{n+1}(N,R,T)$, and we will need to know that these are injective.

More generally, let $M$ and $P$ be compact orientable $3$-manifolds with $P$ prime, and let $\del_0M$ be a component of $\del M$.  We stabilize via an inclusion $ M \to M\#P$ obtained by gluing $P^0$ to $M$ by identifying a disk in $\del_0P^0$ with a disk in $\del_0M$. Let $R$ be any compact subsurface of $\del(M\# P)$ containing the component $\del_0(M\# P)$ corresponding to $\del_0 M$, and let $T$ be any compact subsurface of $\del(M\# P)-R$. We also use $R$ and $T$ to denote the restrictions of these surfaces to $\del M$, where $\del_0M \subset R$ in this case. Then we have stabilization maps $\G(M,R)\to\G(M\#P,R)$ and $A(M,R,T)\to A(M\#P,R,T)$.

\begin{prop}\label{stabilizationmap}
The stabilization map $\G(M,R)\to\G(M\#P,R)$ is injective, and the same is true for the quotient stabilization $A(M,R,T)\to A(M\#P,R,T)$.
\end{prop}

As shown in Proposition~\ref{kernel}, $A(M,R,T)$ can often be identified with a subgroup of $\Aut(\pi_1M)$ (and similarly for
$M\#P$). In this case, it is easy to see that the stabilization $A(M,R,T)\to A(M\#P,R,T)$ is injective since the stabilization $\Aut(\pi_1M)\to\Aut(\pi_1M * \pi_1P)$ is obviously injective. Some instances of this are the groups in Corollaries~\ref{SymAut} and \ref{free}, as well as $\Aut(F_n)$. Injectivity for the stabilizations $A^s_{n,k}\to A^s_{n+1,k}$ and $A^s_{n,k}\to A^s_{n,k+1}$ for $s\ge 1$ can also be deduced algebraically from injectivity on $\Aut(\pi_1)$.

\begin{proof}  We first prove the proposition in the case that $\del_0 M$ is a sphere. This corresponds to a $D^3$
  summand of $M$.  If the other prime summands are $P_1,\cdotss,P_n$ then we can construct $M$ from $D^3_n$, a
 ball with the interiors of $n$ disjoint balls in its interior removed, by attaching the punctured prime manifolds
 $P^0_1,\cdotss,P^0_n$ to $D^3_n$ by identifying their boundary spheres $\del_0P^0_i$ with the corresponding
 interior boundary spheres $S^2_i$ of $D^3_n$.   

Represent an element of the kernel of $\G(M,R)\to \G(M\#P,R)$ by a diffeomorphism $f$ of $M$ fixed
on $R$. The spheres $S^2_i=\del_0P^0_i$ are then isotopic in $M\#P$ to their images
$f(S^2_i)$. Since $M$ is a retract of $M\#P$, these isotopies in $M\#P$ can be replaced by
homotopies in $M$. Then by Laudenbach's homotopy-implies-isotopy theorem \cite[Thm.~III.1.3]{L}
there is an isotopy of $\cup_i f(S^2_i)$ to $\cup_iS^2_i$ in $M$. After extending this isotopy to an
isotopy of $f$ we may assume that $f(S^2_i)=S^2_i$ for each~$i$, and hence also
$f(P^0_i)=P^0_i$. After a further isotopy we can arrange that $f$ is the identity on $D^3_n$ since
any orientation-preserving diffeomorphism of $D^3_n$ that is fixed on the outer boundary sphere and
takes each of the other boundary spheres $S^2_i$ to itself is isotopic, through such
diffeomorphisms, to the identity. Thus we have shown that the kernel of the stabilization map
$\G(M,R)\to\G(M\#P,R)$ is in the image of the natural map $\Psi\cln\prod_i \G(P^0_i,R^0_i) \to
\G(M,R)$, where $R^0_i=(R\cap P^0_i)\cup\del_0P^0_i$. 

Now we come to the key step in the proof, the fact that $\Psi$ is injective. When $R=\del M$ this is an immediate
consequence of the main theorem in \cite{HL}.  To deduce injectivity of $\Psi$ for general $R$ from the special case
that $R=\del M$, consider the fibrations obtained by restricting diffeomorphisms to the boundary:  
$$  
\xymatrix{\prod_i \Dif(P_i^0,\del P_i^0) \ar[r]\ar[d] & \prod_i \Dif(P_i^0,R_i^0) \ar[r]\ar[d] & \prod_i\Dif(\del P_i^0,R_i^0) \ar[d]\\
\Dif(M,\del M) \ar[r] & \Dif(M,R) \ar[r] & \Dif(\del M,R)}
$$
The vertical map on the right is an inclusion onto a union of components, the components that take $\del P^0_i$ to itself for each ~$i$. Hence this vertical map induces an injection on $\pi_0$ and an isomorphism on $\pi_1$ in the following diagram of exact sequences of homotopy groups: 
$$\xymatrix{\prod_i\pi_1\Dif(\del P^0_i,R^0_i) \ar[r] \ar[d] & \prod_i\G(P^0_i,\del P^0_i)\ar[r]\ar[d]^\Psi & \prod_i\G(P^0_i,R^0_i)  \ar[r] \ar[d]^\Psi & \prod_i \G(\del P^0_i,R^0_i) \ar[d]\\
\pi_1\Dif(\del M,R) \ar[r] & \G(M,\del M) \ar[r] & \G(M,R) \ar[r] & \G(\del M,R) }$$
By exactness, injectivity of the first $\Psi$ then implies injectivity of the second $\Psi$.

In the rest of the proof all diffeomorphisms and isotopies will be understood to fix $R$, whether we mention this explicitly or not. 

Injectivity of $\Psi$ can be restated as saying that if two diffeomorphisms of $M$ that are the identity on $D^3_n\cup
R$ are isotopic fixing $R$, then they are isotopic fixing $D^3_n\cup R$. The corresponding statement for $M\#P$ also holds.  
We apply this to
the stabilization $f\# id_P$ of the earlier diffeomorphism $f$ of $M$: From above, we have that $f\# id_P$ fixes $D^3_{n+1}$
and is isotopic to the identity. Hence it is isotopic to the identity fixing $D^3_{n+1}$. Restricting this isotopy to the summands $P^0_i$ of
$M$ shows that $f$ is also isotopic to the identity (fixing $D^3_n$), which 
proves injectivity of $\G(M,R)\to\G(M\#P,R)$ in the special case $\del_0M=S^2$.

The injectivity of $A(M,R,T)\to
A(M\#P,R,T)$ is equivalent to the statement that if a diffeomorphism $f$ of $M$
stabilizes to a diffeomorphism of 
$M\#P$ that is isotopic to a product of twists, then $f$ itself is isotopic to a product of twists.  
Recall from the beginning of the section that twists act trivially on isotopy classes of embedded spheres, so that $f(S^2_i)$ is again
isotopic to $S^2_i$ for each $i$. Hence, by the same argument as before,  
we can again assume $f$ fixes $D^3_n$.  
Then we can use injectivity of $\Psi$
just as before:  If the stabilization of $f$ is isotopic to a product of twists in $M\#P$, then these twists can be chosen to
be supported in the complement of $D^3_{n+1}$ and the isotopy can be taken to be the identity on $D^3_{n+1}$. 
By restriction to the summands $P^0_i$ of $M$, this implies that $f$ itself is isotopic to a product of twists.

Next we treat the case that $\del_0M$ is not a sphere but there is some other component $\del_1 M$ of $\del M$ that is a
sphere, with $\del_1M\subset R$. 
As before, we can start with a diffeomorphism $f$ of $M$ that is the identity on $D^3_n$, where $\del_1M$ is the outer
boundary sphere of $D^3_n$. We want to arrange that $f\# id_P$ is also the identity on $D^3_{n+1}$. 
We index the summands $P_i$ so that $\del_0M$ is a component of $\del P_n$.  Choose an
arc $a$ in $P^0_n$ from a point in $\del_0P^0_n=S^2_n$ to a point in the interior of the disk in $\del_0M$ where $P^0$
attaches (see Fig.~\ref{bdySph}). Since $f$ fixes the spheres $S^2_i$ and
induces the identity on $\pi_1M$, the arcs $a$ and $f(a)$ are homotopic in $P^0_n$, fixing their
endpoints. These arcs have one endpoint on a boundary sphere of $P^0_n$, the sphere $\del_0P^0_n$, so it follows from
the `lightbulb trick', which is explained after the conclusion of the proof, 
that $a$ and $f(a)$ are in fact isotopic
in $P^0_n$, fixing their endpoints. Extending this isotopy to $P^0_n$, we may thus assume that $f$ is the identity on
$a$. 

\vskip-4pt
\begin{figure}[htp]
\includegraphics[width=0.88\textwidth]{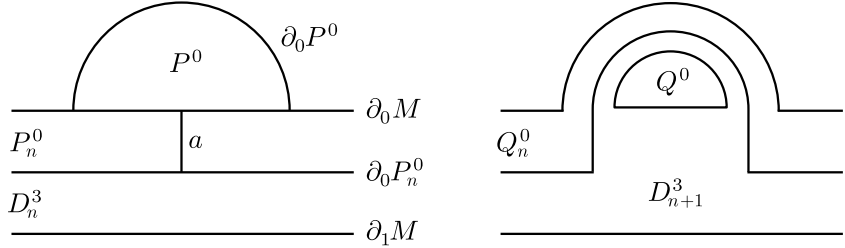}
\vskip-6pt
\caption{Stabilization}\label{bdySph}
\vskip-15pt
\end{figure}

It is possible that $f$ might twist the normal bundle to $a$. This would happen for example if $f$ was a twist along $\del_0P^0_n$ with $a$ as the axis of the twist.  Call this twist~$g$.  After composing $f$ with some power of $g$ we can isotope this composition to be the identity in a tubular neighborhood $V$ of $a$ in $P^0_n$, as well as on $D^3_n$.  Since $g$ has order $2$ in the mapping class group, this means that we can assume that either $f$ or $fg$ is the identity in $V$. Call these two cases (i) and (ii). We can choose $V$ to include the disk where $P^0$ attaches to $\del_0M$, so $V$ is the part of $P^0_n$ shown in the figure. 

We enlarge $D^3_n$ to a copy of $D^3_{n+1}$ in $M\#P$ by adding a smaller tubular neighborhood of $a$ in $V$ together with a tubular neighborhood of a disk in $P^0$ parallel to the disk $\del_0P^0 - \del_0M$ as shown in the second half of the figure. One complementary component of $D^3_{n+1}$ in
$M\#P$ is then a smaller copy $Q^0$ of $P^0$ and another complementary component is a new copy $Q_n^0$ of $P^0_n$. 
In case (i) $f$ is
the identity on $D^3_{n+1}$ and the argument in the earlier situation that $\del_0M=S^2$ applies to finish the proof.  
  For case~(ii), as $f\# id_P$ is isotopic to the identity, it follows by composition with $g$ that $fg\# id_P$ is isotopic to $g\# id_P$. 
The diffeomorphism $g\# id_P$ is not itself the identity on $D^3_{n+1}$, but it is isotopic to a diffeomorphism $h$
which fixes $D^3_{n+1}$. Indeed, $g$ is defined as a twist along the sphere $\del_0P^0_n$ in $M$, which after stabilization is isotopic
to the product $h$ of the twists along $\del_0Q^0$ and along $\del_0Q_n^0$, using the thrice-punctured $S^3$
formed by the union of $D^3_{n+1} - D^3_n$ with a neighborhood of  $\del_0P^0_n$. By injectivity of $\Psi$,
$fg\# id_P$ is thus isotopic, fixing $D^3_{n+1}$, to $h$. 
This isotopy from $fg\#id_P$ to $h$, restricted to $Q^0_n$, is conjugate, via a diffeomorphism $P^0_n\approx Q^0_n$ that is the identity outside $V$, to an isotopy from $fg$ to $g$ on $P^0_n$ fixing $\del_0P^0_n$.  Hence $f$ is isotopic to the identity on $P^0_n$ fixing $\del_0P^0_n$.  On the other $P^0_i$'s the isotopy from $fg\#id_P$ to $h$ gives an isotopy from $f$ to the identity, so altogether $f$ is isotopic to the identity on $M$.

This finishes the proof for $\G(M,R)\to\G(M\# P,R)$ when $R$ contains a sphere component.
The injectivity of the stabilization $A(M,R,T)\to A(M\#P,R,T)$ follows in the same way as in the case where $\del_0M$
was a sphere since if twists along spheres are factored out there is no difference between cases (i) and (ii).

There remains the case that no component of $R$ is a sphere. Let $f$ be a diffeomorphism of $M$ that is isotopic to the identity or a product of twists after stabilization. Both $f$ and the isotopy are fixed on $R$, so we may assume they are fixed on a neighborhood of $R$.  Choose a ball $B$ in this neighborhood, away from where $P^0$ is attached to $M$. Deleting the interior of $B$ from $M$ gives a submanifold $M'\subset M$ on which $f$ is a diffeomorphism whose stabilization to $M'\#P$ is isotopic to the identity or a product of twists, fixing the new boundary sphere $\del B$. Having this boundary sphere, we can deduce from the previous case that the restriction of $f$ to $M'$ is isotopic, fixing $R$ and $\del B$, to the identity or a product of twists. Filling $B$ back in, it follows that $f$ itself is isotopic to the identity or a product of twists.
\end{proof}

\smallskip
\noindent{\em The Lightbulb Trick.} This is a classical technique for avoiding issues of knotting and linking of arcs in $3$-manifolds in one very special situation.  Consider arcs in a $3$-manifold that connect one boundary component to another boundary component that is a sphere. The statement is that if two such arcs have the same endpoints and are homotopic fixing the endpoints, then they are isotopic fixing the endpoints. It is a standard fact that a homotopy of an arc $a$ can be replaced by a deformation that is an isotopy except for finitely many transverse crossings where one subarc $\alpha$ of $a$ passes through another subarc $\beta$, intersecting it transversely in one point at one instant of time.  If one end of $a$ lies on a boundary sphere, such a crossing can be replaced by an isotopy where the middle portion of $\alpha$ is dragged along $a$ to get near the sphere, then is dragged around the sphere to the other side of $a$, and then is dragged back to the other side of $\beta$ without crossing $a$ at any time during the process, assuming that $\beta$ is closer to the sphere along $a$ than $\alpha$ is, which can always be arranged by interchanging $\alpha$ and $\beta$ if necessary.  The same technique can be used to improve homotopies to isotopies for a collection of disjoint arcs $a_i$ joining distinct boundary spheres $S_i$ to boundary points not on any $S_i$.

A variant of this construction, with a weaker hypothesis and a weaker conclusion, will be used in the proof of Lemma~\ref{action}. To distinguish this from the lightbulb trick we call it the {\em balloon trick}.  As before, consider arcs $a$ with one endpoint on a boundary sphere and the other endpoint on another boundary component. Then for any two such arcs having the same endpoints, there is a diffeomorphism of the manifold, fixed on the boundary, taking one arc to the other arc.  To prove this, first fill in the boundary sphere with a ball. Then in the resulting enlarged manifold one can shrink the arc to its endpoint in the other boundary component, dragging the ball along, until the arc and the ball are in a standard configuration near the boundary of the manifold.  Such an isotopy can be extended to an ambient isotopies of the enlarged manifold.  Then if one first performs the ambient isotopy for one arc followed by the reverse of the ambient  isotopy for another arc, the result is a diffeomorphism of the original manifold taking one arc to the other.  A similar construction works also for collections of arcs.

\begin{rem}{\rm
We can describe more explicitly the subgroup of $\G(M,R)$ generated by twists along spheres. (We will not make essential use of this information in the rest of the paper.)  This subgroup is contained in the subgroup $\prod_i\G(P^0_i,R^0_i)$ of $\G(M,R)$ so it suffices to see what the subgroup of  $\G(P^0_i,R^0_i)$ generated by twists along spheres is for each $P^0_i$. For notational convenience we drop the subscript $i$. First assume that $P\ne S^1 \x S^2$ so the only nontrivial sphere in $P^0$ is $\del_0P^0$, up to isotopy. There is a fibration  $\Dif(P^0,R^0) \to \Dif(P,R)\to E$ where $E$ is the space of embeddings of $D^3$ in $P$. Part of the long exact sequence of homotopy groups for this fibration is
$$
\pi_1\Dif(P,R)\to \pi_1 E\to \G(P^0,R^0)
$$
Since orientable $3$-manifolds have trivial tangent bundle, $E$ is homotopy equivalent to the product $P \x O(3)$, with fundamental group $\pi_1P\x \Z/2$. The $\pi_1P$ factor maps to $\G(P^0,R^0)$ as the inner automorphisms of $\pi_1P^0=\pi_1P$ and the $\Z/2$ maps to the twist along the sphere $\del_0P^0$, so the twist is trivial in $\G(P^0,R^0)$ exactly when this $\Z/2$ in $\pi_1E$ comes from an element of $\pi_1\Dif(P,R)$. One situation where there is such an element is if there is an action of $S^1$ on $P$ fixing a circle pointwise and rotating the normal bundle of this circle. This happens for example if $P=S^1\x D^2$ with $R=\emp$. It also happens if $P$ is a lens space since a lens space is the union of two solid tori and the rotation of one solid torus about its core circle always extends to a rotation of the other solid torus. The group $\pi_1\Dif(P,R)$ is trivial if $R\ne\emp$, and when $R=\emp$ it is trivial for hyperbolic manifolds \cite{G2} and for Haken manifolds which are not Seifert-fibered \cite{H1}. By the geometrization conjecture, this leaves only Seifert manifolds to consider. For those which are Haken manifolds the group $\pi_1\Dif(P)$ is known in very explicit terms \cite{H1}, so it would be an exercise to compute the map $\pi_1\Dif(P)\to\pi_1E$ in these cases. As far as we know, the group $\pi_1\Dif(P)$ has not yet been computed for all the remaining small Seifert manifolds, for example for the Poincar\'e homology sphere. }

{\rm For $P=S^1\x S^2$ the argument for lens spaces shows that the twist along $\del_0 P^0$ is trivial in $\G(P^0,R^0)$. The only other twist to consider is the twist along the $S^2$ factor, and it is classical that this is nontrivial in $\G(P)$, hence also in $\G(P^0,R^0)$.
}
\end{rem}

\section{combinatorial connectivity results}\label{connectivity}

In this purely combinatorial section we define the notion of a join complex over a simplicial complex $X$ and prove the main technical result of the paper, Theorem~\ref{XtoY}, which says that a join complex over $X$ is highly connected whenever $X$ and certain of its subcomplexes are highly connected. The proof uses a coloring lemma, which we present first. 

\medskip

Let $E$ be a finite set and consider simplicial complexes $K$ with vertices labeled by elements of $E$. 
We call a simplex of $K$ {\em good\/} if all its vertices are labeled by different elements of $E$, and {\em bad\/} if each of its vertex labels occurs at least twice.  Each simplex is in a unique way the join of a good simplex and a bad simplex, where these two simplices have disjoint sets of labels.

\begin{lem}[Coloring lemma]\label{coloring}
Let a triangulation of $S^k$ be given with its vertices labeled by elements of a set $E$ having at least $k+2$ elements. Then this labeled triangulation can be extended to a labeled triangulation of $D^{k+1}$ whose only bad simplices lie in $S^k$, and with the triangulation of $S^k$ as a full subcomplex. The labels on the interior vertices of $D^{k+1}$ can be chosen to lie in any subset $E_0\subset E$ with at least $k+2$ elements.
\end{lem}

\begin{proof}
We prove the lemma by induction on $k$, starting with the trivial case $k=-1$ when $D^{k+1}$ is a point that can be labeled by any element of $E_0$. Suppose we have proved the lemma in all dimensions less than $k$ and consider a triangulation of $S^k$ labeled by $E$. Triangulate the disk $D^{k+1}$ bounded by $S^k$ by putting a vertex at the center of $D^{k+1}$ and coning off the triangulation of $S^k$ to this vertex. Label this vertex by any element of $E_0$. We are going to modify the triangulation in the interior of $D^{k+1}$, adding vertices labeled by elements of $E_0$, until all the bad simplices lie in $S^k$. 

If there are bad simplices in $D^{k+1}$ not contained in $S^k$, let $\s$ be one of maximal dimension $p$. We must have $p>0$ since vertices are never bad. Denote by $E_\s\subset E$ the set of labels occurring in $\s$. 
The link of $\s$ is a sphere $S^{k-p}$ since $\s$ is not contained in $S^k$.
 The simplices of the link are labeled by elements of $E-E_\s$ by maximality of $\s$. 
As $k-p<k$ we may apply the induction hypothesis to the link using labels from $E_0-E_\s\subset E-E_\s$ since $|E_0-E_\s|\ge k-p+2$ as $|E_\s|\le p$. This gives a triangulation of the disk $D^{k-p+1}$ bounded by the link of $\s$. The star of $\s$ is the join of $\s$ with its link. 
The triangulation of $D^{k-p+1}$ extends to a new triangulation of the star of $\s$ by joining with $\del\s$, and this new triangulation agrees with the old one on the boundary of the star. 
A simplex of the new triangulation has the form $\tau * \mu$ with $\tau$ a face of $\s$ and $\mu$ a simplex of the disk $D^{k-p+1}$, with $\tau$ or $\mu$ possibly empty. 
If such a simplex $\tau *\mu$ is bad, $\mu$ must be empty, since otherwise, as the labels on $\mu$ and $\tau$ are disjoint, $\mu$ would be a bad simplex, hence contained in $S^{k-p}$, contradicting the maximality of $\sigma$.  Thus we have reduced the number 
of bad interior simplices of $D^{k+1}$ of top dimension. This gives the induction step.

The fullness condition holds since it held after the initial coning, and for the induction step, if the vertices of a simplex $\tau *\mu$ in the star of $\sigma$ as above all lie in $S^k$ then $\mu$ must lie in $\del D^{k-p+1}$ by the inductive hypothesis of fullness, hence $\tau *\mu$ must lie in the boundary of the star where the induction step does not change the triangulation so $\tau *\mu$ lies in $S^k$ by induction. 
\end{proof}

\begin{Def}\label{joincpx}
{\rm A {\em join complex\/} over a simplicial complex $X$ is a simplicial complex $Y$ together with a simplicial map $\pi\cln Y\to X$, thought of as a projection, satisfying the following properties:}
\end{Def}
\vspace{-7pt}
\begin{enumerate}[{(1)}]
\item $\pi$ is surjective.

\item $\pi$ is injective on individual simplices.

\item For each $p$-simplex $\s =\lgl x_0,\cdotss,x_p\rgl$ of $X$ the subcomplex $Y(\s)$ of $Y$ consisting of all the $p$-simplices that project to $\s$ is the join $Y_{x_0}(\s)*\cdotss *Y_{x_p}(\s)$ of the vertex sets $Y_{x_i}(\s)=Y(\s)\cap \pi^{-1}(x_i)$. 

\end{enumerate}
Note that $Y(\s)$ need not be equal to $\pi^{-1}(\s)$.  If all the inclusions $Y_{x_i}(\s)\subset \pi^{-1}(x_i)$ are equalities we call $Y$ a {\em complete\/} join complex over $X$.

A reformulation of (3) which we will sometimes use in place of (3) is the following condition:
\begin{list}{}{\setlength{\leftmargin}{22pt}\setlength{\labelwidth}{16pt}\setlength{\labelsep}{5pt}\setlength{\itemsep}{3pt}}

\item[(3$'$)]  A collection of vertices $(y_0,\cdotss,y_p)$ of $Y$ spans a $p$-simplex of $Y$ if and only if for each $y_i$ there
exists a $p$-simplex $\s_i$ of $Y$ such that $y_i\in \s_i$ and $\pi(\s_i)=\lgl \pi(y_0),\cdotss,\pi(y_p)\rgl$. 
\end{list}
\smallskip
\noindent
Clearly (3) implies (3$'$), and the converse follows from the fact that for a vertex $x$ of a $p$-simplex $\s$ of $ X$,  $Y_{x}(\s)$ is the set of vertices $y\in \pi^{-1}(x)$ that are a vertex of at least one $p$-simplex of $Y$ projecting onto $\sigma$.

\begin{ex}[Labeling systems] 
{\rm Given a simplicial complex $X$, define a {\em
    labeling system} for $X$ to be a
  collection of nonempty sets $L_x(\s)$ for each simplex $\s$ of $X$ and each vertex $x$ of $\s$, satisfying $L_x(\tau)\supset L_x(\s)$
  whenever $x\in\tau\subset \s$. 
We can use $L$ to define a new simplicial complex $X^L$ having vertices the pairs $(x,l)$ with
  $x\in X$ and $l\in L_x(\lgl x\rgl)$. A collection of pairs 
  $((x_0,l_0),\cdotss,(x_p,l_p))$ then forms a $p$-simplex of $X^L$ if and only if
  $\s=\lgl x_0,\cdotss,x_p\rgl$ is a $p$-simplex of $X$ and $l_i\in L_{x_i}(\s)$ for each $i$. The natural map $\pi\cln X^L\to X$ with
  $\pi(x,l)=x$ represents $X^L$ as a join complex over $X$.   The set $\pi^{-1}(x)\cong L_x(\lgl x\rgl)$ is viewed as a set of `labels'
  of $x$, and $(X^L)_{x_i}(\s)\cong L_{x_i}(\s)$ is the set of labels of $x_i$ `compatible with $\s$'.  In fact, any join complex is
  isomorphic to one of this form. 
}\end{ex}

Before stating our first result about join complexes we need a couple definitions. A complex is called {\em $d$-spherical} if it is of dimension $d$ and is $(d-1)$-connected. 

\begin{Def}
{\rm A simplicial complex $X$ is {\em Cohen-Macaulay of dimension $n$} if $X$ is $n$-spherical and the link of each $p$-simplex of $X$ is $(n-p-1)$-spherical. More generally, we say $X$ is {\em weakly Cohen-Macaulay of dimension~$n$} if $X$ is $(n-1)$-connected and the link of each $p$-simplex of $X$ is $(n-p-2)$-connected. }
{\rm We often shorten {\em Cohen-Macaulay\/} to {\em CM} and {\em weakly Cohen-Macaulay\/} to {\em wCM}. }
\end{Def}

This definition of a Cohen-Macaulay complex coincides with \cite[Sec.~8]{Qu} and \cite[II.4.2]{St}. The connectivity condition on the links in a wCM complex $X$ of dimension $n$ implies that $\dim X\ge n$ (using the fact that `$k$-connected implies `nonempty' when $k\ge -1$ since it includes the condition that every map from $S^{-1}=\emp$ extends to a map from $D^0$). 
The further dimension condition in the strengthening of wCM to CM is irrelevant for us, but most of the complexes we consider with the CM-connectivity conditions do also satisfy the dimension condition.

\begin{prop}\label{label}
If \,$Y$ is a complete join complex over a CM {\rm (}resp.~ wCM{\rm \/)} complex $X$ of dimension $n$, then $Y$ is also CM {\rm (}resp.~wCM{\rm \/)} of dimension $n$.
\end{prop}

\begin{proof} 
Let $\pi\cln Y\to X$ be the projection. For each vertex $x$ of $X$ choose a lift $s(x) \in \pi^{-1}(x)$. Call a simplex of $Y$ `good' if all its
vertices  are lifts $s(x)$, and call it `bad' if none of its vertices are of that type. The good
simplices form a copy $X^{s}$ of $X$ in $Y$ giving a cross section of the projection $Y\to X$. To show that $Y$ is $(n-1)$-connected, it suffices to show that every map $f\cln S^k\to Y$ with
$k\le n-1$ can be deformed by a homotopy to have only good simplices in its
image. We can assume $f$ is simplicial with respect to some triangulation of $S^k$.  Let $\s$ be a simplex in this triangulation of maximal dimension $q$ such that $f(\s)$ is a bad simplex of $Y$, say of dimension $p\le q$. By maximality, $f$ maps the
link of $\s$ to $X^s$ and in fact to the link of $\pi f(\s)$ in
$X\cong X^{s}$. This link is $(n-p-2)$-connected by assumption and
$\link_{S^k}(\s)\simeq S^{k-q-1}$. As $k-q-1\le n-p-2$, there exists a map
$F\cln D^{k-q}\to \link_X\bigl(\pi f(\s)\bigr)$ extending
$f\big|_{\link(\s)}$. Define $F^s$ to be $F$ followed by the inclusion $X\cong
X^s\inc Y$. As $Y$ is a complete join complex over $X$, we have $F^s\cln D^{k-q}\to
\link_{Y}\bigl(f(\s)\bigr)$, hence $H^s=F^s*f\big|_\s$ maps the ball $D^{k-q}*\s$ to
simplices of $Y$. This ball has boundary $\del D^{k-q}*\s\cup D^{k-q}*\del \s$ and $H^s$ defines a homotopy between $f$ restricted
to $\operatorname{Star}(\s)=\del D^{k-q}*\s$ 
and $F^s* f\big|_{\del\s}$ on
$D^{k-q}*\del\s$. This defines a new map $f'$ homotopic to $f$ and with fewer bad simplices of maximal dimension in its image. Finitely many iterations of this step finish the proof that $Y$ is $(n-1)$-connected.

For each $p$-simplex $\s$ of $Y$ we have $\link_{Y}(\s) =\pi^{-1} (\link_X (\pi(\s)))$.  This complex is a complete join complex over $\link_X(\pi(\s))$ so it is $(n-p-2)$-connected by the above and the wCM assumption on $X$. Thus $Y$ is wCM. As the complexes $X$ and $Y$ have the same dimension, it follows that $Y$ is CM if $X$ is CM.
\end{proof}

The main technical result of this paper gives a weaker connectivity bound for a join complex when the completeness condition is replaced by a more subtle condition:

\begin{thm}\label{XtoY}
Let $Y$ be a join complex over a wCM complex $X$ of dimension $n$. Suppose that  the projected link $\pi(\link_Y(\s))$ of each
$p$-simplex $\s$ of $Y$ is wCM of dimension $n-p-2$.  Then $Y$ is $(\frac{n}{2}-1)$-connected. 
\end{thm}

Both the hypothesis and the conclusion are weaker than what one might have expected, which would have the numbers $n-p-1$ instead of $n-p-2$ in the hypothesis and $n-1$ instead of $\frac{n}{2}-1$ in the conclusion.  In our applications the stronger hypothesis will in fact be satisfied. But even with the stronger hypothesis the conclusion cannot be improved, as the next example shows.

\begin{ex}
{\rm When $n=1$, take $Y$ to be the union of two disjoint $1$-simplices with $\pi\cln Y\to X$ the quotient map identifying an endpoint
  of one $1$-simplex with an endpoint of the other $1$-simplex. The stronger hypotheses of the
  theorem are satisfied (the projections of the links of vertices of $Y$ are nonempty, hence $(-1)$-connected), and $Y$ has connectivity $\frac{n}{2}-1=-\frac{1}{2}$ but not $n-1=0$. The unreduced suspension  $\Sigma Y\to \Sigma X$ gives a $2$-dimensional example with $\Sigma Y$ having connectivity $\frac{n}{2}-1=0$ but not $n-1=1$. For higher-dimensional examples, take the  join of $k+1$ copies of the $1$-dimensional example with itself. This gives spaces $Y^{*k}$ and $X^{*k}$ of dimension $n=2k+1$ with $Y^{*k}$ homotopy equivalent to $S^k$ since $Y$ is homotopy equivalent to $S^0$, so $Y^{*k}$ has connectivity $\frac{n}{2}-1=k-\frac{1}{2}$ but no more.  The stronger hypotheses of the theorem are satisfied because of the following easily-verified facts which imply that the stronger hypotheses are preserved under joins: (1) The join $Y_1*Y_2\to X_1*X_2$ of two join complexes is a join complex. (2) The projected links of a join are the joins of the projected links. (3) The join of wCM complexes of dimensions $l$ and $m$ is a wCM complex of dimension $l+m+1$. (In checking these statements it is convenient to regard the empty set as a simplex of dimension $-1$.) Thus we obtain examples in all odd dimensions. For even-dimensional examples one can take the unreduced suspensions of the odd-dimensional examples, since $\Sigma(Y^{*k})$ is homotopy equivalent to $S^{k+1}$ so it has connectivity exactly $\frac{n}{2}-1=k$.  }
\end{ex}

\smallskip
For the proof of the theorem we will be applying the coloring lemma where the `colors' are vertices of $X$. It will be advantageous to maximize the number of colors available at a given time, so we make a preliminary digression to prove a lemma that will aid in this maximization.

Let $X$ be a simplicial complex. The barycentric subdivision of $X$ is the simplicial complex associated to the poset of
simplices of $X$ ordered by inclusion.  Thus the $p$-simplices of the barycentric subdivision correspond to chains of
inclusions of $p+1$ simplices of $X$.  Let $X_m$ be the subcomplex of the barycentric subdivision corresponding to the
poset of simplices of $X$ with at least $m$ vertices. In particular $X_1$ is the barycentric subdivision of $X$, and
$X_m$ for $m\ge 2$ is homotopy equivalent to the complement of the $(m-2)$-skeleton of $X$.

\begin{lem}\label{genlargement} 
If $X$ is wCM of dimension $n$ then $X_m$ is $(n-m)$-connected. 
\end{lem}

\begin{proof} This will be proved by induction on $m$ using a links argument. The induction starts with the case $m=1$
  which is true by assumption. For $m>1$ suppose we start with a map $f\cln S^k\to X_m$ representing an element of
  $\pi_kX_m$ with $k\le n-m$. By the induction hypothesis this extends to a map $f\cln D^{k+1}\to X_{m-1}$ which we may
  take to be simplicial for some triangulation of $D^{k+1}$. We wish to eliminate the simplices $\s$ of $D^{k+1}$ which
  are `bad' in the sense that $f$ maps each of their vertices to simplices of $X$ with only $m-1$ vertices. As a
  $p$-simplex in $X_{m-1}$ is a chain of inclusions of $p+1$ simplices of $X$ each of which has at least $m-1$ vertices,
  $f$ must in fact be constant on bad simplices. Moreover, such simplices must lie strictly in the interior of
  $D^{k+1}$. If $\s$ is a bad simplex of maximal dimension $p$ then $f$ maps the link of $\s$, a sphere $S^{k-p}$, to
  the subcomplex of $X_m$ of systems containing $f(\s)$ as a strict subsystem. We can regard the restriction of $f$ to
  the link as a map $f_\s\cln S^{k-p}\to (\operatorname{Link}_X(f(\s)))_1$ where the latter complex is the first
  subdivision of the link of $f(\s)$ in $X$. By assumption, this link is $(n-m)$-connected as $f(\s)$ is an
  $(m-2)$-simplex of $X$. Since $k\le n-m$ we have $k-p\le n-m$ 
 so we can extend $f_\s$ to a map $g_\s\cln D^{k-p+1}\to (\operatorname{Link}_X(f(\s)))_1$. This allows us to redefine $f$ in the interior of the star of $\s$ by rewriting $\operatorname{Star}(\s)=\operatorname{Link}(\s)*\s$ as $D^{k-p+1}*\del \s$ and replacing $f$ in this join by the join of $g_\s\cup f(\s)$ 
on $D^{k-p+1}$ and $f$ on $\del\s$.
This eliminates $\s$ as a bad simplex of maximal dimension without introducing any other bad simplices of this dimension or larger. Finitely many repetitions of this process yield a new $f$ with image in $X_m$, without changing the original $f$ on $S^k$.
\end{proof}

\begin{proof}[Proof of Theorem~\ref{XtoY}]  
We wish to show that $\pi_kY=0$ if $k \le (n-2)/2$, so suppose we are given a map $F\cln S^k \to Y$. By composing with the projection $\pi\cln Y\to X$ we get a map $f\cln S^k\to X$. The first step of the proof will be to construct an `enlargement' of this to a map $g\cln S^k \to X_{k+2}$. By Lemma~\ref{genlargement}, this extends to $g\cln D^{k+1}\to X_{k+2}$ if $k\le n - (k+2)$, which is equivalent to $k\le (n-2)/2$.  Using this $g$ and the coloring lemma we will then construct a map $G\cln D^{k+1}\to Y$ whose restriction to $S^k$ is homotopic to $F$, thus showing that $\pi_kY =0$. 

\medskip

\noindent 
{\em Step 1: Construction of $g\cln S^k\to X_{k+2}$.}
The given map $F\cln S^k\to Y$ may be taken to be simplicial with respect to some triangulation $\T_0$ of $S^k$. We denote by $\T_0'$ the barycentric subdivision of $\T_0$. A $p$-simplex of $\T_0'$ is thus a chain  $[\s_0\!<\!\cdotss\!<\!\s_p]$ of simplices of $\T_0$. We want to construct a simplicial map $g$ from a subdivision $\T_1$ of $\T_0'$ to $X_{k+2}$ with the following additional property: 

\smallskip
\begin{list}{}{\setlength{\leftmargin}{22pt}\setlength{\labelwidth}{16pt}\setlength{\labelsep}{5pt}\setlength{\itemsep}{3pt}}

\item[($*$)]  For any vertex $v$ of $\T_1$ such that $v$ lies in the interior of a simplex $[\s_0\!<\!\cdotss\!<\!\s_p]$ of $\T_0'$, there exists a lift of $g(v)$ to $Y$ containing $F(\s_0)$ as a face.
\end{list}

\smallskip
\noindent
We will construct $g$ inductively over the skeleta of $\T'_0$. 

For each vertex $[\sigma]$ of $\T'_0$, with $f(\s)$ a $p$-simplex of $X$, choose  a $(k-p)$-simplex $\tau\in\pi(\link(F(\sigma)))$ and  let $g([\sigma])=f(\sigma)*\tau$, a vertex of $X_{k+2}$. Such a $\tau$ exists since  $(\pi\link(F(\sigma)))_{k-p+1}$ is $(n-p-2-(k-p+1)=n-k-3)$-connected by Lemma~\ref{genlargement}, and $k\le (n-2)/2$ implies that $n-k\ge 2$ so that this complex is nonempty. Since $\tau\in\pi(\link(F(\sigma)))$, there exists a lift of $g([\s])$ containing $F(\s)$ as a face. Thus the property ($*$) is satisfied for the vertex $[\s]$.
For the inductive step we wish to extend $g$ over a $p$-simplex $\tau=[\sigma_0 \!<\! \cdotss \!<\! \sigma_p]$ of $\T_0'$, assuming we have already defined $g$ on $\del\tau$ so that ($*$) is satisfied. In particular, for any vertex $v$ of $\del\tau$, there exists a lift of $g(v)$ to $Y$ containing $F(\s_0)$ as a face as there exists a lift containing $F(\s_i)\ge F(\s_0)$ for some $i$. So the restriction of $g$ to $\del\tau$ has the form $g_\tau * f(\s_0)$ for $g_{\tau}\cln \del\tau \to (\pi\link(F(\sigma_0)))_{k-q+2}$, where $q$ is the number of vertices in $f(\sigma_0)$. This projected link has connectivity $n-(q-1)-2-(k-q+2)=n-k-3$ by hypothesis. 
We need it to be $(p-1)$-connected in order to extend $g_{\tau}$ over the whole simplex $\tau $, and hence to extend $g$ as $g_\tau * f(\s_0)$. Thus we need the inequality $p-1\le n-k-3$, which holds if $k\le (n-2)/2$ as $p\le k$. This gives the induction step in the construction of $g$ on $S^k$, and the property ($*$) still holds after this step.

\smallskip

The map $g$ extends to $g\cln D^{k+1}\to X_{k+2}$ by Lemma~\ref{genlargement}. We may take this map to be simplicial with respect to a triangulation $\T_1$ of $D^{k+1}$ which extends $\T_1$ on~$S^k$.

\medskip

\noindent
{\em Step 2: Construction of $G\cln D^{k+1}\to Y$.}
Since $\pi$ is surjective we can choose, for each vertex $v$ of $\T_1$, a lift of $g(v)$ to a simplex $\overline{g(v)}$ of $Y$. If $v$ is in the interior of a simplex $[\s_0\!<\!\cdotss\!<\!\s_p]$ of $\T'_0$ in $S^k$, we choose a lift containing $F(\s_0)$ as a face.
This is possible by property ($*$) in Step 1. 

\smallskip

The construction of $G\cln  D^{k+1}\to Y$ will be inductive over the skeleta of $\T_1$. For a vertex $v$ of $\T_1$ we let $G(v)$ be some vertex of $\overline{g(v)}$. In case $v$ is in the interior of a simplex  $[\s_0\!<\!\cdotss\!<\!\s_p]$ of $\T_0'$ in $S^k$, we choose $G(v)$ to be a vertex of $F(\s_0)$. This means in particular that if $v$ is in a simplex $\sigma$ of $\T_0$, then $G(v)=F(w)$ for some vertex $w$ of $\sigma$. As $\T_1$ is a subdivision of $\T_0$, it follows that $G$ extends over $S^k$ linearly on simplices of $\T_1$, and this extension is linearly homotopic to $F$.

\smallskip

We define $G$ on the higher skeleta in the interior of $D^{k+1}$ so that it is simplicial on a subdivision $\T_2$ of $\T_1$ that equals $\T_1$ on $S^k$, and so that the following property is satisfied: 

\smallskip
\begin{list}{}{\setlength{\leftmargin}{22pt}\setlength{\labelwidth}{18pt}\setlength{\labelsep}{5pt}\setlength{\itemsep}{3pt}}

\item[($**$)]  If $w$ is a vertex of $\T_2$ which lies in the interior of a simplex $\tau=\lgl v_0,\cdotss,v_p\rgl$ of $\T_1$ with $g(v_0)\le\cdotss\le g(v_p)$ in $X_{k+2}$, then $\pi G(w)$ is a vertex in $g(v_0)$ and $G(w)$ is the lift of this vertex in $\overline{g(v_p)}$.
\end{list}
\smallskip

\noindent
Note that this is satisfied by the above definition of $G$ on the 0-skeleton of $\T_1$ since $v_p=v_0$ in this case. In particular the condition ($**$) is satisfied on the whole boundary sphere $S^k$, where $\T_2=\T_1$.

Suppose that we have extended $G$ over the $(p-1)$-skeleton of $\T_1$.  Let $\tau$ as in ($**$) be a $p$-simplex of $\T_1$ not contained in $S^k$. For each vertex $w$ of the new triangulation $\T_2$ of $\del\tau$, we have defined $G(w)\in\overline{g(v_j)}$ with $\pi G(w)\in g(v_i)$ for some $i,j$ with $i\le j$. Let $E$ be the set of vertices in $g(v_p)$ and $E_0$ the subset of vertices in $g(v_0)$. We apply the coloring lemma (Lem.~\ref{coloring}) to the sphere $\del\tau$ with vertices labeled by $E$ via $\pi\circ G$. This gives an extension of the triangulation $\T_2$ over $\tau$ with $\del\tau$ a full subcomplex, and with the vertices interior to $\tau$ labeled by $E_0$ and bad simplices only in $\del\tau$. (Bad simplices may occur in $\del\tau$ if a face of $\tau$ is included in $S^k$.) For a vertex $w$ interior to $\tau$ labeled by a vertex $x$ of $g(v_0)$, we define $G(w)$ to be the lift of $x$ in $\overline{g(v_p)}$. 

To see that this definition is valid we need to check that for each simplex $\lgl w_0,\cdotss,w_q\rgl$ of $\T_2$ in $\tau$ the vertices $G(w_r)$ span a simplex of $Y$. We may assume the vertices $G(w_r)$ all have distinct $\pi$-images since if $\pi G(w_r)=\pi G(w_s)$ with $r\ne s$ then $w_r$ and $w_s$ lie in $\del\tau$, and therefore the edge $\lgl w_r,w_s\rgl$ also lies in $\del\tau$ since $\del\tau$ is a full subcomplex of $\tau$ in the $\T_2$ triangulation; hence $G(w_r)=G(w_s)$ since $G$ is a well defined simplicial map on $\del\tau$ by induction, and the fibers of $\pi$ are discrete. All the images $\pi G(w_r)$ are vertices of $g(v_p)$, so by condition (3$'$) following the definition of a join complex it now suffices to check that each vertex $G(w_r)$ is a vertex of some simplex of $Y$ projecting to $\lgl \pi G(w_0),\cdotss,\pi G(w_q)\rgl$. 
We have defined $G(w_r)\in\overline{g(v_j)}$ with $\pi G(w_r)\in g(v_i)$ for some $i,j$ with $i\le j$. We claim that $g(v_j)$ contains
$\lgl \pi G(w_0),\cdotss,\pi G(w_q)\rgl$.  If this is so, then the face of $\overline{g(v_j)}$ lying over $\lgl \pi G(w_0),\cdotss,\pi
G(w_q)\rgl$ provides the desired lift of $\lgl \pi G(w_0),\cdotss,\pi G(w_q)\rgl$ containing $G(w_r)$. To see that $g(v_j)$ contains
$\lgl \pi G(w_0),\cdotss,\pi G(w_q)\rgl$, suppose this fails, so there exists some $w_s$ with $\pi G(w_s)$ not contained in
$g(v_j)$. Then we would have $G(w_s)\in\overline{g(v_m)}$ with $\pi G(w_s)\in g(v_l)$ for some $l,m$ with $l\le m$ and $ j < l$. In
particular, $j\neq p$ and $l\neq 0$ and the edge $\lgl w_r,w_s\rgl$ would have its endpoints contained in two disjoint
faces of $\tau$. Since $\del\tau$ is a full subcomplex of $\tau$ in the $\T_2$ triangulation, this would force $\lgl w_r,w_s\rgl$ to be
contained in $\del\tau$, hence in a $(p-1)$-dimensional face of $\tau$. By the inductive construction, this face would have its
boundary a full subcomplex in the $\T_2$ triangulation, and the same argument forces $\lgl w_r,w_s\rgl$ to lie in a $(p-2)$-dimensional face of $\tau$. Iterating this argument, we eventually reach a contradiction. 
\end{proof}

\section{the complexes}\label{complexes}

In this section, we define the complexes $X^A$ and $X^{FA}$ needed for the proof of Theorem~\ref{main}, and we deduce from
Section~\ref{connectivity} that these complexes are highly connected. The case $P=S^1\x S^2$ differs from the general case and is
treated separately. (Though we will use the same notation, we note that the complexes $X^A$ and $X^{FA}$ will be defined
differently in the two cases $P\neq S^1\x S^2$ and $P=S^1\x S^2$.)  
We end the section with some properties of the action of $\G(M,R)$ on the complexes.

\subsection{Prime summands $\bf P\neq S^1\x S^2$}

Let $M$ be a compact connected oriented 3-manifold and let $\Sph(M)$ denote the simplicial complex whose vertices are isotopy classes of embedded spheres in $M$ which neither bound a ball nor are isotopic to a sphere of $\del M$. A set of vertices of $\Sph(M)$ spans a simplex when the corresponding spheres can be isotoped to be all disjoint. 
It was shown in \cite[Thm~3.1, statement (1)]{HW} that $\Sph(M)$ is contractible if $M$ is not irreducible or the connected sum of an irreducible manifold with copies of $D^3$.

If $P$ is a nontrivial connected summand of $M$, let $P^0$ be the manifold obtained from $P$ by deleting the interior of a ball, and let $\del_0P^0\subset\del P^0$ be the sphere bounding this deleted ball. Consider orientation-preserving embeddings $f\cln P^0\to M$ with $f(\del P^0 - \del_0P^0)\subset\del M - R$ for $R$ a given compact subsurface of $\del M$, possibly empty.  Let $X=X(M,P,R)$ be the simplicial complex whose vertices are isotopy classes of images of such embeddings, where a set of $k+1$ vertices spans a $k$-simplex of $X$ if the vertices are represented by embeddings with disjoint images.

\begin{prop}\label{sphP}
If $P\ne S^1 \x S^2$ 
then $X$ is CM of dimension $n-1$, where $n$ is the number of $P$-summands of $M$ disjoint from $R$. 
\end{prop}

\begin{proof}
We first dispose of the easy special case that $P=D^3$. Then $P^0=S^2\x I$ and vertices of $X$ just correspond to boundary spheres of $M$ disjoint from $R$, hence $X$ is a simplex $\Delta^{n-1}$ which is certainly CM.  Thus we can assume $P\neq D^3$ from now on.

Another easy case is when $n=1$ so $X$ is $0$-dimensional, since the CM condition is then automatic. 

When $n\ge2$ there is a simplicial map $X \to \Sph(M)$ sending a collection of disjoint copies of
$P^0$ in $M$ to their boundary spheres $\del_0P^0$, which are nontrivial since $n\ge 2$. This map is injective except in the special case that $M=P\#P$ and $R=\emp$. In this case $\Sph(M)$ is $0$-dimensional since $P\neq S^1\x S^2$, and it follows that $\Sph(M)$ is a single point since it is contractible. The complex $X$ is therefore a $1$-simplex, which is CM. Excluding this special case from now on, we can view $X$ as a subcomplex of $\Sph(M)$.

Next we show that $X$ is $(n-2)$-connected by an argument that proceeds by induction on the `complexity' of $M$, the
number of spheres defining a maximal simplex of $\Sph(M)$. This is the same number for all maximal simplices.  The case of complexity $1$ is covered by the special cases already considered.  For the induction step, let $f\col S^k\to X$ be 
any map with $k\le n-2$. We can extend $f$ to a map $\hat f\col D^{k+1}\to \Sph(M)$ as the latter complex is
contractible. We can assume that 
$\hat f$ is simplicial with respect to some triangulation of $D^{k+1}$. We are going to modify $\hat f$ in the interior of $D^{k+1}$ so that its image lies in $X$.  Let $\s$ be a maximal simplex of $D^{k+1}$ with the property that none of its vertices map to $X$. Let $M_1, \cdotss, M_d$ be the manifolds obtained by splitting $M$ along the spheres of $\hat f(\s)$. If $\s$ has dimension $p$ then we have $\hat f\col \link(\s) = S^{k-p}\to X(M_1)*\cdotss * X(M_d)$. Each $M_i$ has smaller complexity than $M$. 
Hence by induction $X(M_i)$ is $(n_i-2)$-connected where $n_i$ is the number of $P$-summands in $M_i$ disjoint from $R$. Moreover, $\sum_i n_i=n$ as summands $P\ne S^1\x S^2$ cannot disappear when we cut along spheres. Thus $X(M_1)*\cdotss * X(M_d)$ is $((\sum_i n_i)-2)=(n-2)$-connected and the restriction of $\hat f$ to the link of $\s$ can be extended to a map of the disk $D^{k-p+1}$. We take the join of this map with the restriction of $\hat f$ to $\del\s$ and use this join to modify $\hat f$ in the interior of the star of $\s$. This process reduces the number of maximal simplices such as $\s$. Finitely many repetitions of this step finish the argument that $X$ is $(n-2)$-connected.

The link of a $p$-simplex of $X$ is isomorphic to $X(M')$ where $M'$ is obtained from $M$ by removing the submanifolds $P^0$ corresponding to the $p$-simplex. The manifold $M'$ has $(n-p-1)$ $P$-summands disjoint from $R$, so $X(M')$ is $(n-p-2)$-spherical. Thus $X$ is CM. 
\end{proof}

We assume now that $M$ has nonempty boundary and that $R$ is nonempty, and we choose a basepoint $x_0\in R$. 
For the proof of our main theorem we will need two enhancements of the complex $X$ that include more data than just the submanifolds $f(P^0)$.  Let $I\vee P^0$ be obtained from the disjoint union of $I$ and $P^0$ by identifying $1\in I$ with a basepoint $p_0$ in $\del_0P^0$.
Consider embeddings $f\cln I\vee P^0 \to M$ whose restriction to $P^0$ is orientation-preserving, with $f(P^0)$ giving a vertex of $X$, and such that $f(0)=x_0$.  Let $X^A=X^A(M,P,R,x_0)$ be the simplicial complex whose vertices are isotopy classes of such embeddings $f\cln I\vee P^0\to M$, with simplices corresponding to sets of embeddings with images that are disjoint except at $x_0$.  More generally, if we specify a subsurface $T$ of $\del M -R$ as in Section~\ref{results}, we can define a complex $X^A=X^A(M,P,R,T,x_0)$ whose vertices are isotopy classes of embeddings $f\cln I\vee P^0\to M$ as above, modulo twists along disks in $f(P^0)$ with boundary in $T$.  There is a natural projection $X^A\to X$ induced by sending an embedding $f\cln I\vee P^0 \to M$ to the submanifold $f(P^0)$.

Note that for a collection of embeddings $f_i\cln I\vee P^0 \to M$ defining vertices of $X^A$, it is easy to make the arcs
$f_i(I)$ disjoint except at $x_0$ just by general position, so the condition for these embeddings to span a simplex
reduces to the existence of isotopies making the submanifolds $f_i(P^0)$ disjoint and the arcs $f_i(I)$ disjoint from
these submanifolds (except at their endpoints $f_i(1)$). 
It follows that $X^A$ is a join complex over $X$. In the
notation of Definition~\ref{joincpx}, given a simplex $\s=\lgl x_0,\cdotss,x_p\rgl$ of $X$, the set $(X^A)_{x_i}(\s)$ is the
set of embeddings $f_i$ such that $f_i(P^0)=x_i$ and $f_i(I)$ is disjoint from $x_j$ for each $j\neq i$. 

We note that the equivalence relation of isotopy on the arcs $f_i(I)$ is the same as homotopy, by the lightbulb trick.

As a further refinement of $X^A$, consider pairs $(f,\tau)$ where $f$ is as above and $\tau$ is a framing of the normal
bundle of the arc $f(I)$ which agrees with a fixed normal framing of $x_0$ in $\del M$ at one end and with the image
under $f$ of a fixed normal framing of $p_0$ in $\del_0P^0$ at the other end. We assume these two fixed framings are
chosen according to some orientation convention that allows the framings at the endpoints of $f(I)$ to extend over
$f(I)$. Taking collections of  isotopy classes of such pairs $(f,\tau)$ with the same disjointness conditions as before
yields a complex $X^{FA}$ with a projection to $X^A$ and hence also to $X$. The complex $X^{FA}$ is a join complex over
$X$ and a complete join complex over $X^A$. (For $X^{FA}$ we will not need the generalization involving factoring out twists along disks with boundary in $T$.)

The choice of the framing $\tau$ is {\it a priori\/} parametrized by $\Z=\pi_1SO(2)$, but actually the choice lies in $\Z/2 = \pi_1SO(3)$ since twice a twist along the sphere $f(\del_0P^0)$ is isotopically trivial as a diffeomorphism of $P^0$ fixing $\del_0P^0$. For some prime manifolds $P$ such as lens spaces the twist along $\del_0P^0$ is itself trivial, so in these cases the choice of framing is unique.

\begin{prop}\label{XA}
If $P\ne S^1 \x S^2$ then 
$X^A$ and $X^{FA}$ are $(\frac{n-3}{2})$-connected, where $n$ is the number of $P$-summands in $M$ disjoint from $R$. 
\end{prop}

\begin{proof}
Consider $X^A$ first. The projection $X^A\to X$ expresses $A$ as a join complex over $X$. Now we check that the hypotheses of Theorem~\ref{XtoY} are satisfied, with the `$n$' there replaced by $n-1$. By Proposition~\ref{sphP}, $X$ is CM of dimension $n-1$.  Also, the projection of the link of a $p$-simplex $\s$ of $X^A$ is an $(n-p-2)$-dimensional CM subcomplex of $X$. This follows from the fact that the projection of the link of $\s$ is isomorphic to $X(M_\s)$ where $M_\s$ is the submanifold of $M$ obtained by deleting a neighborhood of the union of the images $f_i(I\vee P^0)$ for the vertices of $\s$. Thus the number of $P$-summands in $M_\s$ disjoint from $R$ is $n-p-1$ and hence $X(M_\s)$ is a CM complex of dimension $n-p-2$.

The complex  $X^{FA}$ is treated in exactly the same way.
\end{proof}

In the special case $P=D^3$ the projections $X^{FA}\to X^A$ is an isomorphism, and we can view simplices of $X^A$ as isotopy classes of systems of arcs from $x_0$ to boundary spheres of $M$.  The connectivity result in this case can be improved:

\begin{prop}\label{D3}
When $P=D^3$ the complex $X^A$ is $(n-2)$-connected, where $n$ is the number of boundary spheres of $M$ disjoint from $R$ 
\end{prop}

\begin{proof}
In this case $X^A$ is a complete join complex over $X=\Delta^{n-1}$ so the result follows from Proposition~\ref{label}.
\end{proof}

\begin{rem}\label{pi1}{\rm
 The lightbulb trick is applicable in this context, so isotopy classes of systems of arcs in $M$ defining simplices of
 $X^A$ are the same as homotopy classes. When $P=D^3$, vertices of $X^A$ over a vertex $x$ of $X$ correspond to elements
 of $\pi_1M$ once one lift of $x$ to $X^A$ is chosen.
 Thus the complex $X^A$ in this case is isomorphic to the join of $n$ copies of $\pi_1M$.  
 } \end{rem}

\subsection{Prime summands $\bf P=S^1\x S^2$}

When $P=S^1\x S^2$ we replace the earlier complex $X$ by $\Sph_c(M)$, the complex of nonseparating sphere systems in $M$. It was shown in \cite[Prop.~3.2]{HW} that this complex is $(n-2)$-connected, where $n$ is the number of $S^1\x S^2$ summands in $M$. 
It follows easily that $\Sph_c(M)$ is Cohen-Macaulay of dimension $n-1$. 

We will also change to a new complex $X^A$. To define this, let $I + S^2$ denote the quotient space of the disjoint union of $I$ and $S^2$ obtained by identifying the midpoint of $I$ with a basepoint in $S^2$. Consider embeddings $f\cln I+S^2 \to M$ which are smooth on $I$ and $S^2$, with $f(I)$ and $f(S^2)$ intersecting transversely, such that $f(0)=x_0$ and $f(1)=x_1$ for chosen points $x_0,x_1 \in R$, and where $f(S^2)$ is a nonseparating sphere in $M$. The possibility $x_0=x_1$ is allowed. We also impose an orientation condition on $f$:  The standard orientations of $I$ and $S^2$ give orientations of $f(I)$ and $f(S^2)$ and we require that at $f(I)\cap f(S^2)$ these orientations combine to give the orientation of $M$.  Let $X^A=X^A(M,x_0,x_1)$ be the complex whose vertices are the isotopy classes of such embeddings $f\cln I+S^2\to M$. Simplices of $X^A$ are given by collections of such embeddings $f_i$ whose images are all disjoint except at $x_0$ and $x_1$, and such that the spheres $f_i(S^2)$ form a nonseparating system. Thus there is a natural projection $X^A\to \Sph_c(M)$.  The lightbulb trick again applies, so we can vary the two halves of the arcs $f_i(I)$ on either side of $f_i(S^2)$ by homotopy as well as isotopy.

We can refine $X^A$ to a complex $X^{FA}$ by taking isotopy classes of pairs $(f,\tau)$ where $f$ is as above and $\tau$ is a framing of $f(I)$ that agrees with fixed framings at $x_0$ and $x_1$.

\begin{prop}\label{S1xS2}
When $P=S^1\x S^2$ the complexes $X^A$ and $X^{FA}$ are $(\frac{n-3}{2})$-connected, where $n$ is the number of $S^1\x S^2$ summands in $M$. 
\end{prop}

\begin{proof} The argument for Proposition \ref{XA} applies here as well.
\end{proof}

\subsection{Action of the mapping class group}

The group $\Dif(M,R)$ of diffeomorphisms of $M$ that restrict to the identity on $R$ acts on the complexes $X^A$ and $X^{FA}$ in the preceding three propositions by applying diffeomorphisms to the various types of objects defining simplices and to isotopies of these objects. Elements of $\Dif(M,R)$ in the path-component of the identity act trivially, so there are induced actions of $\G(M,R)$. 
Twists along spheres act trivially on the complexes $X^A$: recall from Section~\ref{mcg} that such a twist is isotopic
to a composition of twists along spheres disjoint from any given sphere system. So one just has to consider how twists
act on arcs, where they preserve homotopy classes and hence isotopy classes.  Similar reasoning shows that twists along disks with boundary in $T$ act trivially on $X^A$. Hence there are induced actions of $A(M,R,T)$ on $X^A$.

\medskip\noindent
Here are three properties of these actions that will be needed in the proof of the main theorem:

\begin{lem}\label{action} 
For each of the cases $P\neq S^1\x S^2$ and $P=S^1\x S^2$ the actions of $\G(M,R)$ on $X^{FA}$ and $A(M,R,T)$ on $X^A$ satisfy\/{\rm :}

\begin{enumerate}[{\rm (i)}]
\item The action is transitive on simplices of a given dimension.

\item  The stabilizer of a simplex contains elements realizing arbitrary permutations of the vertices of the simplex, using diffeomorphisms supported in a neighborhood of the object representing the simplex.

\item An element of the pointwise stabilizer of a simplex can be represented by a diffeomorphism 
supported outside a neighborhood of an object representing the simplex.
\end{enumerate}
\end{lem}

\medskip

In each case there is a natural map from $\G(M_\s,R_\s)$ or $A(M_\s,R_\s,T)$ to the pointwise stabilizer of a simplex $\s$, where $M_\s$ is the complement of a neighborhood $N(\s)$ of an object representing $\s$ and $R_\s=(R\cap M_\s)\cup \del N(\s)$. 
By (iii) this map is surjective, and from Proposition~\ref{stabilizationmap} it  follows that it is also injective.

\begin{proof}
Let us start with statement (i) in the case $P\neq S^1\x S^2$.  Fix a maximal simplex $\lgl[f_0],\cdotss,[f_m]\rgl$ of
$X^A$, where each $f_i$ is an embedding of $I\vee P^0$.  Given any simplex $\lgl[g_0],\cdotss,[g_k]\rgl$ of $X^A$, the
proof of the uniqueness of the factorization of an oriented $3$-manifold into prime manifolds gives an
orientation-preserving diffeomorphism $h$ of $M$ taking each $g_i(P^0)$ to some $f_j(P^0)$, and $h$ can be chosen to be
fixed on $\del M$ (see for example \cite{Hem}). Any permutation of the submanifolds $f_i(P^0)$ can be realized by a
diffeomorphism of $M$ fixing $R$ as $R$ does not intersect these submanifolds, so we can modify $h$ keeping $R$ fixed so
that $h$ takes $g_i(P^0)$ to $f_i(P^0)$ for $0\le i\le k$.  Another modification of $h$ allows us to assume that
$g_i=f_i$ on $P^0$ for each $i$ by composing $h$ with a diffeomorphism supported in a neighborhood of $g_i(P^0)$, using
the fact that $g_i$ and $f_i$ preserve orientation on $P_0$. We can further arrange that $h$ takes $g_i(I)$ to $f_i(I)$
for each $i$ by the balloon trick explained after the proof of Proposition~\ref{stabilizationmap}. 
For the action on $X^{FA}$, we can also arrange that $h$ preserves the framing of the arcs by changing the framings if
necessary by composing $h$ with twists along the spheres $f_i(\del_0P^0)$. 
This proves (i) in the case $P\neq S^1\x S^2$. 
When $P=S^1\x S^2$ one can proceed in a similar fashion, but the construction is more elementary since it suffices just
to cut $M$ along the spheres $f_i(S^2)$ or $g_i(S^2)$, capping the resulting boundary spheres off with balls, then slide
these balls around in the resulting manifold to obtain the required diffeomorphisms. 

Using (i) it suffices to prove (ii) for any particular simplex, so one can choose a nice model where the existence of the desired diffeomorphism permuting the vertices is obvious.

For (iii), consider first the case of $X^{FA}$ with $P\neq S^1\x S^2$. Let a $k$-simplex $\s$ of $X^{FA}$ be represented by a $(k+1)$-tuple of framed embeddings $(f_i,\tau_i)$ of $I\vee P^0$. If $g\in\Dif(M,R)$ represents an element of the pointwise stabilizer of $\s$ then there is an isotopy of the restriction  $g\!\bigm|\!\cup_if_i(I\vee P^0)$ to the identity.  Let us call this isotopy $G$. We would like to extend $G$ to an isotopy of $g$ itself, staying fixed on $R$. If $G$ were stationary on a half-ball neighborhood $B$ of $x_0$ in $M$ then the standard isotopy extension technique for extending isotopies of submanifolds could be applied, first on the submanifolds $f_i(P^0)$ and then on the arcs, to give an isotopy of $g$, fixing $B$ and $R$.  Let us show how to modify $G$ to be stationary on $B$.
We can assume each $f_i(I\vee P^0)$ intersects $B$ in a radial line segment and the isotopy $G$ deforms these segments
through radial segments. The outer endpoints of these segments trace out a pure braid on the boundary hemisphere of
$B$. The isotopy $G$ can be modified to an isotopy $G'$ which is constant in $B$ and ends with a new collection of
embeddings $f'_i$ that differ from $f_i$ only by a pure braid just outside $B$. This collection of $f'_i$'s is isotopic
to the $f_i$'s by an isotopy $H$ fixed on $B$ using the lightbulb trick, undoing the braid by lifting the strands over
the spheres $f_i(\del_0 P^0)$. The composite isotopy consisting of first $G'$, then $H$, is a new isotopy of
$g\!\bigm|\!\cup_if_i(I\vee P^0)$ to the identity that is fixed inside $B$. By our earlier remarks we can then isotope
$g$ to be the identity on $\cup_if_i(I\vee P^0)$. The framings $g(\tau_i)$ can be carried along during this isotopy of
$g$, yielding framings that differ from the $\tau_i$'s by even numbers of full twists. These twists can be eliminated by composing with even powers of twists along the spheres $f_i(\del_0 P^0)$, as these even powers are isotopic to the identity. Then we can further isotope $g$ to be the identity in a neighborhood of $\cup_if_i(I\vee P^0)$.

In the case of $X^A$ with $P\neq S^1\x S^2$, the framing data is absent, so after isotoping $g$ to a diffeomorphism $g'$ that fixes $\cup_if_i(I\vee P^0)$ it might be necessary to compose $g'$ with twists along some of the spheres $f_i(\del_0 P^0)$ to make it preserve framings, in order to isotope it further to be the identity in a neighborhood of $\cup_if_i(I\vee P^0)$. Composing with these twists does not change the class of the diffeomorphism in $A(M,R,T)$. 

The argument for (iii) when $P=S^1\x S^2$ is similar.
\end{proof}

\section{spectral sequence arguments i}\label{spectralseq}

In this section, we prove the first part of Theorem~\ref{main} and Corollary~\ref{free}. This requires a spectral sequence argument, which we state in enough generality so that it can be applied directly to prove other stability results in later sections. This argument has been used many places in the literature to prove stability theorems, but we don't know of a reference for a general statement that applies to all the cases we need. 

\medskip

Suppose we have a sequence of groups $G_1\subset G_2\subset\cdotss$ with
actions of each $G_n$ on a complex $X_n$ of dimension at least $n-1$ 
 such that:

\begin{enumerate}[{(1)}]
\item  The action is transitive on simplices of each dimension, and either 
\begin{enumerate}[{(a)}]
\item the stabilizer of each simplex contains elements that give all permutations of its vertices, or 
\item the stabilizer of each simplex fixes the simplex pointwise.
\end{enumerate}

\item  The subgroup of $G_n$ fixing a $p$-simplex pointwise is conjugate to $G_{n-k-1}$ for some $0\le k\le p$ (where $G_i$ is trivial if $i<1$).

\item  For each edge of $X_n$ with vertices $v$ and $w$ there exists an element of $G_n$ that takes $v$ to $w$ and commutes with all elements of $G_n$ that leave the edge fixed pointwise. 
\end{enumerate}

\noindent
An alternative to condition (1) is  
\begin{enumerate}[{(1$'$)}]
\item The action is transitive on vertices, the stabilizer of each simplex fixes the simplex pointwise, and 
$H_i(X_n/G_n)=0$ for $1\le i\le n-2$. 
\end{enumerate}

\begin{thm}\label{range}
Under these conditions, 
if each $X_n$ is 
 $(n-2)/2$-connected then the inclusion $G_n\to G_{n+1}$ induces an isomorphism on $H_i$ if $n\ge 2i+1$ and a surjection if $n=2i$.
\end{thm}

Note that the theorem also applies if $X_n$ is $(n-2)$-connected since $n \ge 1$. With a shift in indices the theorem takes the following form, which is what we will use most often:

\begin{cor}\label{range2}
Under the same conditions, if $X_n$ is $(n-3)/2$-connected for each $n\ge 2$ then the inclusion $G_n\to G_{n+1}$ induces an isomorphism on $H_i$ if $n\ge 2i+2$ and a surjection if $n=2i+1$.
\end{cor}

\begin{proof} Apply the theorem to the sequence of groups $G'_n = G_{n+1}$ acting on the spaces $X'_n =X_{n+1}$. The hypotheses are easily checked, and the conclusion is that $G_{n+1} \to G_{n+2}$ is an isomorphism on $H_i$ for $n+1\ge 2i+2$ and a surjection for $n+1 = 2i+1$, which is equivalent to what the Corollary claims.
\end{proof}

There is also a stable version of the above theorem:  Consider a group $G_\infty$ with a self inclusion $\lambda\cln G_\infty \inc G_\infty$. Suppose that $G_\infty$  acts on a contractible infinite dimensional complex $X_\infty$, satisfying (1-3) above (or (1$'$-3)) with $n=\infty$ and $\lambda$ replacing the inclusion $G_n\inc G_{n+1}$ for each $n$. In particular, condition (2) requires now that the inclusion of the stabilizer of a $p$-simplex $\St(\s_p)\inc G_\infty$ is conjugate to $\lambda^{k+1}$.

\begin{thm}\label{stable-stability}
Under these conditions, the map $\lambda$ induces a homology isomorphism $\lambda_*\cln H_*(G_\infty)\to H_*(G_\infty)$. 
\end{thm}

\begin{proof}[Proof of Theorems~\ref{range} and \ref{stable-stability}]
(For Theorem~\ref{stable-stability}, take $n=\infty$ in the proof --- this actually simplifies many of the statements.) 

We first reduce (1.a) and (1.b) to (1$'$). 
For (1.a), we replace $X_{n}$ by $Y_{n}=\Delta X_{n}$, the simplicial set (or $\Delta$-complex) whose $p$-simplices are the simplicial maps from the standard $p$-simplex $\Delta^p$ to $X_{n}$.  This simplicial set is homotopy equivalent to $X_{n}$.  The action of $G_{n}$ on $X_{n}$ induces an action on $Y_{n}$ which is no longer transitive on the set of $p$-simplices for $p\ge 1$, but is such that the stabilizer of a simplex fixes the simplex pointwise. Indeed, the stabilizer of a simplex $\s$ of $Y_{n}$ is the subgroup of $G_{n}$ fixing its image $\operatorname{Im}(\s)$ pointwise, and this subgroup fixes $\s$ pointwise.  By assumption, the stabilizer $\St(\s)$  is conjugate to $G_{n-k-1}$ for some $k\le q$ if $q+1$ is the number of distinct vertices of $X_{n}$ in the image of $\s$. If $\s$ is a $p$-simplex, we have $0\le q\le p$, so that condition (2) is still satisfied when $X_n$ is replaced by $Y_n$. Note also that there are two orbits of edges in $Y_n$, one degenerate and one nondegenerate and condition (3) holds trivially on the degenerate orbit, and by assumption on $X_n$ for the other orbit. The quotient $\Delta X_n/G_n$ is $(n-2)$-connected by \cite[Lemma~3.5]{H3}, using the assumption that $X_n$ has dimension at least $n-1$. Hence condition (1.a) for $X_n$ is replaced by condition (1$'$) for $Y_n$, while preserving conditions (2) and (3). 

Condition (1.b) is a special case of (1$'$) since $X_n/G_n$ in that case is a $\Delta$-complex with a single simplex in each dimension $k\le n-1$. This complex 
has the homology of a point below dimension $n-1$, as the boundary maps in its simplicial chain complex are alternately the identity and the $0$-map. 

\smallskip

It remains to prove the theorem assuming (1$'$),(2) and (3). 
 We consider the double complex $E_*G_{n+1}\otimes_{G_{n+1}}\widetilde{C}_*(X_{n+1})$, where $E_*G_{n+1}$ is a free resolution of $\Z$ over $\Z G_{n+1}$ and $\widetilde{C}_*(X_{n+1})$ is the augmented chain complex of $X_{n+1}$. This gives two spectral sequences, one of which has $E^\infty_{p,q}=0$ for $p\le \frac{n-1}{2}$ due to the connectivity of $X_{n+1}$. Thus the other spectral sequence has $E^\infty_{p,q}=0$ for $p+q\le \frac{n-1}{2}$. It has $E^1$-term given by 
$$E^1_{p,q}=\bigoplus_{\rm orbits} H_q(\St(\s_p),\Z)$$
where the sum runs over representatives $\s_p$ of the orbits of $p$-simplices. The coefficients are not twisted because the stabilizer of a simplex of $X$ fixes the simplex pointwise. (This uses Shapiro's lemma. See \cite[VII.7]{B} for more details.)
The differentials are induced by the alternating sum of the face maps:
$$H_q(\St(\s_p),\Z)\stackrel{d_i}{\rar}H_q(\St(d_i\s_p),\Z)\stackrel{c_h}{\rar}H_q(\St(\s_{p-1}),\Z)$$ 
where $c_h$ is conjugation by an element $h\in G_{n+1}$ which takes $d_i\s_p$ to the representative $\s_{p-1}$ of its orbit. We want to show that the differential 
$$d^1\cln E^1_{0,i}=H_i(G_{n})\rar E^1_{-1,i}=H_i(G_{n+1})$$ 
is surjective when $n\ge 2i$ and injective when $n\ge 2i+1$.

We prove this by induction on $i$, the case $i=0$ being trivial. We start with surjectivity,  
so assume that $n\ge 2i$. 
Surjectivity of the $d^1$ above follows from: 
\begin{list}{}{\setlength{\leftmargin}{32pt}\setlength{\labelwidth}{16pt}\setlength{\labelsep}{5pt}\setlength{\itemsep}{3pt}}
\item[(1)] $E^\infty_{-1,i}=0$;
\item[(2)] $E^2_{p,q}=0$ for $p+q=i$ with $q<i$. 
\end{list}
Condition (1) is verified as $E^\infty_{p,q}=0$ when $p+q\le \frac{n-1}{2}$ and $i-1\le \frac{n-1}{2}$ when $n\ge 2i$.

For condition (2), we first show that for $q<i$, the inclusion of stabilizers induces an isomorphism 
$$E^1_{p,q}=\oplus_{\rm orbits} H_q(\St(\s_p),\Z)\stackrel{\cong}{\rar} \oplus_{\rm orbits} H_q(G_{n+1},\Z)$$ 
when $p+q\le i$, and a surjection when $p+q=i+1$. Indeed, for a $p$-simplex $\s_p$, $\St(\s_p)$ is conjugate to $G_{n-k}$ for some $0\le k\le p$ and by the induction assumption, the inclusion $G_{n-k}\inc G_{n+1}$ induces an isomorphism in $H_q$ if  
$n-k\ge n-p\ge 2q+1$,  i.e. if $p+2q\le n-1$. Now $p+2q\le 2i-p\le n-p$ by assumption. If $p+q=i$, we have $p\ge 1$ and the inequality is satisfied, and if $p+q<i$, then we have strict inequalities in the above reasoning, also leading to the desired inequality. The surjectivity when $p+q=i+1$, $q<i$ is checked in the same way. (In that case, we need $p+2q\le n$. Now $p+2q=2i+2-p\le n+2-p\le n$ as $p\ge 2$.) 

The diagram
$$\xymatrix{H_q(\St(\s_p),\Z) \ar[d]\ar[r]^-{d_i}& H_q(\St(d_i\s_{p}),\Z)\ar[dl] \ar[r]^-{c_h} & H_q(\St(\s_{p-1}),\Z) \ar[d]\\
H_q(G_{n+1},\Z) \ar[rr]^-{\rm id}& & H_q(G_{n+1},\Z)
}$$
commutes because $c_h$ acts as the identity on $ H_q(G_{n+1},\Z)$. Thus we have a chain map from the chain complex in the $q$th line of the $E^1$-term to the augmented chain complex of $X_{n+1}/G_{n+1}$ with constant coefficients $H_q(G_{n+1},\Z)$, and this map is an isomorphism for $p+q\le i$ and a surjection for $p+q=i+1$ by the previous paragraph. The homology of $X_{n+1}/G_{n+1}$ is trivial in degree $*\le n-1$ by assumption, which 
 implies condition (2) since $i\le n-1$ as $n\ge 2i$ and $i\ge 1$.

To prove injectivity of the map $d^1\cln E^1_{0,i}=H_i(G_{n})\to E^1_{-1,i}=H_i(G_{n+1})$, we will show that when $n\ge 2i+1$
\begin{list}{}{\setlength{\leftmargin}{32pt}\setlength{\labelwidth}{16pt}\setlength{\labelsep}{5pt}\setlength{\itemsep}{3pt}}
\item[(1)] $E^\infty_{0,i}=0$; 
\item[(2)] $E^2_{p,q}=0$ for $p+q=i+1$ with $q<i$; 
\item[(3)] $d^1\cln E^1_{1,i}\to E^1_{0,i}$ is the $0$\,-map. 
\end{list}
\noindent
Conditions (1) and (2) follow from the same argument as above: for (1), we need $i\le \frac{n-1}{2}$, which is equivalent to $n\ge 2i+1$. For (2), we now need  $n-p\ge 2q+1$  when $p+q\le i+1$ and $n-p\ge 2q$ when $p+q=i+2$ for all $q<i$. This is satisfied in the first case as $2q+p\le 2i+2-p\le n+1-p$, which is smaller than $n-1$ if $p+q=i+1$ as $p\ge 2$ in that case, and if $p+q<i+1$ because we get strict inequalities earlier. In the second case as $2q+p=2i+4-p\le n+3-p\le n$ as $p\ge 3$.

For condition (3), 
the  boundary map is $d^1=d^1_1-d^1_0$ on each orbit $\s_1$ with
$d^1_i=c_{h_i}\circ d_i$ for some $h_0,h_1\in G_{n+1}$.  One can choose $h_0$ to be the identity and $h_1$ to be the hypothesized element taking one vertex of $\s_1$ to the other and commuting with every element in $\St(\s_1)$.  On the group level, we have a commutative diagram
$$\xymatrix{\St(\s_1) \ar[r]^-{d_1} \ar[d]_{c_{h_1}=\textrm{id}} & h_1\St(\s_0)h_1^{-1}\ar[d]^{c_{h_1}} \\
\St(\s_1)\ar[r]^{d_0} & \St(\s_0) 
}$$
where the horizontal maps are the inclusions of $\St(\s_1)$ to the stabilizer of its two vertices. 
Hence $d_1^1=d_0^1$ and the $d^1$-differential is $0$.
\end{proof}

\begin{proof}[Proof of Theorem~\ref{main}(i) {\rm (}and Proposition~\ref{D3} for $m=3${\rm\/)}] The result for $P$ nonprime follows from the result for each of its prime summands, so we can assume that $P$ is prime.  When $P$ is neither $S^1\x S^2$ nor $D^3$ we use the actions of $\G^P_n(N,R)$ on $X^{FA}$ and $A^P_n(N,R,T)$ on $X^A$, as described near the end of Section~\ref{complexes}. The conditions (1.a), (2), and (3) at the beginning of this section are satisfied by Lemma~\ref{action} and Proposition~\ref{stabilizationmap}.  Corollary~\ref{range2} then gives the result, using Proposition~\ref{XA} to guarantee the connectivity hypothesis.  For $P=S^1\x S^2$ the argument is similar using instead the complexes $X^A$ and $X^{FA}$ in Proposition~\ref{S1xS2}. For $P=D^3$ one uses the complex $X^A$ in Proposition~\ref{D3}. (Here the same complex works for both groups because the twists along the boundary spheres added by the stabilization are trivial.) In this case the improved dimension range given in Proposition~\ref{higherm} is obtained by applying Theorem~\ref{range} directly rather than Corollary~\ref{range2}. 
\end{proof}

To modify these arguments to apply to the subgroup $\SHAut(*_n G)$ of $\SAut(*_n G)$ appearing in Corollary~\ref{free} we first need a few algebraic preliminaries. Let $G_1,\cdotss,G_n$ be copies of the group $G$, so each $G_i$ is in particular canonically isomorphic to $G$. Inside the free product $F$ of these groups $G_i$ we have the collection $C(G)$ of subgroups $gG_ig^{-1}$ conjugate to the factors $G_i$. Each of these subgroups is isomorphic to $G$, and if $G$ is neither $\Z$ nor itself a free product, these are the only subgroups of $*_iG_i$ with this property, by the Kurosh Subgroup Theorem.  We call an isomorphism $G\to gG_ig^{-1}$ a {\it parametrization\/} of $gG_ig^{-1}$. A subgroup $gG_ig^{-1}$ has a canonical parametrization given by composing the canonical isomorphism $G \to G_i$ with conjugation by $g$, but this canonical parametrization is only well-defined up to inner automorphism of $G$ since the element $g$ can be multiplied on the right by any element of $G_i$.  For a subgroup $H$ of $\Aut(G)$ containing the inner automorphisms, define an $H$-parametrization of a subgroup in $C(G)$ to be a parametrization that differs from a canonical parametrization by an element of $H$.

The symmetric automorphism group $\SAut(F)$ permutes the subgroups in $C(G)$ and acts on parametrizations of these subgroups by composition. Define $\SHAut(F)$ to be the subgroup of $\SAut(F)$ taking $H$-parametrizations to $H$-parametrizations, so $\SHAut(F)$ acts on the set of $H$-parametrizations. This definition is equivalent to the one given in Section~\ref{results}.

\begin{proof}[Proof of Corollary~\ref{free}] 
Assume first that $G\ne \Z$, hence $P\ne S^1 \x S^2, S^1\x D^2$. We can also assume we are not in the trivial case
$P=D^3$.  Since we assume $G$ is not a free product and $P$ is not $S^1\x D^2$, it follows that $\del P$ must be
incompressible in $P$ if it is nonempty.  Let $M$ be the connected sum of $N=D^3$ with $n$ copies of $P$. We take
$R=\del D^3$ with the basepoint $x_0$ in $R$ as usual, and we may as well choose $T=\emp$ since $\del M$ is
incompressible. 
By Proposition~\ref{kernel}, $A(M,R)\subset\Aut(\pi_1M)$ and we have $\SHAut(\pi_1M)\subset A(M,R)$ by restricting the
proof of Proposition~\ref{image}. 
 Consider the action of $A(M,R)$ on $X^A$.  The vertices of $X^A$ corresponding to embeddings $f\cln I\vee P^0\to M$ for which the induced map on $\pi_1$ gives an $H$-parametrization of the subgroup $f_*(\pi_1P^0)$ span a subcomplex $X^A_H$ of $X^A$ which is again a join complex over the same $X$.  This is $(\frac{n-3}{2})$-connected by the same argument as for $X^A$. The group $\SHAut(\pi_1M)\subset A(M,R)$ acts on $X_H^A$, and the hypotheses of Corollary~\ref{range2} are satisfied as before, so Corollary~\ref{free} follows when $G\ne\Z$.  If $G=\Z$ the argument is the same except that we take $P=S^1\x D^2$ and use the action of $\SHAut(\pi_1M)\subset \SAut(\pi_1M)=A(M,R,T)$ on $X_H^A$ with $T$ the boundary tori of~$M$. 
\end{proof}

\section{spectral sequence arguments ii}\label{correction}

In this section we prove part~(ii) of Theorem~\ref{main}, given here as Theorems~\ref{alpha} and \ref{delta} and Corollaries~\ref{mu} and \ref{epsilon}.  
This gives a new argument for the stability theorems of~\cite{HW}, improving the dimension ranges and avoiding a gap in the proofs given there (see \cite{HWerr}).

\smallskip

Let $N$ be a  compact, connected, oriented 3-manifold and $R$ a (possibly empty) submanifold of its boundary. 
Following the notation of \cite{HW} we write $M=M^s_{n,k}$ where 
$$M_{n,k}^s=N\ \#\ (\#_n\ S^1\x S^2)\ \#\ (\#_k\ S^1\x D^2)\ \#\ (\#_s\ D^3).$$ 
Also, we let $\Ga^s_{n,k}=\Ga(M_{n,k}^s,Q)$ where  $Q=R\cup (\del M_{n,k}^s-\del N)=R\cup_k (S^1\x S^1)\cup_s S^2$, and we denote by $A^s_{n,k}$ its quotient $A(M_{n,k}^s,Q,T)$ by the subgroup generated by twists along $2$-spheres and disks with boundary in $T\subset \del N-R$. (In \cite{HW} we actually only consider the case $R=\del N$ and $T=\emp$, but there is no reason to restrict to that case.)
The three principal stabilizations for $\Ga^s_{n,k}$ are:

\begin{list}{}{\setlength{\leftmargin}{30pt}\setlength{\labelwidth}{16pt}\setlength{\labelsep}{5pt}\setlength{\itemsep}{3pt}}
\item[(a)]  $\alpha\cln  \Ga_{n,k}^s \to \Ga_{n+1,k}^s$ induced by enlarging $M^s_{n,k}$ to $M^s_{n+1,k}$ by adjoining a punctured $S^1\x S^2$, identifying a disk in its boundary sphere with a disk in a boundary sphere of $M^s_{n,k}$.
\item[(b)] $\mu\cln   \Ga_{n,k}^s \to \Ga_{n,k}^{s+1}$ induced by enlarging $M^s_{n,k}$ to $M^{s+1}_{n,k}$ by adjoining a three-times punctured sphere, identifying one of its boundary spheres with a boundary sphere of $M^s_{n,k}$;
\item[(c)] $\epsilon\cln \Ga_{n,k}^s\to \Ga_{n,k+1}^s$ induced by enlarging $M^s_{n,k}$ to $M^s_{n,k+1}$ by adjoining a punctured $S^1\x D^2$, identifying a disk in its boundary sphere with a disk in a boundary sphere of $M^s_{n,k}$.
\end{list}
These induce corresponding stabilization maps for $A_{n,k}^s$.
All these stabilizations require $s\ge 1$. To cover the case $s=0$ we need the map $\delta\cln \Ga^{s+1}_{n,k}\to \Ga^s_{n,k}$ induced by filling in a boundary sphere with a ball.  The composition $\delta\mu$ is the identity, so except when $s=0$ the map $\delta$ will induce an isomorphism on homology if $\mu$ does.

For the proof that the above maps induce isomorphisms in a stable range, we consider also the map $\eta\cln  \Ga_{n,k}^s \to \Ga_{n+1,k}^{s-1}$ induced by gluing a three-times punctured sphere along two of its boundary spheres and the map $\beta=\mu\eta\cln  \Ga_{n,k}^s \to \Ga_{n+1,k}^s$. Note that $\alpha=\eta\mu$.

\begin{thm}\label{alpha}
The maps $\alpha_*,\beta_*\cln H_i(\Ga_{n,k}^s) \to H_i(\Ga_{n+1,k}^s)$ and $\alpha_*,\beta_*\cln H_i(A_{n,k}^s) \to H_i(A_{n+1,k}^s)$ are isomorphisms when $n\ge 2i+2$ and surjections when $n\ge 2i+1$ {\rm (}with $s\ge 1$ for $\alpha$ and $s\ge 2$ for $\beta${\rm \/)}. 
\end{thm}

\begin{proof}
This follows from Corollary~\ref{range2}
using the complexes $X^A$ and $X^{FA}$ for $P=S^1\x S^2$ from Proposition~\ref{S1xS2}, with the arcs having both boundary points in a boundary sphere $\del_0M$ for $\alpha$ and one in $\del_0M$, one in another boundary sphere $\del_1M$ for $\beta$.
\end{proof}

The map $\alpha$ is actually the case $P=S^1\x S^2$ in Theorem~\ref{main}(i).
As $\beta=\mu\eta$, and $\mu$ is always injective, we deduce: 

\begin{cor}\label{mu}
For $s\ge 1$, the maps $ \mu_*\cln H_i(\Ga_{n,k}^s) \to H_i(\Ga_{n,k}^{s+1})$ and  $ \mu_*\cln H_i(A_{n,k}^s) \to H_i(A_{n,k}^{s+1})$  are isomorphisms when $n\ge 2i+2$. 
\end{cor}

\begin{thm}\label{delta}
The maps $ \delta_*\cln H_i(\Ga_{n,k}^1) \to H_i(\Ga_{n,k}^0)$ and $ \delta_*\cln H_i(A_{n,k}^1) \to H_i(A_{n,k}^0)$ are isomorphisms when $n\ge 2i+4$. 
\end{thm}

In the case of the groups $A_{n,k}^s$ it is possible to give a proof of this theorem that is very similar to the earlier proofs in the paper.  One uses the complex $X^A$ whose vertices are isotopy classes of pairs $(S,a)$ where $S$ is a sphere in $M_{n,k}^0$ and $a$ is an oriented circle in $M_{n,k}^0$ that intersects $S$ in one point transversely, and simplices are represented by disjoint collections of such pairs $(S_i,a_i)$. (The existence of the dual circles $a_i$ implies that the spheres $S_i$ form a nonseparating system.) The complex $X^A$ is $(n-3)/2$-connected since the circles $a_i$ can be regarded as labeling the spheres $S_i$. The group $A_{n,k}^0$ acts on $X^A$, and the pointwise stabilizer of a $p$-simplex is $A_{n-p-1,k}^{p+1}$. The conditions (1) and (2) at the beginning of Section~\ref{spectralseq} are satisfied but not (3). One cannot interchange two pairs $(S_i,a_i)$ representing an edge of $X^A$ by a diffeomorphism of $M_{n,k}^0$ supported in a neighborhood of the union of the two pairs, but one can do this in a neighborhood of the union of the two pairs with an arc joining their two circles $a_i$. This neighborhood is diffeomorphic to $S^1\times S^2\ \#\  S^1\times S^2$ with a ball removed. The complement of this enlarged neighborhood is $M_{n-2,k}^1$, and this is included in the complement $M_{n-2,k}^2$ of the neighborhood of the two pairs $(S_i,a_i)$. By Corollary~\ref{mu} the map $H_i(A_{n-2,k}^1)\to H_i(A_{n-2,k}^2)$ is surjective under the given restrictions on $n$ and $i$, and this suffices to replace (3) at the point where (3) is used in the proof of Theorem~\ref{range}. (For more details of this type of argument see the proof of the last theorem in \cite{HV3}.) 

One could try to apply this argument for the $\Ga_{n,k}^s$ groups using a version of $X^A$ that includes mod $2$ framing data on the circles $a_i$. However, the stabilizers for the action of $\Ga_{n,k}^s$ on this complex are not quite what one would want them to be, so we will instead use a different sort of argument that works for both the $\Ga_{n,k}^s$ and $A_{n,k}^s$ groups. 

\begin{proof}
We follow here the argument for the proof of Theorem 1.9 in \cite{Iv}.
The map $\delta\cln M_{n,k}^1\to M_{n,k}^0$ induces an equivariant map $X^1=\Sph_c(M_{n,k}^1)\to X^0=\Sph_c(M_{n,k}^0)$ on the complexes of nonseparating sphere systems in these manifolds. In what follows, let $G_{n,k}^s$ denote either $A_{n,k}^s$ or $\Ga_{n,k}^s$. The group $G_{n,k}^s$ acts transitively on the set of $p$-simplices of $X^s$ for each $p$ and  $s=0,1$ as the complements of any two nonseparating systems of $p+1$ spheres are diffeomorphic.   

The double complex $E_*G_{n,k}^s\otimes_{G_{n,k}^s} C_*(X^s)$, with $E_*G_{n,k}^s$ a free resolution of $\Z$ over $\Z G_{n,k}^s$ and $C_*(X^s)$ the (nonaugmented) chain complex of $X^s$, defines two spectral sequences as in Section~\ref{spectralseq}. One of the two spectral sequences has $E^1$-term 
$${^s E}^1_{p,q}=H_q(\St^s(\s_p),\Z_{\s_p}) \ \Longrightarrow\  H_{p+q}^{G^s_{n,k}}(X^s,\Z)$$
where the spectral sequence now converges to the equivariant homology of $X^s$, and where $E^1_{p,q}$ has only one term as the action is transitive. The action of the stabilizer $\St^s(\s_p)$ on $\Z_{\s_p}$ is through the action on the orientation of $\s_p$ (which is not pointwise fixed by the stabilizer).
As $X^s$ is $(n-2)$-connected \cite[Prop.~3.2]{HW}, $H_{i}^{G^s_{n,k}}(X^s,\Z)\cong H_{i}(G_{n,k}^s,\Z)$ when $i\le n-2$. 

The map $\delta$ induces a map of spectral sequences $\delta^1_{p,q}\cln  {^1E^1_{p,q}}\to {^0E^1_{p,q}}$ converging to 
$$\delta_*\cln H_{p+q}^{G^1_{n,k}}(X^1,\Z) \to H_{p+q}^{G^0_{n,k}}(X^0,\Z)$$ which is the map we are interested in when $p+q\le n-2$. 
By \cite[Thm.~1.2]{Iv}, $\delta_*$ is an isomorphism in the range $p+q\le N$ if $\delta_{p,q}^1$ is an isomorphism when $p+q\le N$ and a surjection when $p+q= N+1$ with $p\ge 1$.

The stabilizer $\St^s(\s_p)$ of a $p$-simplex $\s_p$ fits into a short exact sequence 
$$\widetilde{\St^s}(\s_p) \rar \St^s(\s_p) \rar \Si_2\wr \Si_{p+1}$$
where $\widetilde{\St^s}(\s_p)$ is the subgroup of $\St^s(\s_p)$ that fixes the vertices of $\s_p$ and their orientation. 
The Hochschild-Serre spectral sequence for this short exact sequence has $E^1$-term 
$$
{^s E}^1_{r,t}=F_r \otimes_{\Si_2\wr\Si_{p+1}}H_t(\widetilde{\St^s}(\s_p),\Z_{\s_p}) \ \Longrightarrow\ H_{r+t}(\St^s(\s_p),\Z_{\s_p})
$$ 
where $F_r$ is a projective resolution of $\Z$ over $\Z[\Si_2\wr \Si_{p+1}]$.
The map $\delta$ induces a map of spectral sequences $\delta^1_{r,t}\cln  {^1E^1_{r,t}}\to {^0E^1_{r,t}}$. 
 Note that the action of $\widetilde{\St^s}(\s_p)$ on $\Z_{\s_p}$ is trivial, so that the coefficients in $H_*(\widetilde{\St^s}(\s_p),\Z_{\s_p})$ are actually untwisted.

\smallskip 

\noindent
{\bf Case 1:} If $G_{n,k}^s=A_{n,k}^s$, then $\widetilde{\St^s}(\s_p)\cong A_{n-p-1,k}^{s+2p+2}$. By Corollary~\ref{mu}, $\delta^1_{r,t}$ is an isomorphism for all $r$ when $t\le \frac{n-p-3}{2}$ as $\delta_*=\mu_*^{-1}$ when the manifold has at least two boundary spheres. By \cite[Thm.~1.2]{Iv}, it follows that $\delta$ induces an isomorphism $H_{q}(\St^1(\s_p),\Z_{\s_p})\to H_{q}(\St^0(\s_p),\Z_{\s_p})$ for all $q\le \frac{n-p-3}{2}$. Hence in the first spectral sequence, $\delta^1_{p,q}$ is an isomorphism in this range. 
The theorem follows in that case from applying \cite[Thm.~1.2]{Iv} to the first spectral sequence with the bound $N=\frac{n-4}{2}$.

\smallskip 

\noindent
{\bf Case 2:} If $G_{n,k}^s=\Ga_{n,k}^s$, there is an additional  short exact sequence 
$$(\Z/2)^{p+1} \rar G^{s+2p+2}_{n-p-1,k} \rar \widetilde{\St^s}(\s_p).$$
The associated Hochschild-Serre spectral sequence has $E^2$-term 
$$
{^s E}^2_{r,t}=H_r\big(\widetilde{\St^s}(\s_p),H_t((\Z/2)^{p+1},\Z)\big) \Rightarrow H_{r+t}(G^{s+2p+2}_{n-p-1,k},\Z)
$$
and $\delta$ induces a map of spectral sequences $\delta^2_{r,t}\cln  {^1E^2_{r,t}}\to {^0E^2_{r,t}}$. 
By Corollary~\ref{mu}, this map converges to an isomorphism $\delta_*\cln  H_i(G^{1+2p+2}_{n-p-1,k},\Z)\to H_i(G^{2p+2}_{n-p-1,k},\Z)$
in the range $n-p-1\ge 2i+2$. We also have that 
$\delta^2_{0,t}$ is an isomorphism for all $t$ because $(\Z/2)^{p+1}$ is central in $G^{s+2p+2}_{n-p-1,k}$. By \cite[Thm.~1.3]{Iv}, it  follows that $\delta^2_{r,0}\cln H_r(\widetilde{\St^1}(\s_p),\Z)\to H_r(\widetilde{\St^0}(\s_p),\Z)$ is an isomorphism when $r\le \frac{n-p-3}{2}$ and we can finish the proof as in Case 1. 
\end{proof}

The remaining stabilization is $\epsilon$. An equivalent inclusion $M^s_{n,k}\to M^s_{n,k+1}$ is obtained by attaching a $1$-handle and a disjoint $2$-handle to a sphere boundary component of $M^s_{n,k}$.  Attaching just the $2$-handle gives an inclusion $M^s_{n,k}\to M^{s+1}_{n,k}$ equivalent to $\mu$. Then attaching the $1$-handle gives an inclusion $M^{s+1}_{n,k}\to M^s_{n,k+1}$ which we denote by $\gamma$. Note that $\gamma$ is defined even when $s=0$. The composition $\gamma\mu$ is $\epsilon$.

\begin{thm}\label{gamma}
For $s\ge 0$, the maps $ \gamma_*\cln H_i(\Ga_{n,k}^{s+1}) \to H_i(\Ga_{n,k+1}^s)$ and $ \gamma_*\cln H_i(A_{n,k}^{s+1}) \to H_i(A_{n,k+1}^s)$ are isomorphisms when $n\ge 2i+1$. 
\end{thm}

Since $\epsilon = \gamma\mu$ the next result follows from the above and Corollary~\ref{mu}:

\begin{cor}\label{epsilon}
For $s\ge 1$ the maps $\epsilon_*\cln H_i(\Ga^s_{n,k})\to H_i(\Ga^s_{n,k+1})$ and $\epsilon_*\cln H_i(A^s_{n,k})\to H_i(A^s_{n,k+1})$ are isomorphisms when $n\ge 2i+2$. 
\end{cor}

\begin{proof}[Proof of Theorem~\ref{gamma}] 
Note that $\gamma$ always induces injections on homology since it has a left inverse given by filling in part of the new boundary torus, so we only need to prove surjectivity. Let $G_{n,k}^s$ again denote either $A_{n,k}^s$ or $\Ga_{n,k}^s$. 

Consider the action of $G_{n,k+1}^{s}$ on $\D_c(M^s_{n,k+1},C)$, the complex of nonseparating systems of disks with boundary on a fixed circle $C$ in the $(k+1)$-st boundary torus of $M^s_{n,k+1}$ (see \cite[Sect.~3]{HW}). The stabilizer of a simplex fixes the simplex pointwise since the disks in a system have an intrinsic ordering near their boundary. The stabilizer $\St(\s_p)$ of a $p$-simplex for $p\ge 0$ is isomorphic to $G_{n-p,k}^{s+p+1}$.
The inclusion $\St(\s_p) \to G_{n,k+1}^{s}$ is the composition $\gamma\,\eta^p$. Since $\alpha=\eta\mu$, we know that $\eta$ induces an isomorphism on homology when $\alpha$ and $\mu$ do, and a surjection when $\alpha$ does.  Thus $H_q(\St(\s_p)) \to H_q(G_{n,k+1}^{s})$ is an isomorphism for all $q<i$ when $n-p\ge 2q+2$ and a surjection when $n-p\ge 2q+1$ by Theorem~\ref{alpha} and Corollary \ref{mu}  and by induction on $i$ for $\gamma$. 

The action of $G_{n,k+1}^{s}$ on $\D_c(M^s_{n,k+1},C)$ determines a double complex 
$$E_*G_{n,k+1}^{s}\otimes_{G_{n,k+1}^{s}}\widetilde{C}_*(\D_c(M^s_{n,k+1},C))$$
and hence two spectral sequences as in Section~\ref{spectralseq}.
The action is transitive on simplices of each dimension, so (when $p\le n$) one spectral sequence has $E^1_{p,q}=H_q(\St(\s_p))$ (with untwisted coefficients as the stabilizer of a simplex fixes the simplex pointwise). So $E_{p,q}^1\cong H_q(G_{n,k+1}^{s})$ when $q<i$ and $n-p\ge 2q+2$. The differentials are given by an alternating 
sum of conjugations $c_g$ by elements $g\in G_{n,k+1}^{s}$, one for each face of the simplex, and these fit into a commutative diagram
$$\xymatrix{H_q(St(\s_p)) \ar[r]^{c_g} \ar[d]^\cong & H_q(St(\s_{p-1}))\ar[d]^\cong \\
H_q(G_{n,k+1}^{s}) \ar[r]^{id} & H_q(G_{n,k+1}^{s}) }$$
Hence the chain complex $(E^1_{*,q},d^1)$ is isomorphic to the augmented singular chain complex of a point with (untwisted) coefficients in $H_q(G_{n,k+1}^{s})$ 
in the range $*\le n-2q-2$ and surjects onto it when $*=n-2q-1$. It follows that $E^2_{p,q}=0$ when $p\le n-2q-2$ and $q<i$. 
In particular, this is the case if $p+q=i$ with $q<i$ as $p+2q+2=2i-p+2$ with $p\ge 1$ and $n\ge 2i+1$.

Now the complex $\D_c(M^s_{n,k+1},C)$ is $(n-1)$-connected \cite[Thm.~3.1]{HW}. Using the other spectral sequence, one gets $E^\infty_{p,q}=0$ when $p+q\le n-1$. In particular, 
$E^1_{-1,i}$ must die if $i\le n$. By the above calculation, it can only be killed by $d^1\cln E^1_{0,i}\to E^1_{-1,i}$, which 
is the map $H_i(G_{n,k}^{s+1}) \to H_i(G_{n,k+1}^s)$ we are interested in. Hence this map is surjective. 
\end{proof}

\begin{ex}{\rm 
Let us describe a family of examples where Corollary~\ref{mu} fails to hold when the stabilization is with respect to a prime manifold $P$ different from $S^1\times S^2$.  The manifold $P$ will be the exterior of an arbitrary (possibly trivial) knot $K$ in $S^3$. Let the manifold $M=M_{n,s}$ be the connected sum of $n$ copies of $P$ and $s+1$ copies of $D^3$. We will think of $M$ as being obtained from the exterior of the link $L$ in $D^3$ consisting of $n$ separated copies of $K$ by deleting the interiors of $s$ disjoint balls $B_i$ in the complement of $L$.  We assume $n\ge 1$ and $s\ge 1$.

For each ball $B_i$ we will construct a homomorphism $\phi_i\cln \Gamma(M,R)\to \Z/2$ where $R$ is the union of the boundary spheres of $M$.  As a preliminary remark, we note that diffeomorphisms of $M$ representing elements of $\Gamma(M,R)$ 
take longitudes of $L$ to longitudes, as these are the only nontrivial circles in $\del M$ that are homologically trivial in $M$.  Now to define $\phi_i$, choose a Seifert surface $S$ for $L$ consisting of $n$ disjoint orientable surfaces $S_j$ bounded by the components $L_j$ of $L$ and disjoint from the balls $B_k$.  Applying a diffeomorphism $f$ representing an element of $\Gamma(M,R)$ gives a new Seifert surface $f(S)$ by the preliminary remark above. 
Define $\phi_i(f)$ to be the class in $H_2(D^3-B_i;\Z/2)=\Z/2$ represented by the mod $2$ cycle $S + f(S)$.  (We use $\Z/2$ coefficients since $f$ may reverse orientations of the $L_j$'s.) It is clear that $\phi_i(f)$ depends only on the isotopy class of $f$ as a diffeomorphism of $M$. To see that $\phi_i(f)$ is independent of the choice of the surface $S$, let $S'$ be another choice for~$S$. We can assume the four surfaces $S$, $S'$, $f(S)$, and $f(S')$ have general position intersections, apart from their common boundary $L$.  Since $S$ and $S'$ represent the same element of $H_2(D^3,L;\Z/2)\cong H_1(L,\Z/2)$, their sum $S+S'$ is the mod $2$ boundary of a chain, which we can take to be a simplicial chain with respect to some triangulation of $D^3$, that is, a region $E$ in $D^3$ whose geometric boundary is $S+S'$. Applying $f$, the region $f(E)$ has boundary $f(S)+f(S')$. Thus we have $\bigl(S+f(S)\bigr)+\bigl(S'+f(S')\bigr)=\del\bigl(E+f(E)\bigr)$ in $D^3$.  The ball $B_i$ is fixed by $f$, so $B_i$ is contained either in both $E$ and $f(E)$ or in neither.  Hence $\bigl(S+f(S)\bigr)+\bigl(S'+f(S')\bigr)$ is a mod $2$ boundary in $D^3-B_i$, which shows that $\phi_i(f)$ does not depend on the choice of $S$.  It follows that $\phi_i$ is a homomorphism since if $g$ is another diffeomorphism then $\phi_i(fg)$ is represented by $S+fg(S)=\bigl(S+g(S)\bigr)+\bigl(g(S)+f(g(S))\bigr)$ which represents $\phi_i(g)+\phi_i(f)$.

To see that $\phi_i$ is surjective, consider the diffeomorphism $f_i$ of $M$ obtained by dragging $B_i$ around a loop in $D^3 -L$ that intersects a given Seifert surface $S$ in one point transversely and is disjoint from the other $B_j$'s. Such a loop exists since $D^3-L$ is connected.  The effect of $f_i$ on $S$ is to change it only by replacing a disk $D\subset S$ by another disk $f_i(D)$ with the same boundary but lying on the other side of $B_i$.  Thus the cycle $S+f_i(S)$ is a sphere isotopic to $\del B_i$. 
This is nonzero in $H_2(D^3-B_i;\Z/2)$ and zero in $H_2(D^3-B_j;\Z/2)$ for $j\ne i$, so $\phi_i(f_i)\ne 0$ and $\phi_j(f_i)=0$ for $j\ne i$. Thus the maps $\phi_i$ are the components of a surjective homomorphism $\Phi\cln \Gamma(M,R)\to (\Z/2)^s$.  There is a commutative square 
$$\xymatrix{\Gamma(M_{n,s},R) \ar[r]^\mu \ar[d]^\Phi & \Gamma(M_{n,s+1},R) \ar[d]^\Phi\\
(\Z/2)^s\ar[r]^\psi & (\Z/2)^{s+1}}$$
for a certain homomorphism $\psi$ which could be made explicit. Abelianizing the diagram, it follows that $\mu_*\cln H_1(\Gamma(M_{n,s},R))\to H_1(\Gamma(M_{n,s+1},R))$ cannot be surjective since $\psi$ cannot be surjective.

To see that this holds also for the quotient groups $A(M,R,T)$ it suffices to show that each $\phi_i$ vanishes on twists along spheres and disks in $M$. For twists along spheres, these are generated by twists along spheres in any maximal collection of spheres, and a Seifert surface can be chosen disjoint from such a maximal collection, so each $\phi_i$ is trivial on these generators.  For twists along disks, the only case when these arise is when the knot $K$ is trivial since otherwise the boundary tori of $M$ are incompressible.  For the trivial knot $K$ the twist disks form a Seifert surface for $L$, so by taking parallel copies of these twist disks as $S$ we see that $\phi_i$ is trivial on disk twists.
}\end{ex}

\section{braid groups, symmetric groups and manifolds of other dimensions}\label{otherdim}

One of the stability statements in Theorem~\ref{main} does not require any 3-dimensional topology for its proof, namely the case of the stabilization with $P=D^3$. Indeed, to study the map  $\Ga(N_n^{D^3},R)\to \Ga(N_{n+1}^{D^3},R)$ we used a complex of arcs isomorphic to a join of copies of $\pi_1(M)$ (see Proposition~\ref{D3}), and the proof of high connectivity of this complex is completely combinatorial.  In this section, we prove Proposition~\ref{higherm} which extends this stabilization result to manifolds $M$  (not necessarily compact or orientable) of any dimension $m\ge 2$, with $\Ga(N_n^{D^3},R)$ replaced by $\Ga(M,n,R)$, the mapping class group of diffeomorphisms of $M$ permuting a set $\Pts_n$ of $n$ points interior to $M$ and fixing a submanifold $R$ of its boundary. We moreover prove Propositions~\ref{wreath} and \ref{braidwreath} concerning symmetric and braid groups, which are both closely related to the mapping class groups of punctured manifolds. 

The extension of puncture-stabilization from dimension $3$ to higher dimensions will be immediate, but dimension $2$ requires extra care due to the fact that arcs in general position in a surface can intersect, and thus the arc complexes are not join complexes as they are in dimension $3$ and above. Most of this section is concerned with proving the appropriate connectivity result in dimension $2$.

\medskip

Fix a connected surface $S$ (orientable or not, compact or not) with a boundary circle $\del_0S$. Let $\Delta_0$ be a finite set of points in $\del_0 S$ and let $\Pts_n=\{p_1,\cdotss,p_n\}$ be $n$ distinct points in the interior of $S$.
We consider first the simplicial complex $\mathcal{F}(S;\Delta_0,\Pts_n)$ whose vertices are isotopy classes of embedded arcs in $S$ with one boundary point in $\Delta_0$ and one in $\Pts_n$, and whose higher simplices are collections of such arcs which do not intersect, except possibly at their endpoints.

\begin{lem}\label{arcsurgery}
$\mathcal{F}(S;\Delta_0,\Pts_n)$ is contractible for all $n\ge 1$. 
\end{lem}

\begin{proof}
The proof uses the surgery argument introduced in \cite{H-triang}, in the variant employed in \cite[Lem~2.5]{NW} with $\Pts_n=\Delta_1$ in the notation there. The circle $\del_0S$ is a {\em pure} boundary component in the sense of \cite{NW} since it only contains points of $\Delta_0$. The fact that the points $p_1,\cdotss,p_n$ do not lie on boundary components does not play any role in the proof.
\end{proof}

We denote the vertices of $\mathcal{F}(S;\Delta_0,\Pts_n)$ as  pairs $(i,a)$, where $i\in \{1,\cdotss,n\}$ and $a$ is an arc from $\Delta_0$ to $p_i$. 
Let $A(S;\Delta_0,\Pts_n)$ be the subcomplex of $\mathcal{F}(S;\Delta_0,\Pts_n)$ of simplices $\lgl(i_0,a_0),\cdotss,(i_p,a_p)\rgl$ such that $i_j\neq i_k$ whenever $j\neq k$, that is, collections of arcs with at most one arc ending at each $p_i$. (This is the 2-dimensional analogue of the complex $X^A$ in Proposition~\ref{D3}.)

\begin{prop}\label{surfaces}
 $A(S;\Delta_0,\Pts_n)$ is Cohen-Macaulay of dimension $(n-1)$. 
\end{prop}

\begin{proof}
The complex $A(S;\Delta_0,\Pts_n)$ has dimension $n-1$ since a maximal simplex has $n$ arcs. We will show that it is $(n-2)$-connected. It will follow that the complex is Cohen-Macaulay as the link of a $q$-simplex $\s=\lgl(i_0,a_0),\cdotss,(i_q,a_q)\rgl$ is isomorphic to  $A(S-\s;\Delta'_0,\Pts_n-\{p_{i_0},\cdotss,p_{i_q}\})$, where $\Delta_0'$ is the image of $\Delta_0$ in $S-\s$. 

The proof of connectivity is analogous to the proof of Lemma~\ref{coloring}. First note that the result is trivially true for $n=0,1$. For larger $n$ we induct on $n$. Fix $n\ge 2$ and $k\le n-2$. Let $f\cln S^k\to A(S;\Delta_0,\Pts_n)$ be a map. By Lemma~\ref{arcsurgery}, we can extend $f$ to a map $\hat f\cln D^{k+1}\to \mathcal{F}(S;\Delta_0,\Pts_n)$, which we can assume to be simplicial for some triangulation of $D^{k+1}$. We will modify $\hat f$ so that its image lies in $A(S;\Delta_0,\Pts_n)$. 

Call a simplex $\s$ of $D^{k+1}$ {\em bad\/} if $\hat f(\s)=\lgl (i_0,a_0),\cdotss, (i_p,a_p)\rgl$  with each $i_j$  occurring at least twice with different $a_j$'s. Let $\s$ be a bad simplex of maximal dimension~$p$. By maximality of $\s$, $\hat f$ maps the link of $\s$ to collections of arcs in the complement of $\cup_ja_j$ connecting points of $\Delta_0$ to points $p_i$ disjoint from $\cup_ja_j$, and in fact $\hat f$ restricts to a map  
$$
\hat f \,\cln\, {\rm Link}(\s)\simeq S^{k-p}\ \rar\ J_\s=A(S_1;\Delta^1_0,\Pts_{n_1})*\cdotss *A(S_c;\Delta^c_0,\Pts_{n_c})
$$
where $S_{1},\cdotss, S_{c}$ are the surfaces obtained by cutting $S$ along $\cup_ja_j$, $\Delta^i_0$ is the image of $\Delta_0$ in $S_{i}$ and $\Pts_{n_i}=\Pts_n\cap \,{\rm int}(S_i)$. Note that  each $\Delta^i_0$ is nonempty and $n_1+\cdotss +n_c=n-d$ for $d$ the number of distinct $i_j$'s in $\hat f(\s)$. We have $d\le \lfloor \frac{p+1}{2}\rfloor\le p$ so that, by induction, $J_\s$  has connectivity $(\sum_j n_j) -2\ge n-p-2$. As $k-p\le n-p-2$, we can extend the restriction of $\hat f$ to the link of $\s$ to a map $g\cln D^{k-p+1}\to J_\s$. We modify $\hat f$ on the interior of the star of $\s$ using $\hat f * g$ on $\del\s * D^{k-p+1}\simeq \operatorname{Star}(\s)$. This construction removes $\s$ and does not add any new bad simplices of maximal dimension. The result follows by induction. 
\end{proof}

We are now ready to prove the puncture-stabilization in all dimensions. 
Let $M$ be an $m$-dimensional manifold with $\del_0M\subset R$ submanifolds of the boundary as before, and let $\Pts_n=\{p_1,\cdotss,p_n\}$ be a set of $n$ distinct points in the interior of $M$.  As in Section~\ref{results}, let $\G(M,n,R)$ denote the mapping class group of diffeomorphisms mapping $\Pts_n$ to itself and fixing $R$. 
To see that the map $\G(M,n,R)\to \Ga(M,{n+1},R)$ induced by gluing a punctured $D^m$ to $\del_0M$ along a disk $D^{m-1}\subset \del D^m$ is injective, consider the following diagram of short exact sequences: 
$$\xymatrix{\G(M,R\cup \Pts_n) \ar@{->}[r] \ar@{->}[d] & \G(M,n,R)\ar[r] \ar[d] & \Si_n \ar@{->}[d]\\
\G(M,R\cup \Pts_{n+1}) \ar@{->}[r] & \G(M,{n+1},R)\ar[r] & \Si_{n+1} 
}$$
The vertical map on the left is injective since it has a left inverse defined by forgetting the new puncture. The vertical map on the right is obviously injective, so the middle vertical map is injective as well.

\begin{proof}[Proof of Propositions~\ref{higherm}, \ref{wreath} and \ref{braidwreath}]
The three results follow from Theorem~\ref{range}. For Proposition~\ref{higherm}, we use the action of $G_{n}=\Ga(M,{n},R)$ on the complex $X_{n}$ of nonintersecting arcs in $M$ from a point $p_0\in \del_0M$ to distinct points of $\Pts_n$. We have $X_{n}\cong(\Delta^{n-1})^{\pi_1(M)}$ if ${\rm dim}(M)\ge 3$ (see Remark~\ref{pi1} in the case that ${\rm dim}(M)=3$) and $X_{n}= A(M;\{p_0\},\Pts_n)$ if $M$ has dimension $2$. For Proposition~\ref{wreath}, we use the action of $G_{n}=G\wr \Si_{n}$ on the complex $X_{n}=(\Delta^{n-1})^G$ defined by $(g_1,\cdotss,g_n;\alpha)\cdot(i,h)=(\alpha(i),g_ih)$, where $\alpha\in \Si_n$, $g_1,\cdotss,g_n,h\in G$ and $i\in \{1,\cdotss,n\}$ is a vertex of $\Delta^{n-1}$. Finally for Proposition~\ref{braidwreath}, we use the action of $G_{n}=G\wr B^S_{n}$ on $X_{n}=A(S;\{p_0\},\Pts_n)^G$, where  $B^S_n=\pi_1\Conf(S,n)$ acts via the inclusion $ B^S_n\inc \Ga(S,\Pts_n,\del S)$ coming from the long exact sequence of homotopy groups associated to the fibration $\Dif(S,\Pts_n,\del S)\to \Dif(S,\del S)\to \operatorname{Emb}(\Pts_n,S)/\Sigma_n=\Conf(S,n)$.

Each of the complexes $X_{n}$ above is of dimension $n-1$ and $(n-2)$-connected by Propositions~\ref{label} and \ref{surfaces} and the action of $G_{n}$ on $X_n$ is transitive on the sets of $p$-simplices for each $p$. 
For the complexes of the type $(\Delta^n)^G$, the stabilizer of a $p$-simplex can realize all the permutations of its vertices, whereas for the complexes of type $A(S;\{p_0\},\Pts_n)^G$, the stabilizer of a simplex fixes the simplex pointwise as the arcs are naturally ordered at $p_0$ in dimension 2, and this order cannot be changed by a diffeomorphism. The subgroup of $G_n$ fixing a $p$-simplex pointwise is in each case conjugate to $G_{n-p-1}$, as is clear in the case $m=2$, while for $m\ge 3$ this follows from the argument used to prove (iii) in Lemma~\ref{action}.  Condition (3) for Theorem~\ref{range} is also easily verified in each case since for an edge of $X_n$ defined by a pair of (labeled) arcs there is a diffeomorphism of $M$ fixed on $\del M$ and supported in a neighborhood of the two arcs which interchanges their endpoints in $\Pts_n$ and takes the first arc to the second.
\end{proof}

\section{connected sums along disks in the boundary}\label{boundary}
\setcounter{figure}{0}

In this section, we prove Theorem~\ref{disk-stabilization} about stabilization by boundary connected sum ($\nat$-sum) on irreducible 3-manifolds. Let $(M,\del_0M)$ be a compact, connected, oriented, irreducible 3-manifold with a chosen component $\del_0M$ of its boundary, and let $R$ be a finite collection of disjoint disks in $\del_0M$. Consider also a pair $(P,\del_0P)$ with $P$ prime with respect to $\nat$-sum.
The stabilization map $M\to M\nat P$ is obtained by identifying half of a disk in $\del_0P$ with half of a disk of $R$. It induces a map $\Ga(M,R)\to \Ga(M\nat P,R)$ which is injective: In an irreducible manifold each component of the space of disks with fixed boundary (or half of the boundary fixed) is contractible \cite{H1}, and in particular simply-connected. So given an isotopy $H$ joining a diffeomorphism $g$ of $M\nat P$ and the identity, with $g$ fixing $P$ as well as $R$, the isotopy $H$ can be deformed  (keeping $R$ fixed) to an isotopy fixing the disk splitting off $P$.

 We assume that $P$ is nontrivial, i.e. $P\neq D^3$. 
For ordinary connected sum, $S^1 \!\times\! S^2$-summands play a special role as they correspond to nonseparating spheres in the $3$-manifold. The analog for boundary connected sum is $S^1\!\times\! D^2$-summands, which correspond to disks in $M$ having boundary a nonseparating curve in $\del_0 M$. 

To prove Theorem~\ref{disk-stabilization}, we will use complexes $Y^A$ analogous to the complexes $X^A$ of Section~\ref{complexes}. 
To prove that the complexes $Y^A$ are highly connected, we follow the same strategy as for $X^A$, namely we deduce it from the connectivity of the related disk complexes. However, we cannot use Theorem~\ref{XtoY} directly here because the dual arcs live in a 2-dimensional surface, namely $\del_0M$, and thus the complexes $Y^A$ fail to be join complexes for the same reason as in the 2-dimensional case in the previous section.
We start the section by proving the necessary connectivity results for complexes of embedded disks. We then prove the connectivity of the complexes $Y^A$ using a combination of the techniques from Sections~\ref{connectivity} and \ref{otherdim}.

\medskip

Let $(M,\del_0M)$ be as above with $R$ a possibly empty finite collection of disjoint disks in $\del_0M$. 
We define $\D(M,R)$ to be the simplicial complex whose vertices are isotopy classes of nontrivial disks in $M$ with boundary in $\del_0 M-R$, and higher simplices are disjoint collections of such. Here a disk is nontrivial if it is neither isotopic to a disk of $R$ nor of $\del_0 M-R$. All isotopies must take place in the complement of $R$. When $R$ is empty we abbreviate $\D(M,R)$ to $\D(M)$.

\begin{prop}\label{DMA}
$\D(M,R)$ is contractible when $\D(M)$ is nonempty. 
\end{prop}

\begin{proof}
We proceed by induction on the number of disks in $R$, starting with the base case of $\D(M)$ when $R$ is empty. Contractibility in
this case is Theorem 5.3 in \cite{Mc}. It appears that the proof there might also work when $R$ consists of one disk, but let us instead show that when $R$ consists of one disk, the map $\D(M,R)\to \D(M)$ that ignores $R$ is a homotopy equivalence.  Choose a nontrivial annulus $A$ in $\del_0 M$ that contains the disk $R$ in its interior as in Fig.~\ref{AB}. (Note that $\del_0M$ is not a sphere since $M$ is irreducible and is not $D^3$ so that such an annulus exists.) Systems of disks defining simplices of $\D(M,R)$ can be isotoped so that their boundary curves meet $A$ minimally. The minimal intersections are unions of arcs as shown in the figure as there are no nontrivial disks with boundary a circle in $A$ by irreducibility of $M$. 
 Let $\D_0(M,R)$ be the subcomplex of $\D(M,R)$ formed by disk systems whose boundary curves cross $A$ from one component of $\del A$ to the other.  
\begin{figure}[htp]
\includegraphics[width=0.31\textwidth]{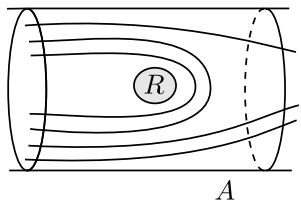}
\vskip-8pt
\caption{Discs intersecting $A$ minimally in $\D(M,R)$}\label{AB}
\vskip-15pt
\end{figure}
The full complex $\D(M,R)$ deformation retracts onto $\D_0(M,R)$ by pushing curves that do not cross $A$ across $R$
according to the following scheme. For a simplex $\s$ of $\D(M,R)$ not contained in $\D_0(M,R)$, let $D$ be the disk
of $\s$ which is not in $\D_0(M,R)$ and which is `innermost' in $A$, lying closest to $R$, and let $D'$ be the disk
obtained from $D$ by pushing it across $R$. Note that such a push preserves nontriviality of disks since $R$ is a
single disk. 
We can isotope $D'$ to be disjoint from $\s$ and to intersect $A$ minimally without moving $\s$.  
Let
$\tau$ be the simplex of $\D(M,R)$ spanned by $\s$ and $D'$, so $\tau=\s$ if $D'$ is already a vertex of $\s$.
Replacing $D$ by $D'$ can be realized by  a linear flow in $\tau$ along lines parallel to the edge from $D$ to $D'$,
giving a deformation retraction of $\tau$ onto its face $\tau - D$.  The restriction of this flow to any face of $\tau$
obtained by deleting a disk other than $D$ or $D'$ is the flow for that face, so we obtain in this way a well-defined
flow on $\D(M,R)$ that is stationary on $\D_0(M,R)$. Each simplex not in $\D_0(M,R)$ flows across finitely many
simplices until it lies in $\D_0(MR)$ since each flow across a simplex decreases the number of arcs of intersection
with $A$.  Thus the flow defines a deformation retraction of $\D(M,R)$ onto $\D_0(M,R)$.

The restricted map $f\cln\D_0(M,R)\to \D(M)$ can now be seen to be a homotopy equivalence because its fibers are contractible. Indeed, for any simplex $\s$ of $\D(M)$, consider the fiber $f_{\le 
\s}$, the subcomplex of $\D_0(M,R)$ of disk systems mapping to a face of $\s$. Suppose $\s=\s_0*\s_1$ with $\s_0$ the largest face of $\s$ disjoint from $A$. Choosing a lift $\hat \s$ of $\s$ to $\D_0(M,R)$ and pushing curves across $R$ gives an isomorphism $f_{\le\s}\cong \s_0*[\R\x\s_1]$, a contractible complex.

In the above arguments, we used the fact that pushing the boundary of a disk across $R$ preserves its nontriviality. This is only true when $R$ consists of a single disk. 
When $R$ has more than one disk, let $R_0$ be one disk of $R$ and $R_1$ the remaining disks. Forgetting $R_0$ does {\em not} define a map $\D(M,R)\rightarrow \D(M,R_1)$ since a disk splitting off a ball from $M$ whose intersection with $\del_0 M$ is a disk containing $R_0$ and only one disk of $R_1$ becomes trivial under this map. Call such a disk {\em special}. 
Note that the link of a vertex of $\D(M,R)$ corresponding to a special disk is isomorphic to $\D(M,R_1)$  and is thus contractible by induction. If $\D_0(M,R)$ denotes the subcomplex of $\D(M,R)$ with no vertices corresponding to special disks, then $\D(M,R)$ is the union of $\D_0(M,R)$ with the stars of the vertices corresponding to special disks. Since the links of these vertices are contractible, $\D(M,R)$ is homotopy equivalent to $\D_0(M,R)$. Furthermore, $\D_0(M,R)$ deformation retracts onto the link in $\D(M,R)$ of a fixed special vertex $D$ by pushing the boundary curves of all disks in $\D_0(M,R)$ off $\del D$ by pushing them across $R_0$ as shown in Fig.~\ref{A0}. This link is again contractible by induction. 
\end{proof}

\begin{figure}[htp]
\vskip-4pt
\includegraphics[width=0.33\textwidth]{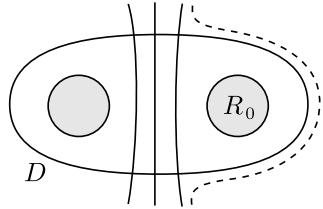}
\vskip-8pt
\caption{Pushing disks away from a special disk $D$ in $\D(M,R)$}\label{A0}
\vskip-15pt
\end{figure}

\subsection{Prime summands $\bf P\neq S^1\x D^2$}

Analogous to the complex $X$ in Section \ref{complexes}, let $Y$ be the complex whose simplices are isotopy classes of collections of
disjoint $P$-summands of $M$ disjoint from $R$. We consider here summands with respect to boundary connected sum along $\del_0M$ and
$\del_0P$, so such summands are cut off by disks in $M$ with boundary in $\del_0 M$. 

\begin{prop}\label{DPMA}
If $P\neq S^1\x D^2$ then $Y$ is Cohen-Macaulay of dimension $n-1$ for $n$ the number of $P$-summands in $M$.
\end{prop}

\begin{proof}
This follows the proof of Proposition~\ref{sphP} very closely so we will not go through the argument in detail. The initial special case $P=D^3$ in Proposition~\ref{sphP} has no analog in the present situation, so the proof begins with the case $n=1$, where the result is automatic. For $n\ge 2$, contractibility of $\D(M,R)$ replaces that of $\Sph(M)$, using the fact that $\D(M)\neq \emp $ as soon as $n\ge 2$.  At the step where the manifolds $M_1,\cdotss,M_d$ arise by cutting $M$ along disks, copies of these disks should be added to the disk-sets $R$ for the manifolds $M_i$.
\end{proof}

\smallskip

 Now we assume $R\ne\emp$ and we define an enhanced version $Y^A$ of the complex $Y$.  Choose a distinguished disk $B\subset \del_0 P$ and let $I\vee P$ be the union of $I$ and $P$, glued along $1\in I$ and $p_0\in\del B$. Consider embeddings $f\cln I\vee P \to M$ such that $f(P)$ is a $P$-summand of $M$ disjoint from $R$, cut off by the disk $f(B)$, and $f(I)$ is an arc in $\del_0M$ intersecting $R$ only in a chosen basepoint $x_0=f(0)\in\del R$. We also assume that $f$ is orientation-preserving on $P$.  Define $Y^A=Y^A(M,P,R,x_0)$ to be the simplicial complex whose vertices are isotopy classes of such embeddings $f$, with simplices corresponding to collections of embedding with images that are disjoint except at $x_0$.  There is a forgetful map $ Y^A\to Y$.

The following result is proved in Section~\ref{theproof} below.

\begin{thm}\label{VPM1}
When $P\neq S^1\x D^2$, the complex $Y^A$ is $(\frac{n-3}{2})$-connected for $n$ the number of $P$-summands in $M$. 
\end{thm}

\medskip

\subsection{Prime summands $\bf P=S^1\x D^2$}

Let $\D_c(M,R)$ denote the subcomplex of $\D(M,R)$ whose simplices are collections of disks $\lgl D_0,\cdotss,D_p\rgl$ whose boundaries form a nonseparating curve system $\lgl \del D_0,\cdotss,\del D_p\rgl$ in $\del_0 M$. 

\begin{prop}\label{DcMA}
$\D_c(M,R)$ is Cohen-Macaulay of dimension $n-1$ for $n$ the number of $S^1\x D^2$-summands in $M$.
\end{prop}

\begin{proof}
Note first that the complex has dimension $n-1$ 
as the neighborhood of $k$ disks with nonseparating boundary curve system and a dual graph to their boundaries in $\del_0M$ is a connected sum of $k$ $S^1\x D^2$-summands. 
Moreover the link of a $p$-simplex $\s=\lgl D_0,\cdotss,D_p\rgl$ in $\D_c(M,R)$ is isomorphic to $\D_c(M_\s,R\cup_i D_i^\pm)$, where $M_\s$, the manifold obtained from $M$ by cutting along $\s$, has $n-p-1$ $S^1\x D^2$-summands. (Here $D_i^\pm$ is the pair of disks in $\del M_\s$ corresponding to $D_i$.)  Hence the CM property will follow from the fact that $\D_c(M,R)$ is $(n-2)$-connected. 

As in the proof of Proposition~\ref{sphP}, we proceed by induction on the complexity of $M$, which is now the number of disks defining a maximal simplex of $\D(M,R)$. Let $k\le n-2$. A map $f\cln S^k\to \D_c(M,R)$ can be extended to a map $\hat f\cln D^{k+1}\to \D(M,R)$ by Proposition~\ref{DMA}, which is applicable here since we can assume $n\ge 2$. We want to modify $\hat f$ so that its image lies in $\D_c(M,R)$. Let $\s$ be a simplex of maximal dimension $p$ in $D^{k+1}$ such that $\hat f(\s)=\lgl D_0,\cdotss,D_q\rgl$ with each $\del D_i$ separating $\del_0M-\{\del D_0,\cdotss,\widehat{\del D_i},\cdotss,\del D_q\}$. Let $M_1,\cdots,M_d$ be the manifolds obtained by cutting $M$ along $\hat f(\s)$. By maximality of $\s$, the restriction of $\hat f$  to the link is a map $\link(\s)\simeq S^{k-p}\to J_\s=\D_c(M_1,R_1)*\cdotss *\D_c(M_d,R_d)$. The $M_i$'s have smaller complexity than $M$ and a total of at least $n-p$ $S^1\x D^2$-summands. Indeed, the dual graph $\Gamma$ of $\hat f(\s)$ in $\del_0 M$ has genus at most $q\le p$. A neighborhood of $\Gamma\cup \hat f(\s)$ is a handlebody of the same genus as $\Ga$ and $M$ is a boundary connected sum of this handlebody with the manifolds $M_1,\cdotss,M_d$. It follows that $J_\s$ has connectivity at least $n-p-2$. Hence we can extend $\hat f|_{\link(\s)}$ to a disk $D^{k-p+1}$ and use this extension to create a new $\hat f$ with fewer bad simplices of top dimension like $\s$. Repeating this process a finite number of times will give a map $\hat f$ with image in $\D_c(M,R)$. 
\end{proof}

Let $x_0, x_1$ be two points in $\del R$ (possibly $x_0=x_1$) and let $I+D^2$ denote the union of $I$ and $D^2$ with the midpoint of $I$ identified with $p_0\in \del D^2$. 
When $P=S^1\x D^2$, define $Y^A=Y^A(M,R,x_0,x_1)$ to be the simplicial complex whose vertices are isotopy classes of embeddings 
$$f\cln (I+D^2,I+\del D^2)\to (M,\del_0M)$$ 
smooth on $I$ and $D^2$ with $f(I)$ transverse to $f(\del D^2)$, 
such that $f(\del D^2)$ does not separate $\del_0M$ and $f(I+\del D^2)$ intersects $R$ only at $f(0)=x_0$ and $f(1)=x_1$. We also place an orientation condition on $f$ like the one for $X^A$ when $P=S^1\x S^2$. 
A $k$-simplex of $Y^A$ is represented by a collection $\lgl f_0,\cdotss,f_k\rgl$ of pairwise disjoint such embeddings such that the union of the images $f_i(\del D^2)$ does not separate $\del_0M$. There is a forgetful map $Y^A \to \D_c(M,R)$. In Section~\ref{theproof} we will prove:

\begin{thm}\label{VPM2}
When $P=S^1\x D^2$ and $x_0=x_1$, the complex $Y^A$ is $(\frac{n-3}{2})$-connected for $n$ the number of $S^1\x D^2$-summands in $M$. 
\end{thm}

\begin{figure}[htp]
\includegraphics[width=0.5\textwidth]{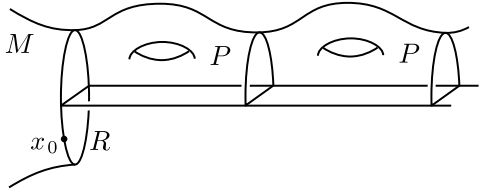}
\caption{Construction of $M_\infty$}\label{Fig4}
\vskip-15pt
\end{figure}

Besides the case $x_0=x_1$ we also need the case that $x_0$ and $x_1$ are distinct and lie in different components of $R$. In this case we consider a stable complex: 
Let $M$ be an irreducible 3-manifold as before and let $M_n$ be the manifold obtained from $M$ by boundary connected sum with $n$ copies of $S^1\x D^2$. 
We consider the stabilization $\alpha\cln M_n\to M_{n+1}$ identifying half of a disk $D$ in $S^1\x \del D^2$ with half of a disk in the
component of $R$ containing $x_0$ and disjoint from $x_0$, as in Figure~\ref{Fig4}. It induces an inclusion of simplicial complexes
$Y^A(M_n,R_n) \to Y^A(M_{n+1},R_{n+1})$ where $R_{n+1}$ is obtained from $R_n$ by replacing the half-disk where $S^1\x D^2$ attaches by
the unattached half of the disk $D$ in $S^1\x\del D^2$. Let $Y^A(M_\infty,R_\infty)$ be the direct limit  
$$
\operatorname{colim}(Y^A(M,R)\stackrel{\alpha}{\rar} Y^A(M_1,R_1)\stackrel{\alpha}{\rar} \cdotss\stackrel{\alpha}{\rar} Y^A(M_n,R_n)\stackrel{\alpha}{\rar}\cdotss)
$$ 
 where $M_\infty = \cup_n M_n$ and $R_\infty=\cup_nR_n \cap \del M_\infty$.

\begin{thm}\label{x0x1}
If $P=S^1\x D^2$ and $x_0$ and $x_1$ lie in different components of $R$ then $Y^A(M_\infty,R_\infty)$ is contractible. 
\end{thm}

\medskip

\subsection{Proof of Theorems~\ref{VPM1}, \ref{VPM2} and \ref{x0x1}}\label{theproof}

To prove the three theorems, we need to expand $Y^A$ to a larger complex $\hat Y^A$ having the same vertices as $Y^A$ but more higher-dimensional simplices. If $P\ne S^1\times D^2$ a simplex of $\hat Y^A$ is represented by a collection of embeddings $f_i\cln I\vee P\to M$ as before, but where the submanifolds $f_i(P)$ are now either disjoint or coincide, and the interiors of the arcs $f_i(I)$ are disjoint and disjoint from the submanifolds $f_j(P)$.  For the case $P=S^1\times D^2$ the definition is a little more complicated:  A simplex of $\hat Y^A$ is represented by a collection of embeddings $f_i\cln I\vee D^2 \to M$ where the disks $D_i=f_i(D^2)$ are disjoint or coincide, but we assume there exist annular neighborhoods $N_i$ of the circles $\del D_i$ in $\del_0 M$, with $N_i=N_j$ if $D_i=D_j$ and $N_i\cap N_j = \emp$ if $D_i\ne D_j$, such that each arc $a_i=f_i(I)$ intersects $\cup_jN_j$ in a single arc that crosses $N_i$ from one circle of $\del N_i$ to the other, and different $a_i$'s are disjoint outside $\cup_jN_j$ except at their endpoints at $x_0$ and $x_1$.  No restrictions are placed on how different $a_i$'s intersect within $\cup_jN_j$. For example, a simplex of $\hat Y^A$ can be constructed from a pair $(D_0,a_0)$ by adding parallel copies $a_1,\cdotss,a_p$ of $a_0$, then modifying each of these copies by applying to it a different power of a Dehn twist along $\del D_0$, and taking $D_0$ as the disk $D_i$ for each $i$.  The expected definition of $\hat Y^A$ would have required that different $a_i$'s intersect only at their endpoints, but this more restrictive definition leads to difficulties, as will be explained at the end of the proof of the next lemma.

We have defined the complex $Y$ for $P\ne S^1\x D^2$, and to unify notation we set $Y=\D_c(M,R)$ when $P=S^1\x D^2$.  Thus in both cases we have a projection $\hat Y^A \to Y$.  We denote by $Y_m$ the subcomplex of the barycentric subdivision of $Y$ corresponding to the subposet of simplices of $Y$ with at least $m$ vertices. Let $\hat Y^A_m$ denote the subcomplex of the barycentric subdivision of $\hat Y^A$ projecting to $Y_m$.

\begin{lem}\label{hatVm}
$\hat Y^A_m$ is $(n-m-1)$-connected, where we assume that $x_0=x_1$ if $P=S^1\x D^2$.
\end{lem}

\begin{proof}
We know that $Y_m$ is $(n-m-1)$-connected by Propositions~\ref{DPMA}, \ref{DcMA} and Lemma~\ref{genlargement}.
Thus to prove the lemma it suffices to show that the projection $\pi\cln \hat Y^A_m \to Y_m$ is a homotopy equivalence.  Consider first the case  $P\neq S^1\x D^2$.  Let $\s=\lgl [f_0(P)],\cdotss,[f_k(P)]\rgl$ be an element of $Y_m$. The fiber ${\pi}_{\ge\s}$ is the subposet of elements $\lgl [g_0],\cdotss,[g_r]\rgl$, with $g_j\cln I\vee P\to M$, such that each  $f_i(P)$ occurs among the $g_j(P)$'s. There is a projection ${\pi}_{\ge\s}\to \pi^{-1}(\s)$ defined by forgetting $g_j$'s for which $g_j(P)$ is not among the $f_i(P)$'s. The fibers of this projection are contractible since they have a minimal element. 
Hence to show that $\pi$ is an equivalence, it is enough to show that $\pi^{-1}(\s)$ is contractible for any $\s$ in $Y_m$. The argument is an adaptation of that used to prove Lemma~\ref{arcsurgery}.

Choose a lift $\hat\s=\lgl [f_0],\cdotss,[f_k]\rgl$ of $\s$ in $\hat Y^A_m$, with $f_i\cln I\vee P\to M$. We are going to define a deformation retraction of $\pi^{-1}(\s)$ to the point $\hat\s$. 
We think of an element $\tau$ of $\pi^{-1}(\s)$ as a family of arcs $b_1,\cdotss,b_r$, $b_j=g_j(I)$, which we can assume have been isotoped into normal form with respect to the arcs $a_i=f_i(I)$, i.e., intersecting the $a_i$'s minimally and transversely; such a normal form is unique up to isotopy through normal forms (except for moving a $b_j$ that is parallel to an $a_i$ from one side of $a_i$ to the other). 
If some $a_i$ meets a $b_j$ at a point other than $x_0$, consider the first such $a_i$ in the ordered list $a_0,\cdotss,a_k$ and look at the first intersection point of this $a_i$ with a $b_j$ as we move along $a_i$ away from $x_0$. Cutting the $b_j$ arc at this intersection point produces two arcs, and we discard the one going to $x_0$.  The one going to a disk $D_\ell$ we keep and extend to rejoin to $x_0$ along an arc parallel to $a_i$, giving a new arc $b'_j$ that is compatible with the original collection $b_1,\cdotss,b_r$ (see Fig.~\ref{Yhat}).
\begin{figure}[htp]
\includegraphics[width=0.68\textwidth]{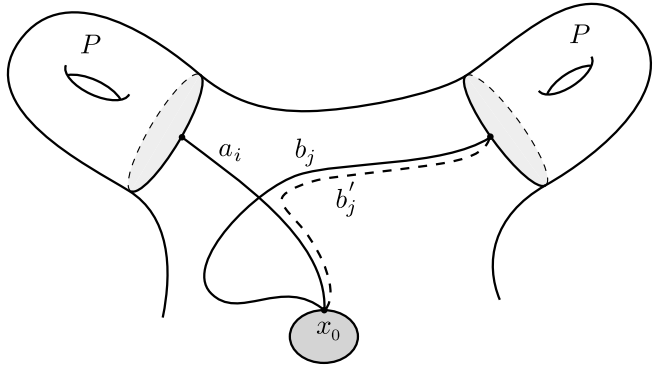}
\vskip-12pt
\caption{Surgery on the arcs $b_j$}\label{Yhat}
\vskip-15pt
\end{figure}
This step can be repeated for each remaining intersection along this $a_i$, then for the intersections along $a_{i+1}, a_{i+2},\cdotss$ until all such intersections are eliminated. This sequence of steps defines a path in $\hat Y^A_m$ by varying the weights from the old arc to the new arc each time, as in the proof of Proposition~\ref{DMA}. This path ends at a new system of arcs $b_j$ that meet the $a_i$'s only at $x_0$, and we can extend the path so that it ends at $\hat\s$ by shifting the weights from the new $b_j$'s to the $a_i$'s. The combined path depends continuously on $\tau$, and hence defines a deformation retraction of $\pi^{-1}(\s)$ to $\hat\s$.

If $P=S^1\x D^2$ the process is similar, doing surgery on the arcs $b_j$ to eliminate their intersections with $a_i$'s outside the annular neighborhoods of the circles $\del D_j$, surgering from first one end of an $a_i$ and then the other.  This is done after the $b_j$'s have first been isotoped to minimize their intersections with the $a_i$'s outside the annular neighborhoods. This minimal position is unique up to isotopy, so the surgery process is well-defined on isotopy classes.  This would not be the case if we also tried to surger away intersections of $b_j$'s with $a_i$'s inside the annuli since such intersections can be pushed from one side of the core circle $\del D_\ell$ of an annulus to the other, and eliminating the intersections by surgery on one side of $\del D_\ell$ or the other can produce different results.
\end{proof}

 Fortunately, allowing intersections in the annuli when $P=S^1\x D^2$ will cause no problems later: In the proof of Theorem~\ref{VPM2}, we will only need that arcs dual to different disks are disjoint. 
 
The above argument does not apply to the case $x_0\neq x_1$ when $P=S^1\x D^2$ because an arc $b_j'$ obtained from an arc $b_j$ by surgery as above could go from $x_0$ to itself or from $x_1$ to itself, even if $a_i$ and $b_j$ go from $x_0$ to $x_1$. It is however enough for our purposes to show that the high connectivity holds stably.

\begin{lem}\label{hatVm2} 
Suppose that $P=S^1\x D^2$ and that $x_0$ and $x_1$ lie in different components of $R$. Then $\hat Y^A_m(M_\infty,R_\infty)$ is contractible.
\end{lem}

\begin{proof} We follow the same general strategy as in the proof of the previous lemma in the case $P=S^1\x D^2$, but with a new {\it rerouting\/} construction replacing surgery. This rerouting technique can be viewed as an analog for surfaces of the lightbulb trick in $3$ dimensions that has played a key role earlier in the paper.

Consider the projection $\pi\col\hat Y^A_m(M_\infty,R_\infty)\to Y_m(M_\infty,R_\infty)$. Since the target space $Y_m(M_\infty,R_\infty)$ is contractible (Prop.~\ref{DcMA}), it suffices, as in the previous lemma, to show contractibility of the pre-image $\pi^{-1}(\sigma)$ of each simplex $\sigma = \lgl D_0,\cdotss,D_k\rgl$. Given a map $f\col S^p \to \pi^{-1}(\sigma)$, this has compact image lying in $\hat Y^A_m(M_n,R_n)$ for some finite~$n$. 
To do the rerouting, we choose a collection of $2k+2$ disks $T_{i,\epsilon}$ for $i=0,\cdotss,k$ and $\epsilon=0,1$, lying in $M_\infty - (M_n\cup R_\infty)$ and forming a nonseparating system with the $D_i$'s. We lift $\sigma$ to a simplex $\hat\sigma =\lgl (D_0,a_0)\cdotss,(D_k,a_k)\rgl$ in $\hat Y^A_m(M_\infty,R_\infty)$ where each $a_i$ is chosen to intersect the disk system of $D_i$'s and $T_{i,\epsilon}$'s transversely in three points: starting at $x_0$, $a_i$ first crosses $T_{i,0}$, then $D_i$, then $T_{i,1}$ before ending at $x_1$.  Such arcs $a_i$ exist since the $D_i$'s and $T_{i,\epsilon}$'s form a nonseparating system. 

We wish to deform $f$ to the constant map with image $\hat\s$. 
A simplex in the image of $f$ is represented by a collection of arcs $b_j$ dual to the $D_i$'s but disjoint from the $T_{i,\epsilon}$'s. After putting this collection into normal form with respect to the $a_i$'s using a neighborhood $N_i$ of each $D_i$ as before, consider an intersection of a $b_j$ with an $a_i$ outside $N_i$ that is closest to the disk $T_{i,\epsilon}$ for $\epsilon$ either $0$ or $1$. We reroute a small segment of $b_j$ near this intersection with $a_i$ so that it travels parallel to $a_i$ to a point near $T_{i,\epsilon}$, then around $\del T_{i,\epsilon}$ to a point on the other side of $a_i$, then parallel to $a_i$ back to the original $b_j$, as shown in Figure~\ref{sigma}. The new $b'_j$ can be isotoped to be disjoint from the original $b_j$, and we obtain a deformation of $f$ by varying the weights to replace $b_j$ by $b'_j$.  The process can be repeated until all intersections of $b_j$'s with $a_i$'s are eliminated (apart from those in the annular neighborhoods of the circles $\del D_i$ that we do not care about). After all these intersections have been eliminated we can do one last deformation of $f$ by shifting the weights on the resulting new arcs $b_j$ to the arcs $a_i$, giving the constant map with image $\hat\s$. 
\begin{figure}[htp]
\vskip-4pt
\includegraphics[width=0.92\textwidth]{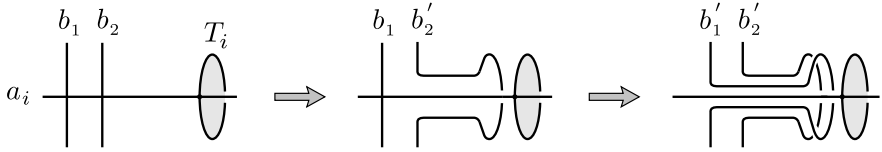}
\vskip-10pt
\caption{The simplex $\s$ and the deformation}\label{sigma}
\vskip-15pt
\end{figure}
\end{proof}

\begin{proof}[Proof of Theorems~\ref{VPM1}, \ref{VPM2} and \ref{x0x1}]
We will show that $Y^A=Y^A(M,R)$ is $(\frac{n-3}{2})$-connected ({\em general case}) when $M$ has $n$ $P$-summands in its prime decomposition, except in the case $P=S^1\x D^2$ with $x_0\neq x_1$ ({\em special case}), where we show that any map $S^k\to Y^A(M,R)$ can be filled with a disk in $Y^A(M_N,R_N)$ for some large $N$. Recall that $M_N$ is the boundary connected sum of $M$ with $N$ copies of $S^1\x D^2$. For notational simplicity we will usually omit the term $R_N$.

Suppose that $k\le \frac{n-3}{2}$ (or $k$ is any number in the special case) and let $f\cln S^k\to Y^A$ be a map. 
The strategy of the proof is similar to the one in the proof of Theorem~\ref{XtoY}: We first `enlarge' $f$ to a map $g\cln S^k\to \hat Y^A_{k+2}$ (after passing to $M_{N_1}$ in the special case). This map can be extended to a map from the disk $D^{k+1}$ by the known connectivity of $\hat Y^A_{k+2}$ in the general case (resp. after passing to $M_{N_1+N_2}$ in the special case). We then use a coloring argument to recover a map $G\cln D^{k+1}\to Y^A$ which is homotopic to $f$ on $S^k$. 

\smallskip

\noindent 
{\em Step 1: Construction of $g\cln S^k\to \hat Y^A_{k+2}$.}
The given map $f\cln S^k\to Y^A$ may be taken to be simplicial with respect to some triangulation $\T_0$ of $S^k$. We denote by $\T_0'$ the barycentric subdivision of $\T_0$. A $p$-simplex of $\T_0'$ is thus a chain  $[\s_0\!<\!\cdotss\!<\!\s_p]$ of simplices of $\T_0$. We want to construct a simplicial map $g$ from a subdivision $\T_1$ of $\T_0'$ to $\hat Y^A_{k+2}$ with the following {\em additional property\/}: For any vertex $v$ of $\T_1$ such that $v$ lies in the interior of a simplex $[\s_0\!<\!\cdotss\!<\!\s_p]$ of $\T_0'$, $g(v)$ contains $f(\s_0)$ as a {\em monic\/} subset, that is, a subset such that each $P$-summand or disk of $f(\s_0)$ intersects only one arc of $g(v)$. 
We will construct $g$ inductively over the skeleta of $\T'_0$. 

Consider first the general case. For each vertex $[\sigma]$ of $\T'_0$, with $f(\s)$ a $p$-simplex of $Y^A$, extend $f(\s)$ to a $(k+1)$-simplex $g([\sigma])=f(\sigma)*\tau$ of $Y^A$. This is always possible as every simplex of $Y^A$ can be extended to a (maximal) $(n-1)$-simplex and $k+1\le n-1$. We consider $g([\s])$ as a vertex of $\hat Y^A_{k+2}$. 
For the inductive step we wish to extend $g$ over a $p$-simplex $\tau=[\sigma_0 \!<\! \cdotss \!<\! \sigma_p]$ of $\T_0'$, assuming we have already defined $g$ on $\del\tau$ so that the additional property in the preceding paragraph is satisfied. In particular, for any vertex $v$ of $\del\tau$, $g(v)$ contains $f(\s_0)$ as a monic subset, and we can consider the restriction of $g$ to $\del\tau$ as a map $g_{\tau}\cln \del\tau \to \hat Y^A_{k+2-q}(M_{f(\s_0)})$, where $q$ is the number of vertices in $f(\s_0)$ and $M_{f(\s_0)}$ is $M$ with the $P$-summands or disks of $f(\s_0)$ removed, the resulting boundary disks being added to $R$. This complex has connectivity $n-q-(k+2-q)-1=n-k-3$ by Lemma~\ref{hatVm}. 
We need it to be $(p-1)$-connected in order to extend $g_{\tau}$ over the whole simplex $\tau $. Thus we need the inequality $p-1\le n-k-3$, which holds when $k\le (n-3)/2$ as $p\le k$. The extension of $g_\tau$ over $\tau$ may involve arcs that intersect arcs of $f(\s_0)$, but such intersections can be eliminated by rerouting the arcs in $g_\tau$ along the arcs of $f(\s_0)$ so that they go around the new disks of $R$.  Alternatively, we could use the surgery technique in the proof of Lemma~\ref{hatVm} to eliminate the intersections.  After this has been done we can combine the extended $g_\tau$ with $f(\s_0)$  on $\tau$ to give the induction step in the construction of $g$ on $S^k$, with the additional property still satisfied.

For the special case, the inductive step extending $g$ from the $(p-1)$-skeleton to the $p$-skeleton requires doing a connected sum with $n_p$ copies of $S^1\x D^2$ by Lemma~\ref{hatVm2}, for some $n_p$ large enough for each of the $p$-simplices of $S^k$ --- a finite number is sufficient by compactness. Taking $N_1=n_0+\cdotss +n_k$, we thus get a map $g\cln S^k\to \hat Y^A_{k+2}(M_{N_1})$. 

\smallskip

In the general case, extend $g$ to a map $g\cln D^{k+1}\to \hat Y^A_{k+2}$, which is possible as $\hat Y^A_{k+2}$ is $(n-(k+2)-1)$-connected. In the special case, extend $g$ to a map $g\cln D^{k+1}\to \hat Y^A_{k+2}(M_{N_1+N_2})$ for some $N_2$ large enough. In both cases, we can assume that $g$ is simplicial with respect to a triangulation $\T_1$ of $D^{k+1}$ which restricts to the already defined $\T_1$ on $S^k$. 

\smallskip 

\noindent 
{\em Step 2: Construction of $G\cln D^{k+1}\to Y^A$.}
We define $G$ inductively on the skeleta of $\T_1$, from a subdivision $\T_2$ of $\T_1$, with the following {\em additional property\/}: if $w$ is a vertex of $\T_2$ which lies in the interior of a simplex $\tau=\lgl v_0,\cdotss,v_p\rgl$ of $\T_1$ with $g(v_0)\le\cdotss\le g(v_p)$ in $\hat Y^A_{k+2}$, we require that 
$G(w)$ is a vertex of $g(v_0)$.  In case $w$ is in $S^k$, in the interior of a simplex  $[\s_0\!<\!\cdotss\!<\!\s_p]$ of $\T_0'$, we choose moreover $G(w)$ to be a vertex of $f(\s_0)$, which is possible as $f(\s_0)\subset g(v_0)$ by step~1. 

There is no obstruction to defining $G$ on the 0-skeleton of $\T_1$ satisfying the additional property. Note that if a vertex $w$ of $\T_1$ lies in a simplex $\sigma$ of $\T_0$, then we must have $G(w)=f(v)$ for some vertex $v$ of $\sigma$. As $\T_1$ is a subdivision of $\T_0$, it follows that $G$ extends over $S^k$ linearly on simplices of $\T_1$. This extension is linearly homotopic to $f$. We are left to define $G$ on higher simplices in the interior of $D^{k+1}$.

Suppose now that we have defined $G$ over the $(p-1)$-skeleton of $\T_1$.  Let $\tau=\lgl v_0,\cdotss,v_p\rgl$ as above be a $p$-simplex of $\T_1$ not contained in $S^k$. For each vertex $w$ of the new triangulation $\T_2$ of $\del\tau$, we have inductively defined $G(w)$ to be a vertex of $g(v_j)$ for some $j$. For the projection $\pi\cln \hat Y^A\to Y$ let $E$ be the set of vertices of $\pi g(v_p)$ and $E_0$ the vertices of $\pi g(v_0)$. 
We apply the coloring lemma (Lem.~\ref{coloring}) to the sphere $\del\tau$ triangulated by $\T_2$ with vertices labeled by $E$ via $\pi\circ G$. The result is an extension of the triangulation $\T_2$ over $\tau$ with the vertices in the interior of $\tau$ labeled by $E_0$ and bad simplices only in $\del\tau$. (Bad simplices may occur in $\del\tau$ if a face of $\tau$ is included in $S^k$.) For a vertex $w$ interior to $\tau$ we define $G(w)$ to be any lift of its label in $E_0$ to a vertex of $g(v_0)$.
To see that this definition is valid we need to check that for each simplex $\s=\lgl w_0,\cdotss,w_q\rgl$ of $\T_2$ in $\tau$ the vertices $G(w_i)$ span a simplex of $Y^A$. But $G(\s)$ is a face of $g(v_p)$ such that $\pi G(w_i)\neq\pi G(w_j)$ for $i\neq j$ unless $w_i,w_j\in S^k$ in which case we could have $G(w_i)=G(w_j)$. Hence $G(\s)$ is a simplex of $Y^A$.
\end{proof}

\subsection{Proof of Theorem~\ref{disk-stabilization}}
The proof uses the action of $\G(M,R)$ on the complex $Y^A$. First we need some information about the stabilizers of simplices. For a
simplex $\s$ of $Y^A$ in the case $P\ne S^1\times D^2$ we let $M_\s$ be the submanifold of $M$ obtained by splitting off the
$P$-summands given by $\s$. In the case $P=S^1\times D^2$ with $\s$ a simplex of $Y^A$ we let $M_\s$ be the result of splitting $M$
along the disks in $\s$. In both cases the collection $R$ of boundary disks in $M$ is contained in $M_\s$ so we can use the same $R$
for $M_\s$ as for $M$.

\begin{lem}\label{actionY}
When $P\neq S^1\x D^2$, the action of $\Ga(M,R)$ on $Y^A$ is transitive on $p$-simplices for any $p$. The stabilizer of a $p$-simplex $\s$  fixes the simplex pointwise and is isomorphic to $\G(M_\s,R)$. 
\end{lem}

\begin{proof}
If $\s$ and $\s'$ are $p$-simplices of $Y^A$ then the manifolds $M_\s$ and $M_{\s'}$ are diffeomorphic by the uniqueness of prime $\nat$-decompositions.  The diffeomorphism can be chosen to be the identity on $R$ and to take the arcs in $\s$ to the arcs in $\s'$. After a further adjustment on the disks splitting off the $P$-summands the diffeomorphism can then be extended over these $P$-summands so as to preserve their parametrizations specified by $\s$ and $\s'$.  Thus the action is transitive on $p$-simplices.

The stabilizer $\St(\s)$ of a simplex $\s$ fixes the simplex pointwise since the arcs of $\s$ have a preferred ordering at $x_0$ and cannot be permuted by a diffeomorphism. There is a natural map $\G(M_\s,R)\to \St(\s)$ since diffeomorphisms of $M_\s$ fixing $R$ can be assumed to fix also the arcs of $\s$ and the disks splitting off the $P$-summands of $\s$. This map $\G(M_\s,R)\to \St(\s)$ is injective since it is equivalent to an iterate of the stabilization map $\G(M_\s,R)\to \G(M_\s\nat P,R)$ which is injective as noted at the beginning of this section. An element $g\in\St(\s)$ fixes the arcs and $P$-summands of $\s$ up to isotopy. The isotopy can be assumed to be the identity in a neighborhood of $x_0$ and then on the rest of $\s$ by the isotopy extension property. Hence the map $\G(M_\s,R)\to \St(\s)$ is also surjective.  
\end{proof}

\begin{lem}\label{actionYc}
When $P= S^1\x D^2$, the action of $\Ga(M,R)$ on $Y^A$ is transitive on vertices. The stabilizer of a $p$-simplex $\s$  fixes the simplex pointwise and is isomorphic to $\G(M_\s,R)$. The quotient space $Y^A/\Ga(M,R)$ is $(n-1)$-spherical for $n$ the number of $S^1\x D^2$-summands in $M$.
\end{lem}

\begin{proof}
The action is transitive on vertices as in the previous lemma, but not on simplices of higher dimension because the pairs of endpoints of the arcs of a $k$-simplex of $Y^A$ can come in all possible orderings at $x_0$ (and $x_1$ if $x_1\neq x_0$), and diffeomorphisms cannot change the ordering pattern of the pairs. The quotient $Q_n=Y^A/\Gamma$ has dimension $n-1$, and we claim that it is $(n-2)$-connected. It has the structure of a $\Delta$-complex since the ordering of the leftmost of the two endpoints of each arc $a_i$ (resp. the position of the arcs at $x_0$ if $x_0\neq x_1$) specifies a preferred ordering for the vertices of a simplex. The location of the rightmost endpoints among the leftmost endpoints (resp. the ordering at $x_1$) specifies the simplex in $Q_n$.  The natural map $Q_{n-1}\to Q_n$ is an inclusion, and its image contains the $(n-2)$-skeleton of $Q_n$. The inclusion $Q_{n-1}\to Q_n$ extends to a map of the cone on $Q_{n-1}$ to $Q_n$ by adjoining an extra arc with both endpoints to the right of all the other endpoints. These properties imply that $Q_n$ is $(n-2)$-connected.  (This argument follows Lemma~3.5 of \cite{H3}, an alternative to Harer's original argument for Lemma~3.3 in \cite{harer}.)
\end{proof}

\begin{proof}[Proof of Theorem~\ref{disk-stabilization}]
In the case $P\ne S^1 \x D^2$, Theorem~\ref{disk-stabilization} follows from Theorem~\ref{VPM1} and Corollary~\ref{range2} using the action of $\Ga(M,R)$ on $Y^A$.  Conditions (1) and (2) at the beginning of Section~\ref{spectralseq} are given by Lemma~\ref{actionY}, and condition (3) follows from the fact that a diffeomorphism taking one vertex of an edge to the other can be chosen to have support in a neighborhood of the union of $R$ and the arcs and $P$-summands representing the edge.

When $P=S^1 \x D^2$ we first use the complex $Y^A$ with $x_0=x_1$. The action of $\Ga(M,R)$ on $Y^A$ satisfies condition (1$'$) and (2)  at the beginning of Section~\ref{spectralseq} by Lemma~\ref{actionYc}, and (3) holds as in the former case. The first part of the theorem for $P=S^1\!\times\! D^2$ then follows from Theorem~\ref{VPM2} and Corollary~\ref{range2}. 

To show independence of the number of disks in $R$, we follow the same strategy as in \cite{HVW}. Consider first a manifold $M$ with a collection $R=D_0\cup D_1\cup\cdots\cup D_d$ of at least two disks in $\del_0M$, with points $x_0\in\del D_0$ and $x_1\in\del D_1$. Let $P=S^1\!\times\! D^2$
and $M_\infty=\operatorname{colim}(M\stackrel{\alpha}{\to} M_1\stackrel{\alpha}{\to} M_2\stackrel{\alpha}{\to}\cdotss)$ as before,
with $\alpha$ the stabilization map which identifies half a disk in $\del P$ with half a disk in $D_0$. There is another map $\beta\cln M\to M'\cong M\nat P$ induced by adjoining a ball with two disks on its boundary by identifying a half disk in each disk with half disks in each of $D_0,D_1$. The induced map on the mapping class groups is also denoted $\beta$. The maps $\alpha$ and $\beta$ commute (see Fig.~\ref{alphabeta}),
\begin{figure}[htp]
\includegraphics[width=0.46\textwidth]{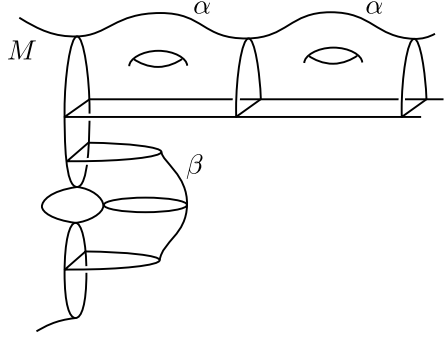}
\vskip-10pt
\caption{The stabilization maps $\alpha$ and $\beta$}\label{alphabeta}
\end{figure}
 so that $\beta$ induces a map $\beta\cln\G(M_\infty,R_\infty)\to \Ga(M'_\infty,R'_\infty)\cong \G(M_\infty,R_\infty)$. The action of $\G(M_\infty,R_\infty)$ on $Y^A(M_\infty,x_0,x_1)$ satisfies the stable conditions (1$'$-3) for Theorem~\ref{stable-stability} with $\lambda=\phi\circ\beta$ where $\phi$ is induced by a diffeomorphism $\phi\cln M'_\infty\to M_\infty$. Contractibility of $Y^A(M_\infty)$ is given by Theorem~\ref{x0x1} and the inclusion of the stabilizer of a vertex is easily seen to be conjugate to $\lambda$. The other conditions are given in Lemma~\ref{actionYc}. We conclude that $\beta$ induces isomorphisms on $H_i(\G(M_\infty,R_\infty))$ for all $i$. 
 
By the first part of the theorem and the following diagram, this implies that $\beta$ induces an isomorphism on $H_i(\G(M,R))$ in the same range as $\alpha$, namely $n\ge 2i+2$, and a surjection when $n=2i+1$:
$$\xymatrix{H_i(\G(M,R))\ar[r]^-\alpha \ar[d]^\beta & H_i(\G(M_1,R_1)) \ar[r]^-\alpha \ar[d]^\beta & H_i(\G(M_2,R_2))\ar[r]^-\alpha \ar[d]^\beta & \cdotss \\
H_i(\G(M',R')\ar[r]^-\alpha  & H_i(\G(M'_1,R'_1))\ar[r]^-\alpha  & H_i(\G(M'_2,R'_2))\ar[r]^-\alpha  & \cdotss}$$
This corresponds precisely to Theorem~\ref{alpha} if we replace $S^1\x D^2$-summands by $S^1\x S^2$-summands. As in Section~\ref{correction}, there are maps $\mu$ and $\eta$, now obtained by gluing a ball with 3 disks along one or two of its disks, and $\alpha=\eta\mu$ and $\beta=\mu\eta$. In particular, we get the analogue of Corollary~\ref{mu} which gives the independence of the number of disks of $R$ when $n\ge 2i+2$ and $R$ is not empty. To show that we can forget all the disks, we need the analogue of Theorem~\ref{delta}. Here we assume that $R$ consists of a single disk and $\delta\cln\G(M,R)\to \G(M)$ forgets that the disk is fixed. The proof of Theorem~\ref{delta}, case~2, applies directly to the new situation, replacing the complexes of spheres by the complexes $X^1=\D_c(M,R)$ and $X^0=\D_c(M)$ of disks with nonseparating boundaries. These complexes are $(n-2)$-connected by Proposition~\ref{DcMA}. We again have short exact sequences  $\widetilde{\St^s}(\s_p) \to \St^s(\s_p)\to \Si_{p+1}\wr \Si_2$ for the stabilizers $\St^0(\s_p),\St^1(\s_p)$ of $p$-simplices in $X^0$ and $X^1$. The second short exact sequence in the proof is replaced by $\Z^{p+1}\to \G^s_p\to\widetilde{\St^s}(\s_p)$, where $\G^0_p=\G(M_{\s_p},\s_p)$ and $\G^1_p=\G(M_{\s_p},R\cup \s_p)$. This comes from the fact that $\pi_1\Dif(D^2)\cong \Z$, replacing $\pi_1\Dif(S^2)\cong \Z/2$. 
\end{proof}

\section{appendix: proof of proposition 2.1}\label{proof}

\begin{proof}[Proof {\rm (}for reducible manifolds{\rm \/)}]
As in the proof of Proposition 2.2, we build $M$ from $S^3_n$ by attaching the manifolds $P^0_i$. We take the basepoint $x$ to lie in $S^3_n$. 

Let us first treat the case that $\del M$ contains no spheres, or equivalently, no $P_i$ is a ball. If an orientation-preserving diffeomorphism $f\cln (M,x)\to(M,x)$ induces the identity on $\pi_1$ then by the lemma below it also induces the identity on $\pi_2$. In particular, each sphere $f(S^2_i)$ is homotopic to $S^2_i$, so by Laudenbach's homotopy-implies-isotopy theorem, which applies to systems of disjoint spheres as well as individual spheres (Theorem III.1.3 and Lemma V.4.2 on p.124 of \cite{L}), we can isotope $f$ to take each $S^2_i$ to itself and hence each $P^0_i$ to itself.  This isotopy might move the basepoint $x$ around a loop in $M$, so the new $f$ might induce an inner automorphism of $\pi_1(M,x)$ rather than the identity. However, this inner automorphism respects the free product decomposition of $\pi_1(M,x)$ given by the $P_i$ summands since the new $f$ takes each $P^0_i$ to itself, so this inner automorphism must be the identity.  The result then follows in the case that $\del M$ contains no spheres from the irreducible case and from the fact that the mapping class group of a punctured sphere fixing the boundary spheres is generated by twists along these spheres. 

When $\del M$ contains spheres we can fill them in with balls to obtain a manifold $\overline M$ with $\pi_1\overline M =\pi_1M$.  There is then a fibration
$$
\Dif(M,x) \to \Dif(\overline M,x)\to B
$$
whose base $B$ is the space of configurations of disjoint balls in $\overline M - x$, the number of balls being equal to the number of boundary spheres of $M$. From the last few terms of the long exact sequence of homotopy groups for this fibration we see that the kernel of $\G(M,x)\to\G(\overline M,x)$ is generated by diffeomorphisms of types (3) and (4).
\end{proof}

\begin{lem}
If $\del M$ contains no spheres then an orientation-preserving diffeomorphism $(M,x)\to (M,x)$ that induces the identity on $\pi_1$ also induces the identity on $\pi_2$.
\end{lem}

This is a result of Laudenbach \cite{L}, Appendix III, top of p.142, although the hypothesis on $\del M$ seems to have been omitted there. We will give a more geometric proof than the one there, which uses a spectral sequence argument.

\begin{proof}
Let $\Sigma$ be the collection of spheres in $M$ consisting of the spheres $S^2_i$ that split off the submanifolds $P^0_i$ with $P_i\neq S^2\x S^2$ together with a nonseparating sphere in each $P^0_i$ with $P_i=S^1\x S^2$.  Splitting $M$ along $\Sigma$ then produces the $P^0_i$'s corresponding to irreducible $P_i$'s, together with $S^3_m$, a $3$-sphere with the interiors of $m$ disjoint balls removed. We can take the basepoint $x$ to lie in the interior of this $S^3_m$.  We may assume $m \ge 2$, otherwise $M$ is irreducible and $\pi_2(M,x)=0$. 

Let $\tild\Sigma$ be the pre-image of $\Sigma$ in the universal cover $\tild M$ of $M$.  Splitting $\tild M$ along $\tild\Sigma$ produces copies of $S^3_m$ and copies of the universal covers $\tild P^0_i$ of the $P^0_i$'s corresponding to irreducible $P_i$'s. Dual to $\tild\Sigma$ is a tree $T$ with a vertex for each component of $\tild M -\tild\Sigma$ and an edge for each component of $\tild\Sigma$. We can view $T$ as a quotient space of $\tild M$. The hypothesis that $M$ has no boundary spheres means that no $P_i$ is $D^3$, so every vertex of $T$ has valence at least $2$ since $m\ge 2$. The spheres of $\tild\Sigma$ generate $H_2(\tild M)$ since attaching balls to the boundary spheres of the $\tild P^0_i$'s produces the universal covers $\tild P_i$, which are either contractible or homotopy equivalent to $S^3$ depending on whether $\pi_1(P_i)$ is infinite or finite. 

Let $f\cln (M,x)\to(M,x)$ be an orientation-preserving diffeomorphism that induces the identity on $\pi_1(M,x)$. Then $f$ has a lift $\tilde f\cln \tild M\to \tild M$ which fixes each lift of $x$. Showing that $f$ induces the identity on $\pi_2(M,x)$ is equivalent to showing that $\tilde f $ induces the identity on $\pi_2(\tild M)$. We will do this by using intersection numbers with properly embedded arcs in $\tild M$. The arcs we consider are arcs whose images under the projection $\tild M\to T$ are proper paths joining distinct ends of $T$. We call such arcs in $\tild M$ admissible. Given an admissible arc $\alpha$ and a map $g\cln S^2\to\tild M$ we can perturb $g$ to be transverse to $\alpha$ and then the algebraic intersection number $I_\alpha(g)\in\Z$ is defined once we orient $\alpha$ and choose a fixed orientation for $\tild M$. This intersection number is an invariant of the homotopy class of $g$. It is defined more generally for maps of not-necessarily-connected closed oriented surfaces into $\tild M$, and is invariant under oriented cobordism of such maps.  The intersection number $I_\alpha$ defines a homomorphism $\pi_2(\tild M)=H_2(\tild M)\to \Z$ which depends only on the proper homotopy class of $\alpha$ since every element of $H_2(\tild M)$ is represented by a linear combination of spheres in $\tild\Sigma$, and the intersection number of $\alpha$ with these spheres is invariant under proper homotopy of $\alpha$.

Letting $\alpha$ vary over proper homotopy classes of admissible arcs, we obtain a homomorphism $H_2(\tild M)\to \Pi_\alpha \Z$ which we claim is injective. To prove this it suffices to show that for each nontrivial class in $H_2(\tild M)$ there is an admissible arc that has nonzero intersection number with it. The homology class can be represented by a cycle $c$ which is a linear combination of spheres in $\tild\Sigma$, once orientations are chosen for these spheres. Orienting the spheres of $\tild\Sigma$ is equivalent to orienting the edges of $T$, and then $c$ can be regarded as a $1$-dimensional simplicial cocycle on $T$ with compact support. There is a minimal finite subtree $T_c\subset T$ containing the support of $c$. Let $e$ be an extremal edge of $T_c$. The cocycle $c$ is then nonzero on $e$. One vertex of $e$ abuts edges of $T-T_c$, and we wish to arrange that this is true for the other vertex of $e$ as well. This is automatic if there are infinitely many edges at this second vertex. If there are only finitely many edges at the second vertex and they are all contained in $T_c$ we can rechoose the cycle $c$ within its homology class by adding a suitable multiple of the coboundary of the second vertex of $e$ so that the new $T_c$ is contained in the old one but does not contain $e$. Iterating this process if necessary, we reach the desired situation that both vertices of the edge $e$ abut edges of $T-T_c$. There is then a bi-infinite edgepath in $T$ that intersects $T_c$ only in the edge $e$. This lifts to an admissible arc $\alpha$ having nonzero intersection number with the given homology class.

Returning to the diffeomorphism $\tilde f\cln \tild M \to \tild M$ that fixes the pre-images of the basepoint $x$, let $\alpha$ be an admissible arc in $\tild M$. We can deform $\alpha$ by a proper homotopy so that in each component of $\tild M -\tild\Sigma$ containing a lift of $x$, $\alpha$ passes through that lift. The image $\tilde f(\alpha)$ is then an admissible arc that also passes through these lifts. Hence $\tilde f(\alpha)$ can be deformed to $\alpha$ by a proper homotopy since after projecting $\alpha$ and $\tilde f(\alpha)$ into $T$ there is a proper homotopy and this lifts to a proper homotopy in $\tild M$. Thus for each $g\cln S^2\to\tild M$ we have $I_\alpha(\tilde f g) =I_{\tilde f(\alpha)}(\tilde f g) = I_\alpha(g)$, the second equality coming from the fact that $\tilde f$ is an orientation-preserving diffeomorphism. Since elements of $H_2(\tild M)$ are determined by their intersection numbers with admissible arcs, we conclude that $\tilde f g$ and $g$ determine the same element of $H_2(\tild M)$, hence the same element of $\pi_2(\tild M)$.
\end{proof}

\end{document}